\documentclass[11pt,reqno]{amsart}
\usepackage{xifthen,setspace}
\usepackage{bbm}
\usepackage{esvect}
\usepackage{amsmath} 
\usepackage{lipsum}  
\usepackage{style}

\newcommand{\wbt}{{\widetilde{\beta}_{p,q}}}
\newcommand{\wht}{{\widetilde{h}_{p,q}}}
\newcommand{\bxd}{{\bar{X}_{\cdot 1}}}

\newcommand{\hf}{{H_{\beta,h}}}
\newcommand{\f}{f_{\beta,h}}

\newcommand{\mxs}{{\cM_{\beta,h}}}

\newcommand{\p}{{\P}}
\newcommand{\e}{{\E}}

\newcommand{\Xrb}{{\bar{X}_{\cdot r}}}
\newcommand{\epm}{{\bar{\bm X}_{N}}}

\allowdisplaybreaks

\title[Tensor Curie-Weiss Potts Model]{Limit Theorems and Phase Transitions in the Tensor Curie-Weiss Potts Model}

\author[Bhowal]{Sanchayan Bhowal}
\address{Statistics and Mathematics Unit, Indian Statistical Institute, Bangalore, India, {\tt sanchayan.bhowal2509@gmail.com}}

\author[Mukherjee]{Somabha Mukherjee} 
\address{Department of Statistics and Data Science, National University of Singapore, Singapore. {\tt somabha@nus.edu.sg}}

\begin{document}

\begin{abstract}
    In this paper, we derive results about the limiting distribution of the empirical magnetization vector and the maximum likelihood (ML) estimates of the natural parameters in the tensor Curie-Weiss Potts model. Our results reveal surprisingly new phase transition phenomena including the existence of a smooth curve in the interior of the parameter plane on which the magnetization vector and the ML estimates have mixture limiting distributions, the latter comprising of both continuous and discrete components, and a surprising superefficiency phenomenon of the ML estimates, which stipulates an $N^{-3/4}$ rate of convergence of the estimates to some non-Gaussian distribution at certain \textit{special} points of one type and an $N^{-5/6}$ rate of convergence to some other non-Gaussian distribution at another special point of a different type. The last case can arise only for one particular value of the tuple of the tensor interaction order and the number of colors. These results are then used to derive asymptotic confidence intervals for the natural parameters at all points where consistent estimation is possible. 
\end{abstract}



\maketitle

\section{Introduction}
The Potts model \cite{fywu}, originally named after Renfrey Potts \cite{pottsorig}, is a generalization of the Ising model \cite{isingorig}, where the \textit{spin} of any particular site can have more than two states, each such state being referred to as a color. It finds broad application in elucidating diverse physical phenomena, including magnetism, phase transitions, and social behavior. This model is related to a number of other well-known models, such as the Heisenberg model, the XY model, and the Ashkin-Teller model (the four-state Potts model), and has found extensive applications in a number of diverse fields including biomedical problems \cite{cellular, gene77}, image processing and computer vision \cite{impro, impro2}, spatial statistics \cite{spatstat}, social sciences \cite{socialsci} and finance \cite{fina66,bornholdt}. The classical Potts model represents pairwise (quadratic) interactions between the sites, which, most often, is not enough to capture the complex dependencies present in real world network data. For example, in a peer group, the behavior of an individual does not depend only on pairwise interactions between his/her friends, but is a function of more complex higher order interactions. In a different context, it is known in chemistry that the atoms on a crystal surface do not interact just in pairs, but in triangles, quadruplets and higher order tuples. A natural extension of the classical Potts model that captures multibody interactions, is the tensor Potts model, and in this paper, we consider the problem of deriving the asymptotics of a natural estimate of the parameters of this model, given only one sample from the model. Obtaining precise asymptotics of the sufficient statistic and the parameter estimates in general tensor Potts models is notoriously difficult, unless one agrees to assume certain special structures on the underlying network. One such natural structural condition is to assume that all tuples of nodes of a fixed order (say, $p$) interact with each other, with a uniform interaction strength. The resulting model is the tensor Potts model on the $p$-uniform complete hypergraph, also referred to as the $p$-tensor Curie-Weiss Potts model. 

A close relative of the Potts model is the Ising model \cite{isingorig}, where there is a huge literature on the problem of consistent parameter estimation. Chatterjee \cite{chatterjee} showed how to estimate the parameters of a general spinglass model consistently, using the idea of \textit{pseudolikelihood estimation}, which was introduced by Besag \cite{besag_lattice, besag_Nl} in the context of spatial statistics. A myriad of works followed in the next few years on the problem of partial and joint estimation of Ising model parameters, some notable ones among them being \cite{BM16, pg_sm, cd_ising_I, cd_ising_II, cd_testing}. In a rather different context, one might be interested in estimating the entire structure (interaction matrix) of a general Ising model, assuming that she has access to multiple samples from such a model. This problem is known as \textit{structure learning}, and has been addressed in details in a series of works \cite{highdim_ising, structure_learning, bresler, lokhov}. The problem of deriving exact asymptotics of the magnetization and parameter estimates in the Curie-Weiss Ising model was addressed in \cite{comets,elcw2}, and in \cite{comets_exp} for Markov random fields on lattices. However, the Ising models in all these works capture only pairwise interactions, which as we discussed above, is often not a practical assumption in many realistic settings involving peer-group effects and multi-particle interactions. A natural substitute for the classical $2$-spin Ising model in such situations, is the $p$-spin Ising model \cite{ab_ferromagnetic_pspin,bovier}. Consistent estimation of the natural parameters in general $p$-spin Ising models was established in \cite{smmpl}, and exact fluctuations of the magnetization and parameter estimates were established for the $p$-spin Curie-Weiss model in \cite{smfl,mlepaper}. However, to the best of our knowledge, nothing is known about the asymptotics of the empirical magnetization vector and the parameter estimates for the closely related $p$-spin Potts model, even for the fully connected case, although the corresponding asymptotics have been established in the $2$-spin case in \cite{ellis, gan, eich}. This is precisely the goal of this paper. We will see that even in this simple case where we have a $p$-spin Curie-Weiss Potts model, many surprising phase transitions arise in the asymptotics of the magnetization vector and the parameter estimates. Some salient features of these surprising phenomena include the appearance of rates of convergence (of the estimates) like $N^{-3/4}$ and $N^{-5/6}$ at some \textit{special} points in the parameter space, and the existence of a smooth curve in the interior of the parameter space, where the estimates have limiting mixture distributions.

\subsection{Model Description} For integers $p\ge 2$ and $q\ge 2$, the $p$-tensor Potts model is a discrete probability distribution on the set $[q]^N$ (here and afterwards, for a positive integer $m$, we will use $[m]$ to denote the set $\{1,2,\ldots,m\}$)  for some positive integers $q$ and $N$, given by:
\begin{equation}\label{gp}
    \p_{\beta,h,N}(\bm X) := \frac{1}{q^N Z_N(\beta,h)} \exp \left(\beta \sum_{1\le i_1,\ldots,i_p\le N} J_{i_1,\ldots,i_p}\mathbbm{1}_{X_{i_1}=\ldots=X_{i_p}} + h\sum_{i=1}^N\mathbbm{1}_{X_i=1}\right) \quad(\bm X \in [q]^N)~,
\end{equation}
where $\beta>0$, $h \geq 0$ and $\bm J := ((J_{i_1, \ldots,i_p}))_{i_1,\ldots,i_p\in [N]}$ is a symmetric tensor. The $p$-tensor Curie-Weiss Potts model is obtained by taking $J_{i_1,\ldots,i_p} := N^{1-p}$ for all $(i_1,\ldots,i_p) \in [N]^p$, whence model \eqref{gp} takes the form:
\begin{equation}\label{cp}
    \p_{\beta,h,N}(\bm X) := \frac{1}{q^N Z_N(\beta,h)} \exp \left(\beta N\sum_{r=1}^q \Xrb^p + Nh\bar{X}_{\cdot 1}\right) \quad(\bm X \in [q]^N)
\end{equation}
where $\Xrb := N^{-1} \sum_{i=1}^N X_{i,r}$ with $X_{i,r}:= \mathbbm{1}_{X_i=r}$. The variables $p$ and $q$ are called the \textit{interaction order} and the \textit{number of states/colors} of the Potts model. A sufficient statistic for the exponential family \eqref{cp} is the empirical magnetization vector:
$$\epm := \left(\bar{X}_{\cdot 1},\ldots,\bar{X}_{\cdot q}\right)^\top~. $$
Note that $\epm$ is a probability vector, i. e. has non-negative entries adding to $1$. In this paper, we give a complete description of the asymptotics of $\epm$ on the entire parameter space:
$$\Theta := \{(\beta,h): \beta>0, h \geq 0\} = (0,\infty)\times [0,\infty)~. $$
We then use these asymptotics to establish limit theorems for the maximum likelihood (ML) estimators of $\beta$ and $h$, which is crucial for constructing asymptotic confidence intervals for these parameters.

\subsection{Maximum Likelihood Estimation} Hereafter, given $\bm X \sim \mathbb{P}_{\beta, h, p}$, we denote by $\hat{\beta}_N$ and $\hat{h}_N$ the marginal maximum likelihood (ML) estimators of $\beta$ and $h$, respectively. It follows from Lemma \ref{ml_exp}, that for fixed $h \in \mathbb{R}, \hat{\beta}_N$ is a solution of the equation (in $\beta$),
\begin{equation}
    \mathbb{E}_{\beta, h, p}\left(\|\epm\|_p^p\right)=\|\epm\|_p^p~,
    \label{b_like}
\end{equation}
and for fixed $\beta\in \mathbb{R}$, $\hat{h}_N$ is a solution of the equation (in $h$),
\begin{equation}
   \mathbb{E}_{\beta, h, p}\left(\bar{X}_{\cdot 1}\right)=\bar{X}_{\cdot 1}.
    \label{h_like}
\end{equation}
The limiting distribution of the ML estimates of $h$ and $\beta$ therefore depend on the fluctuations of the average magnetization $\epm$ across the parameter space $\Theta$. The main features of these asymptotics are highlighted below:

\begin{itemize}
    \item The parameter space $\Theta$ has a subset of \textit{regular} points, where the magnetization vector and the ML estimates are asymptotically normal, their rates of convergence being $N^{-1/2}$.

    \item The complement of the set of regular points contains the so called \textit{critical points}, which forms a continuous curve in the interior of the parameter space, on which the magnetization vector and the ML estimates have limiting mixture distributions, the latter consisting of both continuous and discrete components. 

    \item The remaining portion of the parameter space consists of exactly one \textit{special} point, where the magnetization and the ML estimates have rates of convergence different from the classical $N^{-1/2}$ rate. In case $(p,q)\ne (4,2)$, the magnetization converges at rate $N^{-1/4}$ and the parameter estimates at rate $N^{-3/4}$ to limiting non-Gaussian distributions. On the other hand, if $(p,q)=(4,2)$, the convergence rate of the magnetization at the special point changes to $N^{-1/6}$, whereas the estimates converge at rate $N^{-5/6}$. The estimates are thus \textit{superefficient} at the special points.
\end{itemize}
Note that the $N^{-5/6}$ convergence rate for the ML estimates is a special phenomenon noticed in the $4$-spin, $2$-color Curie-Weiss Potts model, that is never observed in the closely related tensor Curie-Weiss Ising models, or in the classical $2$-spin Curie-Weiss Potts models. In Figures \ref{fig1} and \ref{fig2}, we illustrate the different phase transitions through phase diagrams. 

The rest of the paper is organized as follows. In Section \ref{asymmagn} we describe the asymptotics of the magnetization vector of the $p$-spin Curie-Weiss Potts model. These asymptotics depend on the location of the parameters on one of the several components of a partition induced by the so called \textit{free energy function}, mainly characterized by whether this function has one or multiple global maximizers, and what is the order of the first non-zero derivative at these maximizers. We use the results in Section \ref{asymmagn} to derive limiting distributions of the ML estimators in Section \ref{sec:asmml}. In Section \ref{sec:confint}, we use the results in Section \ref{sec:asmml} to derive asymptotic confidence intervals for the model parameters. In that section, we also summarize the partition of the parameter space into the regular, critical and special points as sketched above, in details. A brief sketch of the proofs of the main results in this paper is given in Section \ref{prskt}. Finally, complete proofs of all the results in the main paper are given in the appendix.

\section{Asymptotics of the Magnetization Vector}\label{asymmagn}
In this section, we state our main results regarding the asymptotics of the magnetization vector. For this, we need a few definitions and notations. For $p,q\ge 2$ and $(\beta,h) \in \Theta$, the \textit{negative free energy} function $H_{\beta,h}: \cP_q \to \R$ is defined as:
$$H_{\beta,h}(\bm t) := \beta \sum_{r=1}^q t_r^p + ht_1 - \sum_{r=1}^q t_r \log t_r$$  
where $\cP_q$ denotes the set of all $q$-dimensional probability vectors. We start by showing that the magnetization vector concentrates around the set $\mathcal{M}_{\beta,h}$ of all global maximizers of the function $H_{\beta,h}$. Actually, this and all the subsequent results in this section are proved under slightly perturbed versions of the model parameters.

\begin{thm}\label{conclem}
    Let $\beta_N \rightarrow \beta$ and $h_N \rightarrow h$. Then, under $\p_{\beta_N,h_N,N}$, the empirical magnetization $\epm$ satisfies a large deviation principle with speed $N$ and rate function $-H_{\beta,h} + \sup H_{\beta,h}$.
    Consequently, for a point $\bm t \in \R^q$ and a set $A\subseteq \R^q$, if we define $d(\bm t, A) := \inf_{\bm a\in A} \|\bm t - \bm a\|_2$, then for every $\varepsilon>0$, there exists a constant $C_{q,\varepsilon} >0$ depending only on $q$ and $\varepsilon$, such that:
    $$\p_{\beta_N,h_N,N}\left(d(\epm,\mxs) \ge \varepsilon\right) \le e^{-C_{q,\varepsilon} N}$$ for all large $N$.
\end{thm}

Theorem \ref{conclem} is proved in Appendix \ref{prclem}. It enables us to derive a law of large numbers of the magnetization vector towards the set $\cM_{\beta,h}$ of global maximizers of $H_{\beta,h}$. We now derive the fluctuations of the magentization vector around $\cM_{\beta,h}$, which depends, among other things, on the location of the point $(\beta,h)$ in the parameter space.
\begin{defn} We partition the parameter space into the following three components:
    \begin{enumerate}
        \item A point $(\beta,h) \in \Theta$ is called \textit{regular}, if the function $H_{\beta,h}$ has a unique global maximizer $\bm m_*$ and the quadratic form
              $$\bm Q_{\bm s,\beta}(\bm t) := \sum_{r=1}^q \left(\beta p(p-1)s_r^{p-2} - \frac{1}{s_r}\right) t_r^2~,$$
              is negative definite on $\cH_q := \{\bm t\in \R^q: \sum_{r=1}^q t_r=0\}$ for $\bm s = \bm m_*$. The set of all regular points is denoted by $\cR_{p,q}. $
        \item A point $(\beta,h) \in \Theta$ is called \textit{critical}, if $H_{\beta,h}$ has more than one global maximizer, and for each such global maximizer $\bm m$, the quadratic form $\bm Q_{\bm m,\beta}$ is negative definite on $\cH_q$. The set of all critical points is denoted by $\cC_{p,q}$.
        \item A point $(\beta,h) \in \Theta$ is called \textit{special}, if $H_{\beta,h}$ has a unique global maximizer $\bm m_*$ and the quadratic form $\bm Q_{\bm m_*,\beta}$ is singular on $\cH_q$ (i.e. $\mathrm{Ker}(\bm Q_{\bm m_*,\beta}) \bigcap \cH_q \ne \{\boldsymbol{0}\}$). The set of all special points is denoted by $\cS_{p,q}$.
    \end{enumerate}
\end{defn}

It is proved in Lemma \ref{partn} in the appendix, that the above three subsets indeed form a partition of the parameter space $\Theta$. From Proposition \ref{oned_reduce}, it follows that the global maximizers of $H_{\beta,h}$ can be reparametrized as permutations of the vector
$$
    \bm x_s=\left(\frac{1+(q-1)s}{q},\frac{1-s}{q},\ldots,\frac{1-s}{q}\right).
$$
for some $s \in [0,1)$, and hence, the problem can be reduced to a one dimenional optimization of the function $f_{\beta,h}(s) : =H_{\beta,h}(\bm x_s)$. Note that the map $s\mapsto \bm x_s$ is one-one, since $s = 1-q x_{s,2}$.

We write $f_{\beta,h}(s)$ as,
$$f_{\beta,h}(s)=(q-1)k\left(\frac{1-s}{q}\right)+k\left(\frac{1+(q-1)s}{q}\right) +\left(\frac{1+(q-1)s}{q}\right)\cdot h,$$
where $k(x)= k_{\beta,p}(x) := \beta x^p-x\log x$.
Hence for $\bm t \in \cH_q$,
\begin{equation}
    \begin{aligned}
        \bm Q_{\bm x_s,\beta}(\bm t) & = k^\dprime\left(\frac{1-s}{q}\right)\sum_{r=2}^{q}t_r^2 +k^\dprime\left(\frac{1+(q-1)s}{q}\right)\left(\sum_{r=2}^{q}t_r\right)^2 \\
    \end{aligned}
    \label{q_form}
\end{equation}
\begin{defn}
    We now further classify the special points into the following two categories:
    \begin{enumerate}[i.]
        \item A special point $(\beta,h)$ is said to be of \textit{type-I}, if the unique global maximizer $\bm m_* =: \bm x_s$ satisfies $f_{\beta,h}^{(4)}(s)<0$. The set of all type-I special points is denoted by $\cS^1_{p,q}$.
        \item A special point $(\beta,h)$ is said to be of \textit{type-II}, if the unique global maximizer $\bm m_* =: \bm x_s$ satisfies $f_{\beta,h}^{(4)}(s)=0$. We denote the set of all type-II special points by $\cS^2_{p,q}$.
    \end{enumerate}
\end{defn}

We now state our results regarding the central limit theorem (CLT) of the magnetization under the $p$-tensor Potts model with perturbed parameters. We begin with the CLT at regular points.
\begin{thm}\label{cltreg}
    Suppose $(\beta, h)$ is regular and let  $\bm m_*= \bm m_*(\beta, h) = \bm x_s$ denote the unique maximizer of $H_{\beta,h}$. Then, for $\bm X \sim \P_{\beta+N^{-\frac{1}{2}} \bar{\beta}, h+N^{-\frac{1}{2}} \bar{h}, p}$ for some $\bar{\beta},\bar{h} \in \mathbb{R}$, as $N \rightarrow \infty$, we have:
    \begin{equation*}
        N^{\frac{1}{2}}\left(\epm - \bm m_*\right) \xrightarrow{D} \cN_q\left(\Sigma (\bar{\beta}p\bm m_*^{p-1}+\bar{h}\bm e_1), \Sigma\right),
    \end{equation*}
    where, $\bm x^\ell:=(x_1^\ell,\ldots,x_d^\ell)$ for $\bm x \in \R^d$, $\bm e_1=(1,0,\ldots,0)$, and
    \begin{equation}
        \label{sigma_def}
        \Sigma:=
        \left(-\frac{q^2}{q-1}f^\dprime_{\beta,h}(s)\right)^{-1}\left(\begin{array}{cccc}q-1 & -1 & \cdots & -1 \\ -1 & 1+(q-2) \frac{k^\dprime\left(\frac{1+(q-1)s}{q}\right)}{k^\dprime\left(\frac{1-s}{q}\right)} & & -\frac{k^\dprime\left(\frac{1+(q-1)s}{q}\right)}{k^\dprime\left(\frac{1-s}{q}\right)} \\ \vdots & & \ddots & \\ -1 & -\frac{k^\dprime\left(\frac{1+(q-1)s}{q}\right)}{k^\dprime\left(\frac{1-s}{q}\right)} & & 1+(q-2) \frac{k^\dprime\left(\frac{1+(q-1)s}{q}\right)}{k^\dprime\left(\frac{1-s}{q}\right)}\end{array}\right).
    \end{equation}
\end{thm}

The proof of Theorem \ref{cltreg} is given in Appendix \ref{prcltreg}. Next, we state the CLT result at the critical points.

  \begin{figure}[ht]
\centering
\centering
\includegraphics[width=4in]
    {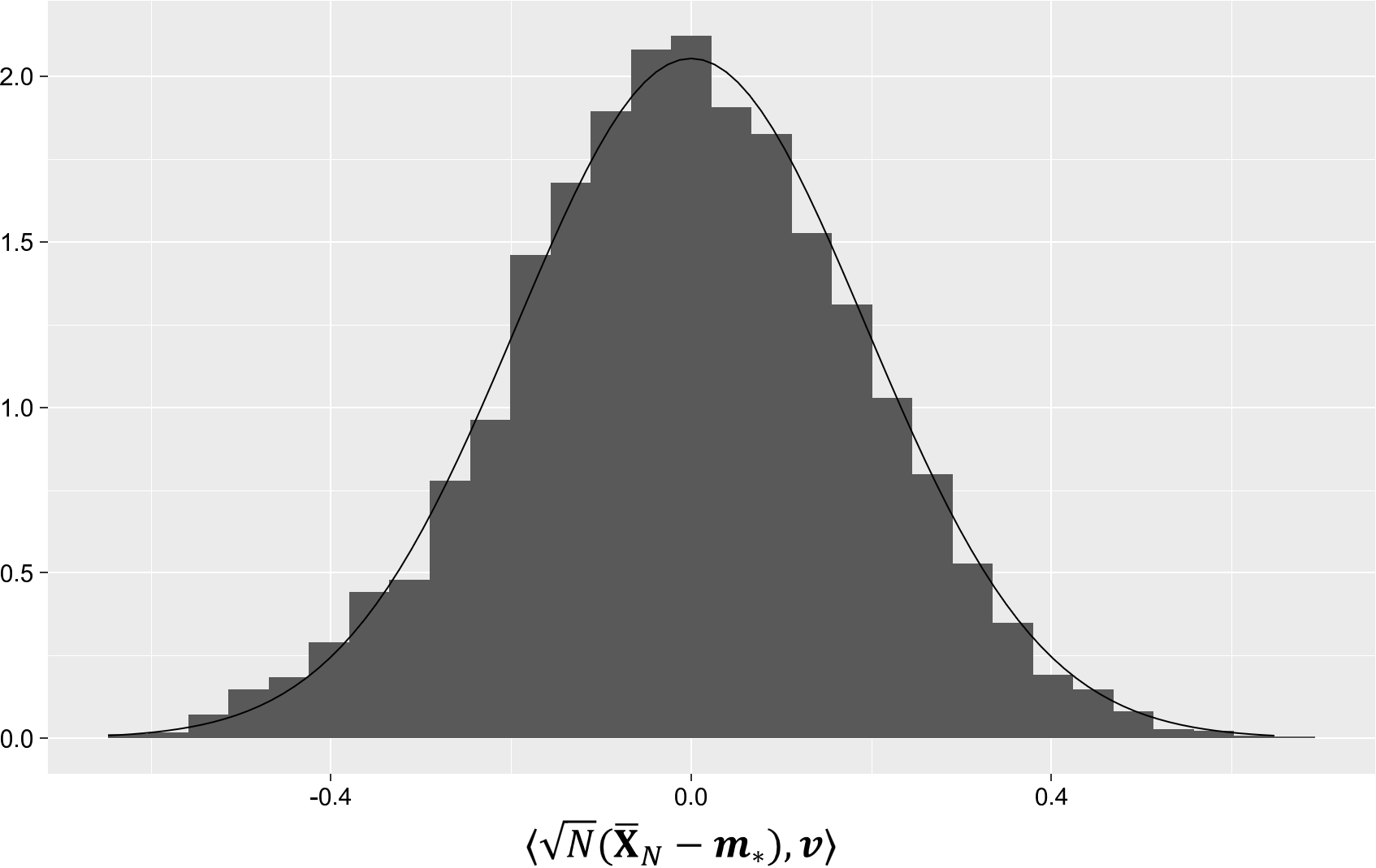}\\
\caption{\small{Histogram and theoretical density curve of $\sqrt{N}(\epm-\bm m_*)$ projected at a random direction $\bm v := (0.157, 0.396, 0.323)$ at a regular point ($\beta=0.616$, $h=0.67$).}}
\label{fig23071}
\end{figure}

\begin{thm}\label{pcr}
    Suppose $(\beta, h)$ is critical and let $\bm m_1:=$ $\bm m_1(\beta, h, p)$, $\ldots$, $\bm m_K:= \bm m_K(\beta, h, p)$ denote the $K$ maximizers of $H_{\beta,h}$. Then, for $\bm X \sim \P_{\beta, h, p}$, as $N \rightarrow \infty$, we have:
    \begin{equation}
        \label{crit_wlln}
        \epm \xrightarrow{P} \sum_{k=1}^K p_k \delta_{m_k},
    \end{equation}
    where
    \begin{equation*}
        p_k:=\frac{\tau(\bm m_k)}{\sum_{i=1}^{K}\tau(\bm m_i)},
    \end{equation*}
    and 
    \begin{equation}
            \tau(\bm m_i):=\sqrt{\frac{f^\dprime_{\beta,h}(s_i)^{-1}\left(-k^\dprime\left(\frac{1-s_i}{q}\right)\right)^{2-q}}{\prod_{r=1}^q m_{i,r}}}~.
            \label{tau_def}
    \end{equation}
    where $\bm m_i$ is a permutation of $\bm x_{s_i}$. Moreover, if $(\beta_N,h_N) = (\beta+\frac{\bar{\beta}}{\sqrt{N}},h+\frac{\bar{h}}{\sqrt{N}}$) for some critical point $(\beta,h)$, then for every $\varepsilon >0$ smaller than the minimum distance between any two global maximizers of $H_{\beta,h}$, we have the following under $\p_{\beta_N,h_N}\left(\cdot \Big| \epm \in B(\bm m_i ,\varepsilon) \right)$:
     
              \begin{equation*}
                  \sqrt{N}\left(\epm - \bm m_i\right) \xrightarrow{D} \cN_q(P\Sigma P^\top  (\bar{\beta}p\bm m_i^{p-1}+\bar{h}\bm e_1), P\Sigma P^\top),
              \end{equation*}
    where $\Sigma$ is as defined in \eqref{sigma_def} and $P$ is the permutation matrix corresponding to the permutation, i.e. $\bm m_i := P \bm x_{s_i}$.
\end{thm}

 Theorem \ref{pcr} is proved in Appendix \ref{proofprpcr}. Finally, we state the CLT result at the special points. We start with the CLT for type-I special points.  

    \begin{figure}[ht]
\centering
\centering
\includegraphics[width=5.5in]
    {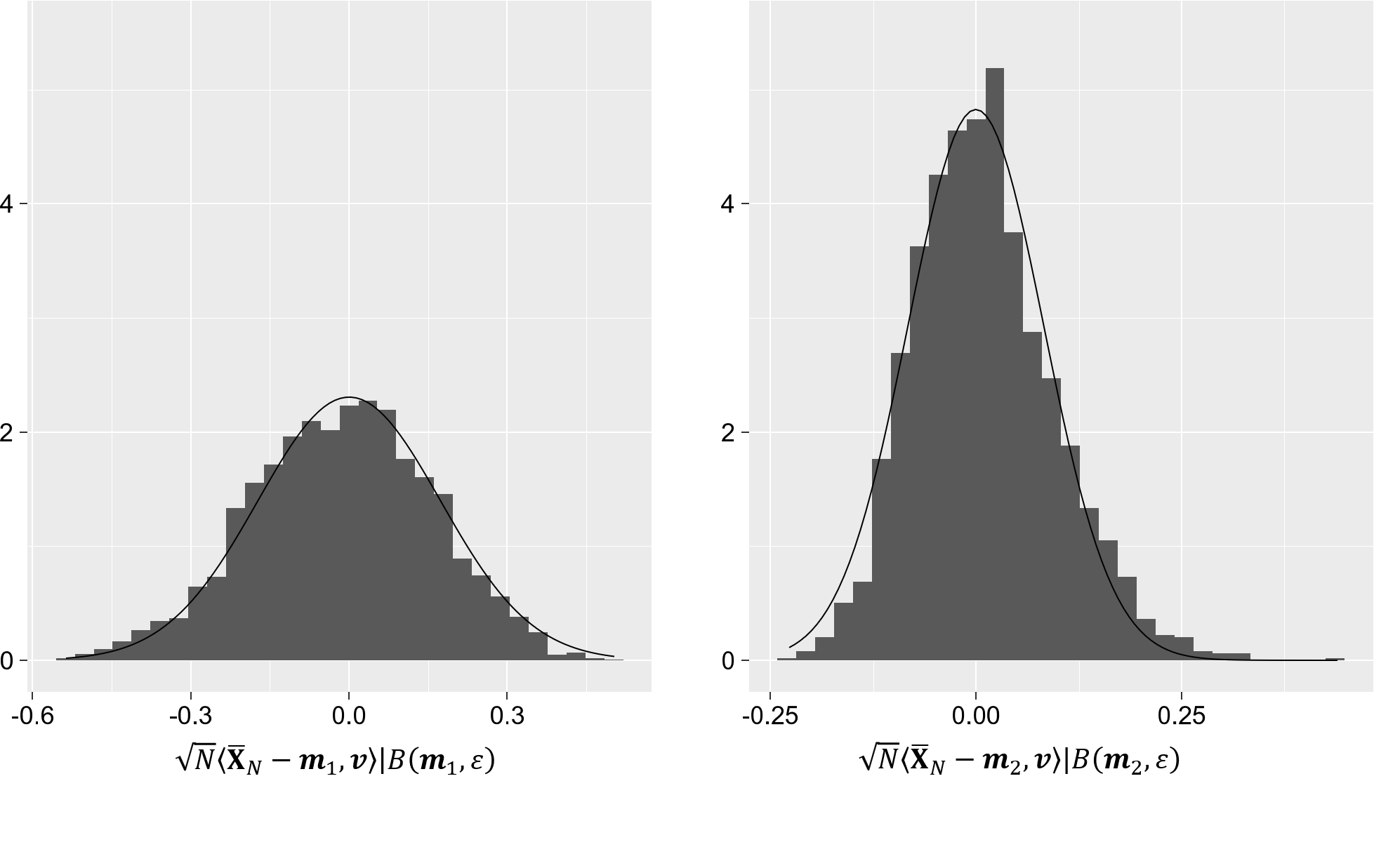}\\
\caption{\small{Conditional histograms and theoretical density curves of $\sqrt{N}(\epm - \bm m_i)$ projected at a random direction $\bm v := (0.157, 0.396, 0.323)$ at a strongly critical point ($\beta=0.965$, $h=0.2$).}}
\label{fig23072}
\end{figure}

\begin{thm}\label{psp1}
    Suppose $(\beta, h)$ is type-I special, and let $m_*=m_*(\beta, h, p) = \bm x_s$ denote the unique global maximizer of $H_{\beta,h}$. Define $\bm u := (1-q,1,\ldots,1)$. Note that there exists unique $T_N$ and $\bm V_N \in \cH_q \cap\;\operatorname{Span}(\bm u)^\perp$ such that $\epm - \bm m_*=N^{\frac{-1}{4}}T_N \bm u+ N^{\frac{-1}{2}}\bm V_N$. Then, for $\bm X \sim \P_{\beta+N^{-\frac{3}{4}} \bar{\beta}, h+N^{-\frac{3}{4}} \bar{h}, p}$, as $N \rightarrow \infty$, we have:
    \begin{equation*}
        T_N \xrightarrow{D} T := T_{\bar{\beta},\bar{h}}
    \end{equation*}
    where $T_{\bar{\beta},\bar{h}}$ is a random variable with density at $x$ proportional to,
    \begin{equation}
        \label{spec1_dist}
        \exp\left(\frac{x^4}{24}q^4f^{(4)}_{\beta,h}(s)+  (\bar{\beta}p\ip{\bm m_*^{p-1}}{\bm u}+\bar{h}(1-q))x\right).
    \end{equation}
    Also, $$\bm V_N \xrightarrow{D} \bm V$$ where $\bm V$ is a multivariate normal random vector in $\R^{q}$ with mean $\bm 0$ and covariance matrix of rank $q-2$, given by:
    $$
        \frac{1}{-(q-1)k^\dprime\left(\frac{1-s}{q}\right)}\begin{bmatrix}
            0      & 0   & \ldots & 0   \\
            0      & q-2 & \ldots & -1  \\
            \vdots &     & \ddots       \\
            0      & -1  & \ldots & q-2
        \end{bmatrix}
    $$
    Further, $T$ and $\bm V$ are independent.
\end{thm}

Theorem \ref{psp1} is proved in Appendix \ref{proofpsp1}. To conclude, we prove the CLT for type-II special points. 

\begin{thm}\label{psp2}
    Suppose $(\beta, h)$ is type-II special, and let $m_*=m_*(\beta, h, p) = \bm x_s$ denote the unique maximizer of $H_{\beta,h}$. Define $\bm u :=(1-q,1,\ldots,1)$. Then, for $\bm X \sim \P_{\beta+N^{-\frac{5}{6}} \bar{\beta}, h+N^{-\frac{5}{6}} \bar{h}, p}$, as $N \rightarrow \infty$,
    \begin{equation*}
        N^{\frac{1}{6}}\left(\epm - \bm m_*\right) \xrightarrow{D} F_{\bar{h}}\bm u.
    \end{equation*}
    where the random variable $F_{\bar{h}}$ has density with respect to the Lebesgue measure is proportional to
    \begin{equation}
        \label{spec2_dist}
        \exp\left(-\frac{32}{15}x^6 - \bar{h}x\right).
    \end{equation}
\end{thm}

\begin{figure}[ht]
\centering
\centering
\includegraphics[width=4in]
    {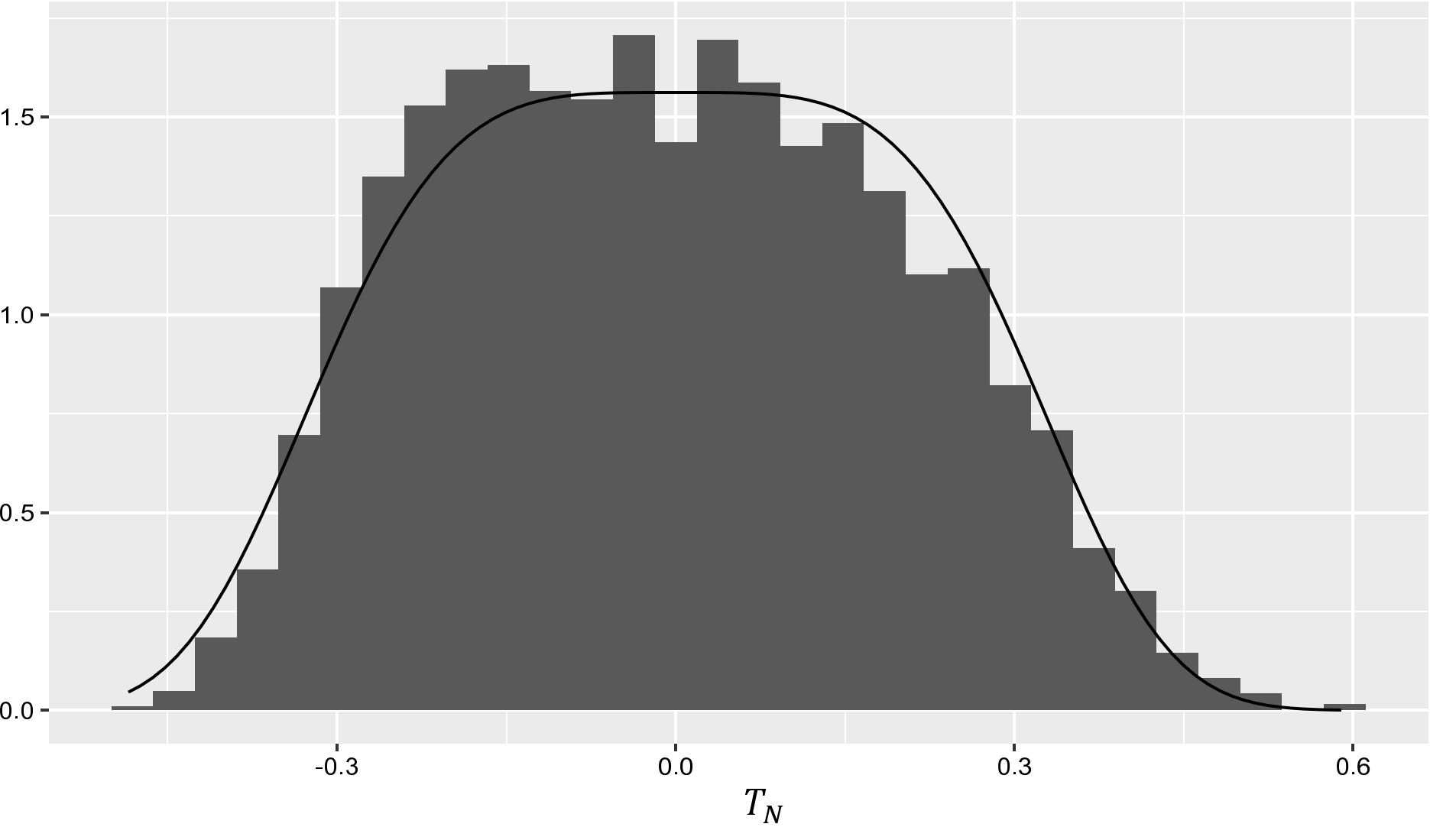}\\
\caption{\small{Histogram and theoretical density of $T_N$ at a (type-I) special point ($\beta=0.778$, $h=0.485$).}}
\label{fig23073}
\end{figure}

Theorem \ref{psp2} is proved in Appendix \ref{prooftype2spc}. In Figures \ref{fig23071}, \ref{fig23072} and \ref{fig23073}, we compare the empirical distrbutions of the magnetization with their corresponding asymptotic theoretical distributions as stated in the above theorems, in each of the three cases where the true parameter is regular, critical and special. The simulations were performed for the case $p=4, q=3$ with $N = 1000$, using MCMC. 

\section{Asymptotics of the Maximum Likelihood Estimates}\label{sec:asmml}
In this section, we prove results about the asymptotics of the maximum likelihood (ML) estimates of the parameters $\beta$ and $h$. We define $u_{N, p}$ and $u_{N, 1}$ to be the functions appearing in the LHS of the equations \eqref{b_like} and \eqref{h_like}, respectively, that is,
$$u_{N,p}(\beta,h,p) := \e_{\hat{\beta}_N,h,p}(\|\epm\|_p^p)\quad\textrm{and}\quad u_{N,1}(\beta,h,p) := \e_{\beta,h,p}(\bar{X}_{\cdot 1})~.$$
It follows from Lemma \ref{ml_exp} that for fixed $h$, the ML estimate $\hat{\beta}$ satisfies the equation:
$$u_{N,p}(\beta,h,p) = \|\epm\|_p^p$$ and for fixed $\beta$, the ML estimate $\hat{h}$ satisfies the equation:
$$u_{N,1}(\beta,h,p) = \bar{X}_{\cdot 1}.$$ We start with the results about the asymptotic distribution of $\hat{h}_N$, which depend on whether the underlying parameters are regular, special or critical. 

\begin{thm}[Asymptotic distribution of $\hat{h}_N$ at regular points]\label{hmlreg}
    Fix $p \geq 2$ and suppose $(\beta, h) \in$ $\Theta$ is regular. Assume $\beta$ is known and $\bm X \sim \P_{\beta, h, p}$. Then, denoting the unique maximizer of $H$ by $\bm m_*=\bm m_*(\beta, h, p)$, as $N \rightarrow \infty$, we have:
    $$
        N^{\frac{1}{2}}\left(\hat{h}_N-h\right) \xrightarrow{D} \cN \left(0,-\frac{q^2}{(q-1)^2}f^\dprime_{\beta,h}(s)\right)
    $$
\end{thm}

Theorem \ref{hmlreg} is proved in Appendix \ref{prhmlreg}. It shows that $\hat{h}_N$ is $N^{\frac{1}{2}}$-consistent and asymptotically normal at the regular points. Before discussing more about the implications of this theorem, we state the result for the asymptotic distribution of $\hat{h}_N$ when $(\beta, h)$ is special.

\begin{thm}[Asymptotic distributions of $\hat{h}_N$ at special points]\label{hmle2}
    Fix $p \geq 2$ and suppose $(\beta, h) \in$ $\Theta$ is special. Assume $\beta$ is known and $\bm X \sim \P_{\beta, h, p}$. Denote the unique maximizer of $H$ by $\bm m_*=\bm m_*(\beta, h, p)$.
    \begin{enumerate}
        \item If $(\beta, h) \in$ $\Theta$ is type I special then, as $N \rightarrow \infty$,
              $$
                  N^{\frac{3}{4}}\left(\hat{h}_N-h\right) \xrightarrow{D} G_1
              $$
              where the distribution function of $G_1$ is given by
              $$
                  G_1(t)=R_{0,0}\left(\int_{-\infty}^{\infty} u \mathrm{~d} R_{0, t}(u)\right),
              $$
              where $R_{\bar{\beta},\bar{h}}$ denotes the distribution function of the random variable $T_{\bar{\beta},\bar{h}}$ as defined in \eqref{spec1_dist}.
              
        \item If $(\beta, h) \in$ $\Theta$ is type II special then, as $N \rightarrow \infty$,
              $$
                  N^{\frac{5}{6}}\left(\hat{h}_N-h\right) \xrightarrow{D} G_2
              $$
              where the distribution function of $G_2$ is given by
              $$
                  G_2(t)=H_{0}\left(\int_{-\infty}^{\infty} u \mathrm{~d} H_{t}(u)\right),
              $$
              where $H_{\bar{h}}$ denotes the distribution function of $F_{\bar{h}}$ as defined in \eqref{spec2_dist}.
    \end{enumerate}
\end{thm}

The proof of Theorem \ref{hmle2} is exactly similar to the proof of Theorem \ref{hmlreg}, so we skip it. It shows that at the type-I and type-II special points, $\hat{h}_N$ is superefficient, and is $N^{3/4}$ and $N^{5/6}$-consistent, respectively, and the limiting distributions are also non-Gaussian. We now state the result on the asymptotics of $\hat{h}_N$ at the critical points. For this, we need a few definitions:
\begin{defn}
    For $\sigma>0$, the positive half-normal distribution $\cN^{+}\left(0, \sigma^2\right)$ is defined as the distribution of $|Z|$, where $Z \sim \cN\left(0, \sigma^2\right)$, and the negative half-normal distribution $\cN^{-}\left(0, \sigma^2\right)$ is defined as the distribution of $-|Z|$, where $Z \sim \cN\left(0, \sigma^2\right)$.
\end{defn}

\begin{defn}
    We partition the set of critical points as follows:
    \begin{enumerate}[i.]
        \item If $(\beta,h)$ is a critical point such that $f_{\beta,h}$ has more than one global maximizer then it is called \textit{strongly critical}. We denote the set of all strongly critical points as $\cC^1_{p,q}$.
        \item If $(\beta,h)$ is a critical point such that $f_{\beta,h}$ has a unique global maximizer then it is called \textit{weakly critical}. We denote the set of all weakly critical points as $\cC^2_{p,q}$.
    \end{enumerate}
\end{defn}

\begin{thm}[Asymptotic distributions of $\hat{h}_N$ at critical points]\label{hml77}
    Suppose that $(\beta, h)$ is a critical point. Let $p_1,\ldots,p_K$ be the weights defined in the statement of Theorem \ref{pcr} for the global maximizers $\bm m_1,\ldots,\bm m_K$, respectively, where these maximizers are arranged in ascending order of their first coordinates. Then, for $\bm X \sim \P_{\beta, h, p}$, as $N \rightarrow \infty$, we have the following:
    \begin{enumerate}
        \item If $(\beta,h)\in\cC^1_{p,q} \setminus \{(\beta_c,0)\}$, then $f_{\beta,h}$ has exactly two global maximizers $s_2>s_1>0$, and
              \begin{equation*}
                  N^{\frac{1}{2}}\left(\hat{h}_N-h\right) \xrightarrow{D} \frac{p_1}{2}\cN^{-} \left(0,-\frac{q^2}{(q-1)^2}f^\dprime_{\beta,h}(s_1)\right)+\frac{1-p_1}{2}\cN^{+} \left(0,-\frac{q^2}{(q-1)^2}f^\dprime_{\beta,h}(s_2)\right)+\frac{1}{2} \delta_0,
              \end{equation*}
        \item If $(\beta,h)\in\cC^2_{p,q}$, then $f_{\beta,h}$ has exactly one global maximizer $s>0$, and
              \begin{eqnarray*}
                  N^{\frac{1}{2}}\left(\hat{h}_N-h\right) &\xrightarrow{D}& \frac{1-p_q}{2}\cN^{-} \left(0,-\frac{q^2 f^\dprime_{\beta,h}(s)}{(q-1)\left(1+(q-2) \frac{k^\dprime\left(\frac{1+(q-1)s}{q}\right)}{k^\dprime\left(\frac{1-s}{q}\right)}\right)}\right)\\ &+&\frac{p_q}{2}\cN^{+} \left(0,-\frac{q^2}{(q-1)^2}f^\dprime_{\beta,h}(s)\right) +\frac{1}{2}\delta_0,
              \end{eqnarray*}
        \item If $(\beta,h)=(\beta_c,0)$, then $f_{\beta,h}$ has exactly two global maximizers, $0$ and $s > 0$, and
              \begin{multline*}
                  N^{\frac{1}{2}}\left(\hat{h}_N-h\right) \xrightarrow{D} \frac{(1-p_q)(q-1)}{2q}\cN^{-} \left(0,-\frac{q^2 f^\dprime_{\beta,h}(s)}{(q-1)\left(1+(q-2) \frac{k^\dprime\left(\frac{1+(q-1)s}{q}\right)}{k^\dprime\left(\frac{1-s}{q}\right)}\right)}\right)\\
                  +\frac{1-p_q}{2q}\cN^{+} \left(0,-\frac{q^2}{(q-1)^2}f^\dprime_{\beta,h}(s)\right)+\frac{1+p_q}{2}\delta_0.
              \end{multline*}
    \end{enumerate}
\end{thm}

Theorem \ref{hml77} is proved in Appendix \ref{hml77pr}. It shows that at the critical points, the limiting distribution of $\hat{h}_N$ is a mixture distribution consisting of half-normal distributions and a point mass at $0$. In particular, $\hat{h}_N$ is always $\sqrt{N}$-consistent at the critical points. We now shift our attention to the asymptotics of $\hat{\beta}_N$.

\begin{thm}[Asymptotic distributions of $\hat{\beta}_N$ at regular points]\label{betaml7}
    Fix $p \geq 2$ and suppose $(\beta, h) \in$ $\Theta$ is regular. Assume $\beta$ is known and $\bm X \sim \P_{\beta, h, p}$. Then denoting the unique maximizer of $H$ by $\bm m_*=\bm m_*(\beta, h, p)$, as $N \rightarrow \infty$
    \begin{enumerate}
        \item If $h>0$, then $\bm m_* \neq \bm x_0$, and
              \begin{equation}
                  N^{\frac{1}{2}}\left(\hat{\beta}_N-\beta\right) \xrightarrow{D} \cN\left(0,-\frac{q^2f^\dprime_{\beta,h}(s)}{p^2(q-1)^2}\left(m_1^{p-1}-m_2^{p-1}\right)^{-2}\right),
                  \label{mle_b_reg1}
              \end{equation}
        \item If $h=0$, then $\bm m_* = \bm x_0$ and
              \begin{equation*}
                  N^{\frac{1}{2}}\left(\hat{\beta}_N-\beta\right) \xrightarrow{D} \gamma_1 \delta_{-\infty} + (1-\gamma_1) \delta_\infty,
              \end{equation*}
              where $\gamma_1:=\P\left(\bm W^\top \bm W \leq \frac{1-q}{k^\dprime\left(\frac{1}{q}\right)}\right)$ with $\bm W \sim \cN_q(\boldsymbol{0},\Sigma)$.
    \end{enumerate}
\end{thm}

Theorem \ref{betaml7} is proved in Appendix \ref{prbetaml7}. It shows that $\hat{\beta}_N$ is $N^{\frac{1}{2}}$-consistent and asymptotically normal at the regular points when the maximizer is not $\bm x_0$, whereas if the maximizer happens to be $\bm x_0$, then $N^{\frac{1}{2}}(\hat{\beta}_N-\beta)$ is inconsistent.

\begin{thm}[Asymptotic distributions of $\hat{\beta}_N$ at special points]\label{asmbetasp}
    Fix $p \geq 2$ and suppose $(\beta, h) \in$ $\Theta$ is special. Assume $\beta$ is known and $\epm \sim \P_{\beta, h, p}$. Denote the unique maximizer of $H$ by $\bm m_*=\bm m_*(\beta, h, p)$.
    \begin{enumerate}
        \item If $(\beta, h) \in$ $\Theta$ is type I special then, as $N \rightarrow \infty$,
              \begin{itemize}
                  \item if $(p,q)\notin \{(2,2)\}\cup \{(3,2)\}$,
                        $$
                            N^{\frac{3}{4}}\left(\hat{\beta}_N-\beta\right) \xrightarrow{D} L_1
                        $$
                        where the distribution function of $L_1$ is given by
                        $$
                            L_1(t)= F_{0,0}\left(-\int_{-\infty}^{\infty} u \mathrm{~d} F_{t, 0}(u)\right),
                        $$
                        with $T_{t, 0}$ as defined in \eqref{spec1_dist} below.
                  \item if $(p,q)= (2,2)$ or $(p,q)= (3,2)$ then,
                        \begin{equation*}
                            N^{\frac{3}{4}}\left(\hat{\beta}_N-\beta\right) \xrightarrow{D} \alpha\delta_{-\infty} + (1-\alpha) \delta_\infty.
                        \end{equation*}
                        where $\alpha := \p(T_{0,0}^2 \le \e T_{0,0}^2)$.
              \end{itemize}
        \item If $(\beta, h) \in$ $\Theta$ is type II special then, as $N \rightarrow \infty$,
              $$
                  N^{\frac{5}{6}}\left(\hat{\beta}_N-\beta\right) \xrightarrow{D} \gamma_2 \delta_{-\infty} + (1-\gamma_2) \delta_\infty
              $$
              where $\gamma_2:= \P(F_0^2\leq \e F_0^2)$.
    \end{enumerate}
\end{thm}

Once again, we skip the proof of Theorem \ref{asmbetasp} due to its very close similarity with the proof of Theorem \ref{betaml7}. Finally, we state the result about the asymptotics of $\hat{\beta}_N$ at the critical points.

\begin{thm}[Asymptotic distributions of $\hat{\beta}_N$ at critical points]\label{btml44}
 Suppose that $(\beta, h)$ is a critical point. Let $p_1,\ldots,p_K$ be the weights defined in the statement of Theorem \ref{pcr} for the global maximizers $\bm m_1,\ldots,\bm m_K$, respectively, where these maximizers are arranged in ascending order of their $L^p$ norms. Then, for $\bm X \sim \P_{\beta, h, p}$, as $N \rightarrow \infty$, we have the following:
 
    \begin{enumerate}
        \item If $(\beta,h) \in \cC_{p,q}^1\setminus \{(\beta_c,0)\}$, then $f_{\beta,h}$ has exactly two global maximizers $s_2>s_1>0$, and
              \begin{eqnarray*}
                  N^{\frac{1}{2}}\left(\hat{\beta}_N-\beta\right) &\xrightarrow{D}& \frac{p_1}{2}\cN^{-} \left(0,-\frac{q^2f^\dprime_{\beta,h}(s_1)}{p^2(q-1)^2}\left(m_{1,1}^{p-1}-m_{1,2}^{p-1}\right)^{-2}\right)\\ &+&\frac{1-p_1}{2}\cN^{+} \left(0,-\frac{q^2f^\dprime_{\beta,h}(s_2)}{p^2(q-1)^2}\left(m_{2,1}^{p-1}-m_{2,2}^{p-1}\right)^{-2}\right)
                  +\frac{1}{2}\delta_0
              \end{eqnarray*}

              \item If $(\beta,h) \in \cC_{p,q}^2$, then $f_{\beta,h}$ has exactly one global maximizer $s>0$, and
              $$N^{\frac{1}{2}}\left(\hat{\beta}_N - \beta\right) \xrightarrow{D} \mathcal{N}\left(0, \frac{q^2 f_{\beta,h}''(s)}{p^2 (q-1)^2} (x_{s,1}^{p-1} - x_{s,2}^{p-2})^{-2}\right)$$

        \item If $(\beta,h) =(\beta_c,0)$, then $f_{\beta,h}$ has exactly two maximizers, 0 and $s>0$, and
              \begin{eqnarray*}
                  N^{\frac{1}{2}}\left(\hat{\beta}_N-\beta\right) &\xrightarrow{D}& p_1\gamma_1\delta_{-\infty}+\frac{1-p_1}{2}\cN^{+} \left(0,-\frac{q^2f^\dprime_{\beta,h}(s)}{p^2(q-1)^2}\left(x_{s,1}^{p-1}-x_{s,2}^{p-1}\right)^{-2}\right)\\
                  &+&\left(\frac{1+p_1}{2}-p_1\gamma_1\right)\delta_0
              \end{eqnarray*}
              where $\gamma_1$ is as defined in the statement of Theorem \ref{betaml7} (2).
    \end{enumerate}
\end{thm}

Theorem \ref{btml44} is proved in Appendix \ref{btml44pr}. It says that as long as $(\beta,h) \ne (\beta_c,0)$, $\hat{\beta}_N$ is $\sqrt{N}$-consistent, and its asymptotic distribution is either a mixture of half-normals and a point mass at $0$, or just a normal, depending on whether the point is strongly or weakly critical, respectively. However, if $(\beta,0)=(\beta_c,0)$, then $\hat{\beta}_N$ is no longer $\sqrt{N}$-consistent, and a portion of the asymptotic mass escapes to $-\infty$. The last phenomenon can be explained by the fact that for $h=0$, if $\beta<\beta_c$, $\sqrt{N}(\hat{\beta}_N-\beta)$ does not have any asymptotic finite mass, and for $\beta>\beta_c$,  $\hat{\beta}_N$ is $\sqrt{N}$ consistent, so at the transition point $\beta_c$, a portion of the asymptotic mass of $\sqrt{N}(\hat{\beta}_N-\beta)$ is finite, and the remaining mass stays at $-\infty$.  

\section{Confidence Intervals for the Model Parameters}\label{sec:confint}
In this section, we start by summarizing the partition of the parameter space into different components, induced by the function $H_{\beta,h}$. This summary is a consequence of the results proved in Appendix \ref{apxf}. The existence of this partition and the different forms of the limiting distributions of the ML estimates on the different components of this partition gives rise to an inherent difficulty in constructing confidence intervals for the model parameters. In this context, there are two different scenarios:

\begin{figure}[ht]
\centering
\centering
\includegraphics[width=4in]
    {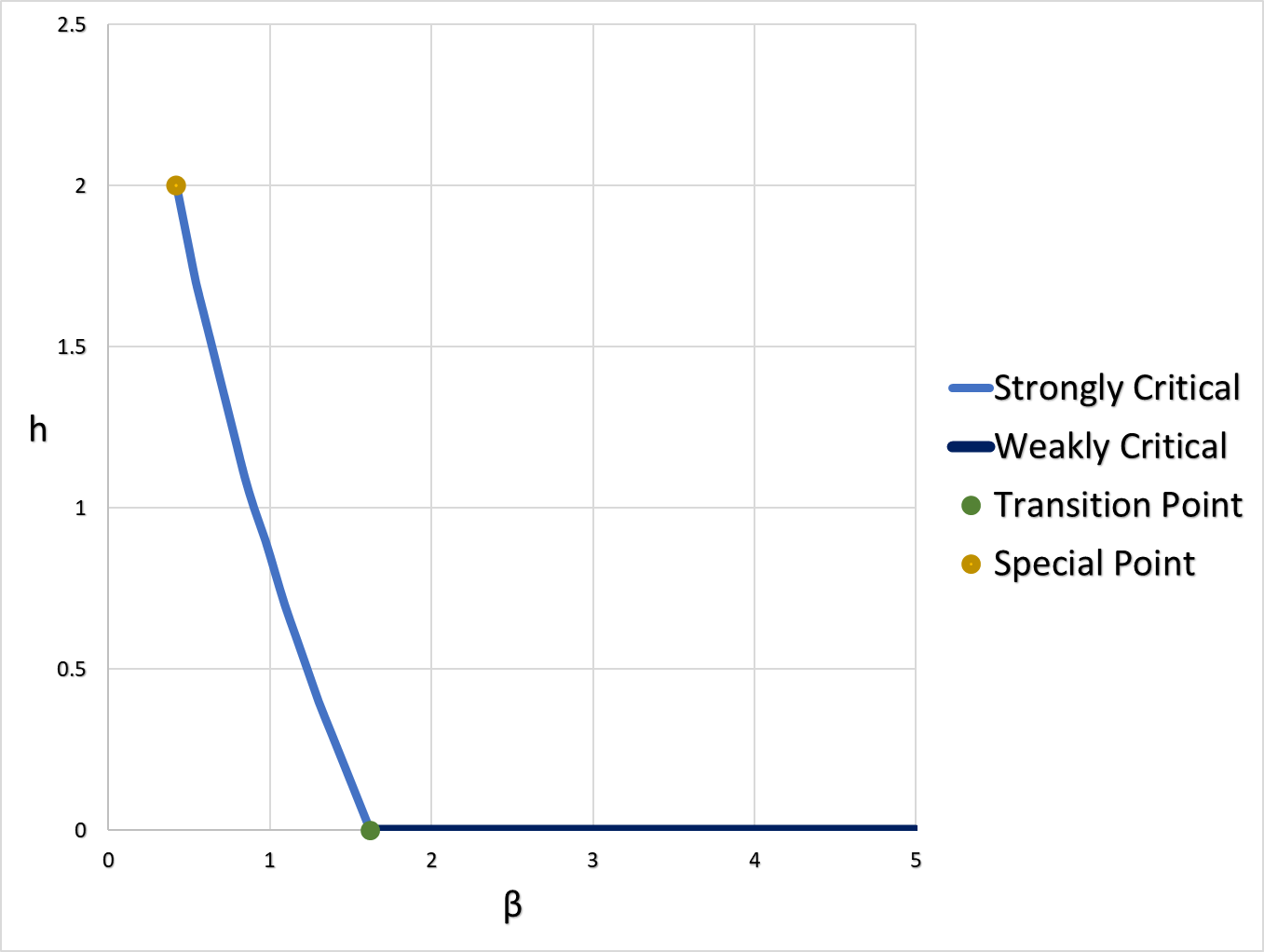}\\
\caption{\small{Phase diagram for the case $(p,q)=(7,5)$. The light blue curve denotes the set of strongly critical points, the deep blue line denotes the set of weakly critical points, the golden point denotes the special point (which in this case is of type-I), and the green point denotes the transition point $\beta_c$. The white region, which is the complement of all these colored curves, lines and points, is the set of regular points.}}
\label{fig1}
\end{figure}

\begin{figure}[ht]
\centering
\centering
\includegraphics[width=4in]
    {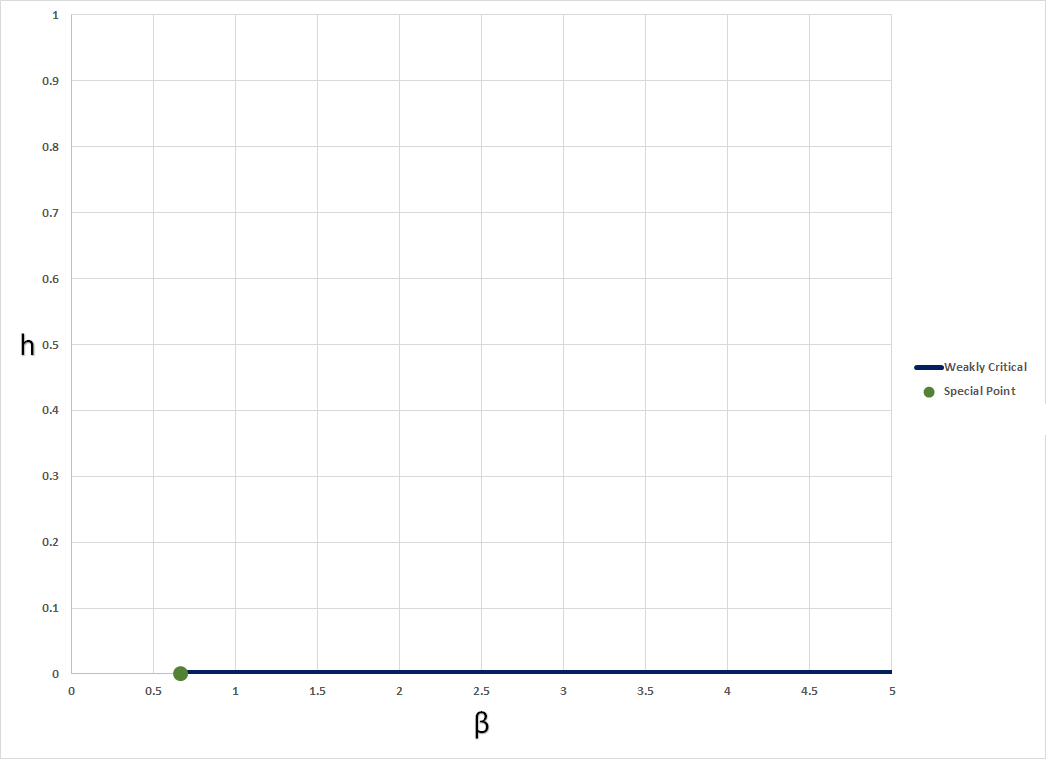}\\
\caption{\small{Phase diagram for the case $(p,q)=(4,2)$. The deep blue line denotes the set of weakly critical points, the green point denotes the special point, which is of type II. The white region, which is the complement of these two sets, is the set of regular points.}}
\label{fig2}
\end{figure}

\begin{enumerate}
    \item \boldsymbol{$p\ge 5, q\ne 2$}: In this case, the only special point in the parameter space $(\wbt,\wht)$ lies in $(0,\infty)\times (0,\infty)$. This point is type-I special. The set $\cC_{p,q}^1$ is a smooth, strictly decreasing curve starting from the point $(\wbt,\wht)$ (excluding it), and continuing till a point $(\beta_c(p,q),0)$ (including it). The set $\cC_{p,q}^2$ is the interval $\{(\beta,0): \beta > \beta_c(p,q)\}$. The remaining portion of the parameter space $\Theta$ is the set of all regular points.
\vspace{0.05in}

    \item \boldsymbol{$p\in \{2,3,4\}, q= 2$}: In this case, any point $(\beta,h)$ with either $h>0$ or $\beta<\beta_c(p,q) = \frac{2^{p-1}}{p(p-1)}$ is a regular point. The point $(\beta_c(p,q),0)$ is the unique special point, which is of type-I if $p\in \{2,3\}$, and type-II if $p=4$. The remaining portion of $\Theta$, i.e. the  interval $\{(\beta,0): \beta > \beta_c(p,q)\}$ is the set $\cC_{p,q}^2$. Consequently, $\cC_{p,q}^1 =\emptyset$ in this case.
\end{enumerate}

In Figure \ref{fig1}, we illustrate this partition for the case $p=7, q=5$, and in Figure \ref{fig2}, for the case $p=4,q=2$, through phase diagrams. 

We now discuss how to construct confidence intervals for the model parameters $\beta$ and $h$, with asymptotic coverage probability $1-\alpha$. This is not a direct task, since the asymptotics of the ML estimates depend upon the exact position of the true $(\beta,h)$ in $\Theta$. However, intuitively speaking, since the complement of the set of regular points has Lebesgue measure $0$, it should be enough to just use the limiting distributions at the regular points to construct the confidence intervals for the model parameters. So, let us imagine that an oracle told us beforehand that the unknown parameter $(\beta,h)$ is regular. Then, the intervals:
$$I:= \left(\hat{h}_N - \frac{q}{q-1}\sqrt{-\frac{f_{\beta,0}''(1-q \bar{X}_{\cdot q})}{N}} ~z_{1-\frac{\alpha}{2}}~,~\hat{h}_N + \frac{q}{q-1}\sqrt{-\frac{f_{\beta,0}''(1-q \bar{X}_{\cdot q})}{N}} ~z_{1-\frac{\alpha}{2}}\right)$$

$$J := \left(\hat{\beta}_N - \frac{qN^{-\frac{1}{2}}\sqrt{-f_{\hat{\beta}_N,0}''(1-q\bar{X}_{\cdot q})}}{p(q-1)\left(\bar{X}_{\cdot 1}^{p-1}-\bar{X}_{\cdot 2}^{p-1}\right)}~z_{1-\frac{\alpha}{2}}~,~\hat{\beta}_N + \frac{qN^{-\frac{1}{2}}\sqrt{-f_{\hat{\beta}_N,0}''(1-q\bar{X}_{\cdot q})}}{p(q-1)\left(\bar{X}_{\cdot 1}^{p-1}-\bar{X}_{\cdot 2}^{p-1}\right)}~z_{1-\frac{\alpha}{2}}\right)$$
are asymptotic $(1-\alpha)$-coverage confidence intervals for $h$ given $\beta$, and $\beta$ given $h\ne 0$, respectively.

We now discuss how to modify the intervals $I$ and $J$ to asymptotically valid confidence sets at all points. Towards this, for every $\beta$, let $S(\beta)$ be the set of all $h$, such that $(\beta,h)$ belongs to the closure of the set $\cC_{p,q}$, and for every $h\ne 0$, let $T(h)$ be the set of all $\beta$, such that $(\beta,h)$ belongs to the closure of the set $\cC_{p,q}$. Note that both $S(\beta)$ and $T(h)$ have cardinality at most $1$. Clearly, $I\bigcup S(\beta)$ and $J\bigcup T(h)$ are asymptotically level $1-\alpha$ confidence sets for $h$ given $\beta$ and $\beta$ given $h\ne 0$, respectively, which have the same Lebesgue measure as the intervals $I$ and $J$, respectively.

There is an alternative, more precise two-step algorithm one can follow, than just uniting the points on the closure of the critical curve to $I$ and $J$ as described above, to get the universally valid confidence intervals. For fixed $\beta$, one can first consistently test the null hypothesis $H_0: h \in S(\beta)$ at level $\alpha$ using the asymptotic distribution of $\hat{h}_N$ at the critical or special points. If this null is rejected, then he can report $I$ as the confidence interval for $h$, and otherwise, he can declare the singleton set $S(\beta)$ as the confidence interval (which is either empty, or just a point). A similar approach can be followed for constructing the confidence interval for $\beta$ also, where this time, one tests the null hypothesis $H_0: \beta \in T(h)$ in the first step, and if this is accepted, reports $T(h)$ as the confidence interval for $\beta$, and $J$ otherwise.

\section{Sketch of Proof}\label{prskt}
In this section, we provide a brief sketch of the proofs of the main results in this paper. We begin with the proof of the asymptotics of the magnetization vector.  The first step towards this, is to show that the magnetization vector concentrates around the set of all global maximizers of the function $H_{\beta,h}$, which makes them natural candidates for centering in the central limit theorems. The next step is to show that conditional on the event that $\epm $ is some neighborhood of a global maximizer $\bm m_*$ whose closure is devoid of any other maximizer, every bounded, continuous function $g:\mathbb{R}^q \to \mathbb{R}$, satisfies:
$$\e\left[g\left(\sqrt{N}(\epm - \bm m_*)\right)\mathbbm{1}_{\{\sqrt{N}(\epm - \bm m_*) \le M\}}\right] \rightarrow \e\left[g(Y) \mathbbm{1}_{\{Y\le M\}}\right]$$ where $Y$ follows the law of the appropriate limiting distribution (which is either a Gaussian, or a fourth-order or sixth-order Gaussian, depending on whether the true parameter is regular/critical or special). A subsequent uniform integrability argument for all moments of $\sqrt{N}(\epm - \bm m_*)$ will now imply its weak convergence and convergence in all moments to $Y$. With the vision of applying these results to derive the asymptotics of the ML estimates, we prove these convergence results under slightly perturbed versions of the true parameters, the perturbations being of the order $N^{-1/2}$.

Next, for proving asymptotics of the ML estimates, using monotonicity of the functions $u_{N,1}$ and $u_{N,p}$, one can express the cumulative distributions of $\sqrt{N}(\hat{h}_N - h)$ and $\sqrt{N}(\hat{\beta}_N-\beta)$ in terms of the cumulative distribitions of $\sqrt{N}(\bar{X}_{\cdot 1} - m_{*1})$ and $\sqrt{N}(\|\epm\|_p^p - \|\bm m_*\|_p^p)$ at their respective expectations under the perturbed parameters. This then enables one to translate the asymptotic results of $\epm$ to asymptotics of the ML estimates. Some care needs to be cautioned at critical points where there are more than one maximizer, but in that case, the leaning of $\epm$ towards some particular maximizers and away from the others, is largely governed by the sign of the perturbation of the true parameters, which is made rigorous through some perturbative concentration results proved in Appendix \ref{sec:partsec}.

\section{Acknowledgment} 
\noindent S. Mukherjee was supported by the National University of Singapore start-up grant WBS A-0008523-00-00 and the FoS Tier 1 grant WBS A-8001449-00-00.

\appendix

\section{Proof of Theorem \ref{conclem}}\label{prclem}
 In this section, we prove Theorem \ref{conclem}. Towards this, let $A_N(\bm v) := \{\bm x \in [q]^N: \bar{\bm x}_N = \bm v\}$, and $S_N := \{0,\frac{1}{N},\frac{2}{N},\ldots, 1\}$. Then, for any set $G \subseteq \R^q$, we have by Lemma \ref{amcard}:
    \begin{eqnarray*}
        &&\p_{\beta_N,h_N,N}\left(\epm \in G\right)\\ &=& \frac{\sum_{\bm v \in S_N^q\bigcap \cP_q\bigcap G} |A_N(\bm v)|  \exp\left\{N \left(\beta_N \sum_{r=1}^q v_r^p + h_Nv_1\right)\right\}}{\sum_{\bm v \in S_N^q\bigcap \cP_q} |A_N(\bm v)|  \exp\left\{N \left(\beta_N \sum_{r=1}^q v_r^p +  h_Nv_1\right)\right\}}\\&= & e^{o(N)}\frac{\sum_{\bm v \in S_N^q\bigcap \cP_q\bigcap G} |A_N(\bm v)|  \exp\left\{N \left(\beta \sum_{r=1}^q v_r^p + hv_1\right)\right\}}{\sum_{\bm v \in S_N^q\bigcap \cP_q} |A_N(\bm v)|  \exp\left\{N \left(\beta \sum_{r=1}^q v_r^p +  hv_1\right)\right\}} \\&\lesssim_q & \frac{e^{o(N)} N^{\frac{1}{2}}(N+1)^q\sup_{\bm v \in S_N^q\bigcap\cP_q \bigcap G} \exp\left\{N H_{\beta,h}(\bm v)\right\} }{\sup_{\bm v \in S_N^q\bigcap\cP_q } \exp\left\{N H_{\beta,h}(\bm v)\right\}}\\&\le& e^{o(N)}\exp\left\{ N \left(\sup_{\bm v \in G\bigcap\cP_q} H_{\beta,h}(\bm v) - \sup_{\bm v \in S_N^q\bigcap\cP_q} H_{\beta,h}(\bm v)  \right)   \right\}~.
    \end{eqnarray*}
    Now, note that:
    $$\limsup_{N\rightarrow \infty}\sup_{\bm v \in S_N^q\bigcap\cP_q} H_{\beta,h}(\bm v) \le  \sup_{\bm v \in \cP_q} H_{\beta,h}(\bm v)~. $$ On the other hand, for any maximizer $\bm m_*(\beta,h)$ of $H_{\beta,h}$, Lemma \ref{denselem} gives a sequence $\bm v_N \in S_N^q\bigcap\cP_q$ such that $\bm v_N \rightarrow \bm m_*(\beta,h)$. This shows that:
    $$\liminf_{N\rightarrow \infty} \sup_{\bm v \in S_N^q\bigcap\cP_q} H_{\beta,h}(\bm v) \ge  \liminf_{N\rightarrow \infty} H_{\beta,h}(\bm v_N) = \hf(\bm m_*(\beta,h)) =  \sup_{\bm v \in \cP_q} H_{\beta,h}(\bm v)~. $$
    Hence, as $N\rightarrow \infty$, we have:
    $$\sup_{\bm v \in S_N^q\bigcap\cP_q} H_{\beta,h}(\bm v) \rightarrow \sup_{\bm v \in \cP_q} H_{\beta,h}(\bm v)~,$$ which implies that:
    \begin{equation*}
        \p_{\beta_N,h_N,N}\left(\epm \in G\right) \lesssim_q e^{o(N)}\exp\left\{ N \left(\sup_{\bm v \in G\bigcap\cP_q} H_{\beta,h}(\bm v) - \sup_{\bm v \in \cP_q} H_{\beta,h}(\bm v) +o(1) \right)   \right\}~.
    \end{equation*}
    Hence, for every $G \subseteq \R^q$, we have:
    \begin{equation}\label{ldpup}
        \limsup_{N\rightarrow \infty} \frac{1}{N} \log \p_{\beta_N,h_N,N}\left(\epm \in G\right) \le \sup_{\bm v \in G\bigcap\cP_q} H_{\beta,h}(\bm v) - \sup_{\bm v \in \cP_q} H_{\beta,h}(\bm v)~.
    \end{equation}

    Next, for any set $G \subseteq \R^q$, we also have by Lemma \ref{amcard}:

    \begin{eqnarray*}
        &&\p_{\beta_N,h_N,N}\left(\epm \in G\right)\\ &=& \frac{\sum_{\bm v \in S_N^q\bigcap \cP_q\bigcap G} |A_N(\bm v)|  \exp\left\{N \left(\beta_N \sum_{r=1}^q v_r^p + h_Nv_1\right)\right\}}{\sum_{\bm v \in S_N^q\bigcap \cP_q} |A_N(\bm v)|  \exp\left\{N \left(\beta_N \sum_{r=1}^q v_r^p +  h_Nv_1\right)\right\}}\\&=& e^{o(N)}\frac{\sum_{\bm v \in S_N^q\bigcap \cP_q\bigcap G} |A_N(\bm v)|  \exp\left\{N \left(\beta \sum_{r=1}^q v_r^p + hv_1\right)\right\}}{\sum_{\bm v \in S_N^q\bigcap \cP_q} |A_N(\bm v)|  \exp\left\{N \left(\beta \sum_{r=1}^q v_r^p +  hv_1\right)\right\}}\\&\gtrsim_q & \frac{e^{o(N)}\sup_{\bm v \in S_N^q\bigcap\cP_q \bigcap G} \exp\left\{N H_{\beta,h}(\bm v)\right\} }{N^{\frac{1}{2}}(N+1)^q\sup_{\bm v \in S_N^q\bigcap\cP_q } \exp\left\{N H_{\beta,h}(\bm v)\right\}}\\&\ge& e^{o(N)}\exp\left\{ N \left(\sup_{\bm v \in G\bigcap\cP_q\bigcap S_N^q} H_{\beta,h}(\bm v) - \sup_{\bm v \in \cP_q} H_{\beta,h}(\bm v)  \right)   \right\}~.
    \end{eqnarray*}
    Once again, note that:
    $$\limsup_{N\rightarrow \infty}\sup_{\bm v \in G\bigcap S_N^q\bigcap\cP_q} H_{\beta,h}(\bm v) \le  \sup_{\bm v \in G\bigcap \cP_q} H_{\beta,h}(\bm v)~. $$ Let $\varepsilon>0$ be given. Then, assuming $G\bigcap \cP_q\ne \emptyset$, one can choose $\bm m \in G\bigcap \cP_q$ such that $H_{\beta,h}(\bm m) > \sup_{\bm v \in G\bigcap \cP_q}  H_{\beta,h}(\bm v)-\varepsilon~. $  Lemma \ref{denselem} gives a sequence $\bm v_N \in S_N^q\bigcap\cP_q$ such that $\bm v_N \rightarrow \bm m$. If $G$ is assumed to be open, then $\bm v_N \in G\bigcap S_N^q\bigcap\cP_q$ eventually, and hence,
    $$\liminf_{N\rightarrow \infty} \sup_{\bm v \in G\bigcap S_N^q\bigcap\cP_q} H_{\beta,h}(\bm v) \ge  \liminf_{N\rightarrow \infty} H_{\beta,h}(\bm v_N) = \hf(\bm m) =  \sup_{\bm v \in G\bigcap \cP_q} H_{\beta,h}(\bm v) - \varepsilon~. $$
    Since $\varepsilon >0$ is arbitrary, we conclude that:
    $$\sup_{\bm v \in G\bigcap S_N^q\bigcap\cP_q} H_{\beta,h}(\bm v) \rightarrow \sup_{\bm v \in G\bigcap \cP_q} H_{\beta,h}(\bm v)~,$$
    which implies that:
    \begin{equation*}
        \p_{\beta_N,h_N,N}\left(\epm \in G\right) \gtrsim_q e^{o(N)}\exp\left\{ N \left(\sup_{\bm v \in G\bigcap\cP_q} H_{\beta,h}(\bm v) - \sup_{\bm v \in \cP_q} H_{\beta,h}(\bm v) +o(1) \right)   \right\}~.
    \end{equation*}
    Hence, for every open set $G \subseteq \R^q$, we have:

    \begin{equation}\label{ldpdown}
        \liminf_{N\rightarrow \infty} \frac{1}{N} \log \p_{\beta_N,h_N,N}\left(\epm \in G\right) \ge \sup_{\bm v \in G\bigcap\cP_q} H_{\beta,h}(\bm v) - \sup_{\bm v \in \cP_q} H_{\beta,h}(\bm v)~.
    \end{equation}
    The large deviation principle of $\epm$ now follows from \eqref{ldpup} and \eqref{ldpdown}.

    Now, take $G:= \{\bm t \in \R^q: d(\bm t, \cM_{\beta,h}) \ge \varepsilon\}$. Then,
    $\sup_{\bm v \in G\bigcap\cP_q} H_{\beta,h}(\bm v) < \sup_{\bm v \in \cP_q} H_{\beta,h}(\bm v)$.  It thus follows from the large deviation principle of $\epm$ (or directly from \eqref{ldpup}), that
    $$\limsup_{N\rightarrow \infty} \frac{1}{N} \log \p_{\beta_N,h_N,N}\left(\epm \in G\right) < 0~,$$
    which completes the proof of Theorem \ref{conclem}.

\section{Proofs of the Asymptotics of the Magnetization}
In this section, we prove the results on the asymptotics of the magnetization when the sample is coming from a $p$-tensor Potts model with perturbed parameters. Some of these proofs closely follow the proofs in \cite{gan}.

\subsection{Proof of Theorem \ref{cltreg}} \label{prcltreg}
In this section, we prove Theorem \ref{cltreg}.
 Towards this, denote $\bm W_N := \sqrt{N}\left(\epm - \bm m_*\right)$. Fix a positive real number $M$ and a bounded, continuous function $g: \R^q \to \R$. For every $\bm v \in \cP_{q,N} := \cP_q\bigcap S_N^q$, define $\bm w(\bm v) = \bm w_N(\bm v) := \sqrt{N}(\bm v - \bm m_*)$. Then, we have by Lemma \ref{densappr},
    \begin{eqnarray}
        &&q^N Z_N(\beta_N,h_N) \e_{\beta_N,h_N,N}\left[g(\bm W_N) \mathbbm{1}_{\|\bm W_N\| \le M}\right]\nonumber\\ &=&  \sum_{\bm v \in \cP_{q,N}} g(\bm w(\bm v)) \mathbbm{1}_{\|\bm w(\bm v)\| \le M} q^N Z_N(\beta_N,h_N) \p_{\beta_N,h_N,N}(\epm = \bm v)\nonumber\\&=& (1+o_N(1)) N^{-\frac{q-1}{2}}  \sum_{\bm v \in \cP_{q,N}} A(\bm v) e^{N H_{\beta_N,h_N}(\bm v)} g(\bm w(\bm v)) \mathbbm{1}_{\|\bm w(\bm v)\| \le M}\nonumber\\&=& (1+o_N(1)) N^{-\frac{q-1}{2}}  \sum_{\bm v \in \cP_{q,N}} A\left(\bm m_* + N^{-\frac{1}{2}} \bm w(\bm v)\right) e^{N H_{\beta_N,h_N}\left(\bm m_* + N^{-\frac{1}{2}} \bm w(\bm v)\right)}\nonumber\\ && g(\bm w(\bm v)) \mathbbm{1}_{\|\bm w(\bm v)\| \le M}\nonumber\\&=& (1+o_N(1)) N^{-\frac{q-1}{2}}  A(\bm m_*) \sum_{\bm v \in \cP_{q,N}} e^{N H_{\beta_N,h_N}\left(\bm m_* + N^{-\frac{1}{2}} \bm w(\bm v)\right)}g(\bm w(\bm v)) \mathbbm{1}_{\|\bm w(\bm v)\| \le M}\label{stp1}~.
    \end{eqnarray}
    By Lemma \ref{l2nm} we get that,
    \begin{multline}
        H_{\beta_{N}, h_{N}}\left(\bm m_*+N^{-\frac{1}{2}} \bm w(\bm v)\right)=H_{\beta, h}\left(\bm m_*+N^{-\frac{1}{2}} \bm w(\bm v)\right)+\frac{\bar{\beta}}{\sqrt{N}}\|\bm m_*\|_{p}^p+\frac{\bar{h}}{\sqrt{N}} m_{N,1}\\
        +\ip{\bar{\beta}p\bm m_*^{p-1}+\bar{h}\bm e_1}{\bm w(\bm v)}\frac{1}{N}+o\left(N^{-1}\right)
        \label{bnhn_to_bh}
    \end{multline}
    By Lemma \ref{sectay}, we have the following on the event $\{\|\bm w(\bm v)\|\le M\}$ for all large $N$:
    \begin{eqnarray*}
        &&H_{\beta,h}\left(\bm m_* + N^{-\frac{1}{2}} \bm w(\bm v)\right)\\&=& H_{\beta,h}(\bm m_*) +\frac{1}{2N} \sum_{r=1}^q \left(\beta_N p(p-1)t_{N,r,\alpha}^{p-2} -\frac{1}{t_{N,r,\alpha}}\right) w_r^2(\bm v) ~,
    \end{eqnarray*}
    where $t_{N,r,\alpha}$ denotes the $r^{\text{th}}$ element of $\bm m_* + \alpha N^{-\frac{1}{2}} \bm w(\bm v)$ for some $\alpha \in [0,1]$ that can depend on $N$. Therefore, on noting that $t_{N,r,\alpha} = m_{*,r} + o_N(1)$ on the event $\{\|\bm w(\bm v)\| \le M\}$, we have the following on the event $\{\|\bm w(\bm v)\| \le M\}$:
    \begin{eqnarray}
        && NH_{\beta,h}\left(\bm m_* +   N^{-\frac{1}{2}} \bm w(\bm v)\right)\nonumber\\ &=& NH_{\beta,h}(\bm m_*) + \frac{1}{2} \bm Q_{\bm m_*,\beta}(\bm w (\bm v)) + o_N(1)\label{mn_to_m}~.
    \end{eqnarray}
    Now putting back \eqref{bnhn_to_bh} and \eqref{mn_to_m} together and using Lemma \ref{l2nm} we get that,
    \begin{equation}
        NH_{\beta_{N}, h_{N}}\left(\bm m_*+N^{-\frac{1}{2}} \bm w(\bm v)\right)=NH_{\beta_N,h_N}(\bm m_*) + \frac{1}{2} \bm Q_{\bm m_*,\beta}(\bm w (\bm v))
        +\ip{\bar{\beta}p\bm m_*^{p-1}+\bar{h}\bm e_1}{\bm w(\bm v)}+o\left(1\right).
        \label{final_reduc_reg}
    \end{equation}
    It thus follows from \eqref{stp1} that:
    \begin{eqnarray*}
        &&q^N Z_N(\beta_N,h_N) \e_{\beta_N,h_N,N}\left[g(\bm W_N) \mathbbm{1}_{\|\bm W_N\| \le M}\right]\nonumber\\ &=& (1+o_N(1)) N^{-\frac{q-1}{2}}  A(\bm m_*) e^{N H_{\beta_N,h_N}\left(\bm m_*\right)}\sum_{\bm v \in \cP_{q,N}} g(\bm w(\bm v)) \mathbbm{1}_{\|\bm w(\bm v)\| \le M} e^{\ip{\bar{\beta}p\bm m_*^{p-1}+\bar{h}\bm e_1}{\bm w(\bm v)}+\frac{1}{2} \bm Q_{\bm m_*,\beta}(\bm w (\bm v))}~.
    \end{eqnarray*}

    Hence by Riemann sum approximation, we get:
    \begin{equation}
        \label{cont_approx}
        \begin{aligned}
             & q^N Z_N(\beta_N,h_N) \e_{\beta_N,h_N,N}\left[g(\bm W_N) \mathbbm{1}_{\|\bm W_N\| \le M}\right]                                                                                                                                                 \\
             & \sim C_{g,M}A(\bm m_*) e^{N H_{\beta_N,h_N}\left(\bm m_*\right)} \int_{\cH_q \bigcap B(\bm 0,M)} g(\bm w) e^{\ip{\bar{\beta}p\bm m_*^{p-1}+\bar{h}\bm e_1}{\bm w}+\frac{1}{2}  \bm Q_{\bm m_*,\beta} (\bm w)}~dw_1dw_2\ldots dw_q.
        \end{aligned}
    \end{equation}
    Therefore, we have:
    \begin{equation*}
        \begin{aligned}
            \e_{\beta_N,h_N,N}\left[g(\bm W_N) \mathbbm{1}_{\|\bm W_N\| \le M}\right] \propto \int_{\cH_q \bigcap B(\bm 0,M)} g(\bm w) e^{\ip{\bar{\beta}p\bm m_*^{p-1}+\bar{h}\bm e_1}{\bm w}+\frac{1}{2}  \bm Q_{\bm m_*,\beta} (\bm w)}~dw_1dw_2\ldots dw_q.
        \end{aligned}
    \end{equation*}
    Hence, under $\p_{\beta_N,h_N,N}$, $\bm W_N$ conditioned on $\|\bm W_N\| \le M$ converges weakly to the density on $\cH_q \bigcap B(0,M)$ with density (with respect to the Lebesgue measure on $\cH_q$) proportional to $$\bm w \mapsto e^{\ip{\bar{\beta}p\bm m_*^{p-1}+\bar{h}\bm e_1}{\bm w}+\frac{1}{2}\bm Q_{\bm m_*,\beta}(\bm w)}~,$$ where $\cH_q := \{\bm t \in \R^q: \sum_{r=1}^q t_r = 0\}$.

    Next, we show that $\bm W_N$ is uniformly integrable under $\p_{\beta_N,h_N,N}$. Let us first break down $\E_{\beta_N,h_N,N}[\nrm{\bm W_N}^r\one_{\|\bm W_N\| \ge K}]$ as,
    \begin{equation*}
        \E_{\beta_N,h_N,N}\left(\nrm{\bm W_N}^r\one_{\|\bm W_N\| \ge K}\right)=Z_1+Z_2,
    \end{equation*}
    where
    \begin{equation*}
        \begin{aligned}
            Z_1&=\E_{\beta_N,h_N,N}\left(\nrm{\bm W_N}^r\one_{\|\bm W_N\| \ge K} \mid \|\bm W_N\| \le \varepsilon\sqrt{N} \right) \p_{\beta_N,h_N,N}\left(\|\bm W_N\| \le \varepsilon\sqrt{N}\right)\\
            Z_2&=\E_{\beta_N,h_N,N}\left(\nrm{\bm W_N}^r\one_{\|\bm W_N\| \ge K} \mid \|\bm W_N\| \ge \varepsilon\sqrt{N} \right) \p_{\beta_N,h_N,N}\left(\|\bm W_N\| \ge \varepsilon\sqrt{N}\right).
        \end{aligned}
    \end{equation*}
    Now, $\nrm{\bm W_N}^r=O_P(N^{\frac{r}{2}})$ and in view of Theorem \ref{conclem}, $\p_{\beta_N,h_N,N}\left(\|\bm W_N\| \ge \varepsilon\sqrt{N}\right)=O(\exp(-C_{q,\varepsilon}N))$. Hence, $Z_1 \rightarrow 0$ as $N\rightarrow \infty$.
    So, it suffices to show that,
    \begin{equation}\label{ntsh}
        \lim_{K\rightarrow \infty} \limsup_{N\rightarrow \infty}\E_{\beta_N,h_N,N}\left(\nrm{\bm W_N}^r\one_{\|\bm W_N\| \ge K} \Big| \|\bm W_N\| \le \varepsilon\sqrt{N}\right) = 0
    \end{equation}
    Towards this, it follows from Lemmas \ref{densappr} and \ref{sectay2}, that for $\varepsilon >0$ small enough and $K>0$,

    \begin{eqnarray*}
        &&\E_{\beta_N,h_N,N}\left(\nrm{\bm W_N}^r\one_{\|\bm W_N\| \ge K} \Big| \|\bm W_N\|\le \varepsilon \sqrt{N}\right)\\ &\le & \frac{\E_{\beta_N,h_N,N}\left(\nrm{\bm W_N}^r\one_{K\le \|\bm W_N\| \le \varepsilon\sqrt{N}}\right)}{\p_{\beta_N,h_N,N}(\|\bm W_N\| \le K)}\\&\le& (1+o_N(1)) \frac{\sum_{\bm v \in \cP_{q,N} : K \le \|\bm w(\bm v)\| \le \varepsilon\sqrt{N}} ~\nrm{\bm w(\bm v)}^re^{-\alpha \|\bm w(\bm v)\|^2}}{\sum_{\bm v \in \cP_{q,N} : \|\bm w(\bm v)\| \le K}  ~e^{\ip{\bar{\beta}p\bm m_*^{p-1}+\bar{h}\bm e_1}{\bm w}+\frac{1}{2} \bm Q_{\bm m_*,\beta}(\bm w (\bm v))}}\\ &\le & (1+o_N(1)) \frac{\int_{\cH_q \setminus B(0,K)} \nrm{\bm w(\bm v)}^r e^{-\alpha \|\bm w\|^2}~d\lambda(\bm w)}{\int_{\cH_q\bigcap B(0,K)} e^{\ip{\bar{\beta}p\bm m_*^{p-1}+\bar{h}\bm e_1}{\bm w}+\frac{1}{2} \bm Q_{\bm m_*,\beta}(\bm w (\bm v))}~d\lambda(\bm w) }
    \end{eqnarray*}
    for every $r\ge 1$, where the last step follows from Riemann Approximation of a sum. Hence,
    $$\limsup_{N\rightarrow \infty} \E_{\beta_N,h_N,N}\left(\nrm{\bm W_N}^r\one_{\|\bm W_N\| \ge K} \Big| \|\bm W_N\|\le \varepsilon \sqrt{N}\right) \le \frac{\int_{\cH_q \setminus B(\bm 0,K)} \nrm{\bm w(\bm v)}^r e^{-\alpha \|\bm w\|^2}~d\lambda(\bm w)}{\int_{\cH_q\bigcap B(\bm 0,K)} e^{\ip{\bar{\beta}p\bm m_*^{p-1}+\bar{h}\bm e_1}{\bm w}+\frac{1}{2} \bm Q_{\bm m_*,\beta}(\bm w (\bm v))}~d\lambda(\bm w) }~. $$
    Since $\bm Q_{\bm m_*,\beta}$ is negative definite, the above ratio goes to $0$ as $K \rightarrow \infty$, which gives \eqref{ntsh}. We thus conclude that $\bm W_N$ converges in moments to the density on $\cH_q$ with density (with respect to the Lebesgue measure on $\cH_q$) proportional to 
    $$\bm w \mapsto e^{\ip{\bar{\beta}p\bm m_*^{p-1}+\bar{h}\bm e_1}{\bm w}+\frac{1}{2}\bm Q_{\bm m_*,\beta}(\bm w)}~. $$

    Now, note that for $\bm w \in \cH_q$, from \eqref{q_form} we have $\bm Q_{\bm x_s,\beta}(\bm w) =-\tilde{\bm w}^\top \Xi_{\beta, s} \tilde{\bm w}$ where
    $$
        \Xi_{\beta, s}=-k^\dprime\left(\frac{1-s}{q}\right)I_{q-1}-k^\dprime\left(\frac{1+(q-1)s}{q}\right)J_{q-1}
    $$
    and $\tilde{\bm w} := (w_2,\ldots, w_q)^\top$.
    The covariance matrix of $\tilde{\bm W_N}$ is thus given by $-\Xi_{\beta,s}^{-1}$,
    $$
        \Xi_{\beta,s}^{-1}=\left(-k^\dprime\left(\frac{1-s}{q}\right)\right)^{-1} \times\left(I_{q-1}-\frac{k^\dprime\left(\frac{1+(q-1)s}{q}\right)}{(q-1)k^\dprime\left(\frac{1+(q-1)s}{q}\right)+k^\dprime\left(\frac{1-s}{q}\right)} J_{q-1}\right),
    $$ 
    where $J_{q-1}$ being the $(q-1)\times (q-1)$ matrix with all entries equal to $1$.
    Using the constraint $W_1 = 1-\sum_{s=2}^q W_s$, we also obtain that:

    \[
        \mathrm{Cov}(W_1, W_r) =
        \begin{cases}
            \frac{1}{(q-1)k^\dprime\left(\frac{1+(q-1)s}{q}\right)+k^\dprime\left(\frac{1-s}{q}\right)}   & \quad\text{if} ~r\ge 2, \\
            \frac{1-q}{(q-1)k^\dprime\left(\frac{1+(q-1)s}{q}\right)+k^\dprime\left(\frac{1-s}{q}\right)} & \quad\text{if} ~r=1~.   \\
        \end{cases}
    \]
    It thus follows that the asymptotic distribution of $W$ is $\cN_q(\bm 0,\Sigma)$, where
    $$
        \Sigma=\left(-\frac{q^2}{q-1}f^\dprime_{\beta,h}(s)\right)^{-1}\left(\begin{array}{cccc}q-1 & -1 & \cdots & -1 \\ -1 & 1+(q-2) \frac{k^\dprime\left(\frac{1+(q-1)s}{q}\right)}{k^\dprime\left(\frac{1-s}{q}\right)} & & -\frac{k^\dprime\left(\frac{1+(q-1)s}{q}\right)}{k^\dprime\left(\frac{1-s}{q}\right)} \\ \vdots & & \ddots & \\ -1 & -\frac{k^\dprime\left(\frac{1+(q-1)s}{q}\right)}{k^\dprime\left(\frac{1-s}{q}\right)} & & 1+(q-2) \frac{k^\dprime\left(\frac{1+(q-1)s}{q}\right)}{k^\dprime\left(\frac{1-s}{q}\right)}\end{array}\right).
    $$
    The proof of Theorem \ref{cltreg} is now complete.

\subsection{Proof of Theorem \ref{pcr}}\label{proofprpcr}
The tightness of $\bm W_{N,i} := \sqrt{N}(\bar{\bm X}_N - \bm m_i)$ conditioned on $\epm \in B(\bm m_i, \varepsilon)$ and the convergence of the law of $\bm W_{N,i}$ on bounded sets (as in \eqref{cont_approx}) imply that for any $\varepsilon>0$ smaller than the distance between any two maximizers of $H_{\beta, h}$,
    $$
        \frac{\P_{\beta, h, N}\left(\epm \in B(\bm m_i, \varepsilon)\right)}{\P_{\beta, h, N}\left(\epm \in B\left(\bm m_j, \varepsilon\right)\right)}=(1+o_N(1))\frac{A(\bm m_i) e^{N H_{\beta, h}\left(\bm m_i\right)} \int_{\cH_q} e^{\frac{1}{2} Q_{\bm m_i, \beta}(\bm w)} dw_1dw_2\ldots dw_q}{A\left(\bm m_j\right) e^{N H_{\beta, h}\left(\bm m_j\right)} \int_{\cH_q} e^{\frac{1}{2} Q_{\bm m_j, \beta}(\bm w)} dw_1dw_2\ldots dw_q}.
    $$
   Let
    $$
        C_{\beta, h, N}(\bm x)=A(\bm x) e^{N H_{\beta, h}\left(\bm x\right)} \int_{\cH_q} e^{\frac{1}{2} \bm Q_{\bm x, \beta}(\bm w)} dw_1dw_2\ldots dw_q.
    $$
    Since, $\bm m_i$ and $\bm m_j$ are maximizers of $H_{\beta,h}$, by Lemma \ref{propep1} they are either equal to $\bm x_{s_i}$ and $\bm x_{s_j}$ or one of their permutations. Assume, $\bm m_i=\bm x_{s_i}$.
    Now, note that for $\bm w \in \cH_q$, from \eqref{q_form} we have $\bm Q_{\bm m_i,\beta}(\bm w)=-\tilde{\bm w}^\top \Xi_{\beta, s_i} \tilde{\bm w}$ where
    $$
        \Xi_{\beta,s_i}=-k^\dprime\left(\frac{1-s_i}{q}\right)I_{q-1}-k^\dprime\left(\frac{1+(q-1)s_i}{q}\right)J_{q-1}
    $$
    and $\tilde{\bm w} := (w_2,\ldots, w_q)^\top$. Therefore,
    $$
        \begin{aligned}
             & \int_{\cH_q} e^{\frac{1}{2} Q_{\bm m_i, \beta}(\bm w)} dw_1dw_2\ldots dw_q= \int_{R^{q-1}} e^{-\frac{1}{2}\tilde{\bm w}^\top \Xi_{\beta,s_i} \tilde{\bm w}} dw_1dw_2\ldots dw_q \\
             & =\sqrt{2 \pi}^{q-1} \sqrt{\operatorname{det} (\Xi_{\beta,s_i}^{-1})}                                                                                                           \\
             & =\sqrt{2 \pi}^{q-1}\frac{\sqrt{q-1}}{q}\sqrt{-f^\dprime_{\beta,h}(s_i)^{-1}\left(-k^\dprime\left(\frac{1-s_i}{q}\right)\right)^{2-q}}
        \end{aligned}
    $$
    If we multiply the above expression with the prefactor $A(\bm m_i)$, we obtain:
    $$
        A(\bm m_i) \int_{\cH_q} e^{-\frac{1}{2} \bm Q_{\bm m_i, \beta}(\bm w)} dw_1dw_2\ldots dw_q =\frac{\sqrt{q-1}}{q}\tau(\bm m_i).
    $$
   Therefore, assuming that $m_i=\bm x_{s_i}$ and $m_j=\bm x_{s_j}$, we have:
    \begin{equation*}
        \lim_{N\rightarrow \infty} \frac{\P_{\beta, h, N}\left(\epm \in B(\bm m_i, \varepsilon)\right)}{\P_{\beta, h, N}\left(\epm \in B\left(\bm m_j, \varepsilon\right)\right)}=\frac{\tau(\bm m_i)}{\tau(\bm m_j)}.
    \end{equation*}

Now, note that if $\bm m_i$ and $\bm m_j$ are some permutations of $\bm x_{s_i}$ and $\bm x_{s_j}$, with at least one of these permutations not being identity, then by Proposition \ref{propep1}, one must have $h=0$, and in this case, $\p_{\beta,h,N}(\epm \in B(\bm m_k,\varepsilon)) = \p_{\beta,h,N}(\epm \in B(\bm x_{s_k},\varepsilon))$, since the measure $\p_{\beta,h,N}$ is permutation invariant if $h=0$.  Since, $\sum_{r=1}^{q}\P_{\beta, h, N}\left(\epm \in B(\bm m_r, \varepsilon)\right)=1$. Therefore,
    $$
        \lim _{N \rightarrow \infty} \P_{\beta, h, N}\left(\epm \in B(\bm m_i, \varepsilon)\right)=\frac{\tau(\bm m_i)}{\sum_{r=1}^{k}\tau(\bm m_r)}.
    $$
    The proof of the CLT part of Theorem \ref{pcr} is exactly similar to the proof of Theorem \ref{cltreg}, so we skip it. One has to only keep in mind that the variables which should be tight here, are $\bm W_{N,i}$ conditioned on the event $\epm \in B(\bm m_i,\varepsilon)$ under $\p_{\beta_N,h_N,N}$. 
    If $m_i=x_s$ for some $s$, then the exact proof \ref{prcltreg} follows. Whereas if $\bm m_i = P \bm x_s$ for some permutation matrix $P$, then the matrix with respect to $\bm Q_{\bm m_*,\beta}$ has permuted rows and columns, and hence, the covariance matrix is $P\Sigma P^T$. Therefore, the mean of the distribution is $P\Sigma P^T(\bar{\beta}p\bm m_*^{p-1}+\bar{h}\bm e_1)$. This completes the proof of Theorem \ref{pcr}.

\subsection{Proof of Theorem \ref{psp1}}\label{proofpsp1}
We now prove Theorem \ref{psp1}. For every $\bm x \in \mathcal{P}_{q,N}$, there exist unique $t(\bm x)$ and $v(\bm x) \in \mathcal{H}_q \cap \mathrm{Span}(\bm u)^\perp$, such that $\bm x = \bm m_* + N^{-1/4} t(\bm x) \bm u + N^{-1/2} v(\bm x)$.
Setting $\bm W_N := T_N\bm u + \bm V_N$ and $w(\bm x) := t(\bm x)\bm u + v(\bm x)$, and by essentially following the first few arguments in the proof of Theorem \ref{cltreg}, we get that,
\begin{equation*}
    \begin{aligned}
        &q^N Z_N(\beta_N,h_N) \e_{\beta_N,h_N,N}\left[g(\bm W_N) \mathbbm{1}_{\|\bm W_N\| \le M}\right]\\
        &= (1+o_N(1)) N^{-\frac{q-1}{2}}  A(\bm m_*) \sum_{\bm x \in \cP_{q,N}} e^{N H_{\beta_N,h_N}\left(\bm m_* +N^{-1 / 4} t(\bm x) \bm u+N^{-1 / 2} \bm v(\bm x)\right)}g(\bm w(\bm x)) \mathbbm{1}_{\|\bm w(\bm x)\| \le M}
    \end{aligned}
\end{equation*}
Now using \eqref{taylor_with_o} and Lemma \ref{l2nm}, we get that,
\begin{eqnarray*}
    &&NH_{\beta_N, h_N}\left(\bm m_*+N^{-1 / 4} t(\bm x) \bm u+N^{-1 / 2} \bm v(\bm x)\right)\\ &=&   NH_{\beta_N, h_N}(\bm m_*)+(\bar{\beta}p\ip{\bm m_*^{p-1}}{\bm u}+\bar{h}(1-q))t(\bm x)
    + \frac{1}{2} k^{\prime \prime}\left(\frac{1-s}{q}\right)\nrm{\bm v(\bm x)}^2+\frac{1}{24} q^{4} \f^{(4)}(s) t(\bm x)^4\\ &+& o_N(1).
\end{eqnarray*}
We also conclude that $T_N\bm u+ \bm V_N$ converges in law to a density, which at the point $t\bm u + \bm v$ (where $\bm u$ and $\bm v$ are orthogonal), is proportional to 
\begin{equation*}
    \exp\left(\frac{1}{2}k^\dprime\left(\frac{1-s}{q}\right)\nrm{\bm v}^2 +\frac{t^4}{24}q^4 f^{(4)}_{\beta,h}(s) +  (\bar{\beta}p\ip{\bm m_*^{p-1}}{\bm u}+\bar{h}(1-q))t\right)
\end{equation*}
 in $\cH_q \cap B(0,M)$ for any $M>0$. Similarly, we also prove that $T_N\bm u+ \bm V_N$ is tight, by \eqref{tight_tay1}. Therefore, $T_N \bm u +\bm V_N$ indeed converges in law to some random vector $T\bm u+ \bm V$. Since, $\bm u$ and $\bm V$ are orthogonal, there exists a one-to-one transformation $T\bm u + \bm V \mapsto (T,\bm V)$, and since the density factorizes into the $t$ and $\bm v$ terms, we conclude that $T$ and $\bm V$ are independent.
    \par
    Now, note that for $\bm v \in \cH_q\cap u^\perp$, using the fact that $v_1 = v_2+\ldots+v_q = 0$, we have: $$\|\bm v\|^2 = -\tilde{\bm v}^\top \Xi_{\beta,s}\tilde{\bm v}$$ where
    $$
        \Xi_{\beta,\bm s}=-\left(k^\dprime\left(\frac{1-s}{q}\right)\right)\left(I_{q-2} + J_{q-2}\right)
    $$
    and $\tilde{\bm v} := (v_3,\ldots, v_q)^\top$.
    Hence, $\tilde{\bm V}$ is Gaussian, with covariance matrix: 
    $$
        -\left(k^\dprime\left(\frac{1-s}{q}\right)\right)^{-1} \times\left(I_{q-2}-\frac{1}{q-1} J_{q-2}\right).
    $$
    Since, $V_1=0$ and hence $\operatorname{Cov}(V_1,V_r)=0$ for all $r=1,\ldots,q$. Moreover, $V_2=-\sum_{r=3}^{q}V_r$.
    Hence,
    \[
        \operatorname{Cov}(V_2, V_r) =
        \begin{cases}
            \frac{1}{(q-1)k^\dprime\left(\frac{1-s}{q}\right)}    & \quad\text{if} ~r\ge 3, \\
            -\frac{q-2}{(q-1)k^\dprime\left(\frac{1-s}{q}\right)} & \quad\text{if} ~r=2~.   \\
        \end{cases}
    \]
The proof of Theorem \ref{psp1} is now complete.

\subsection{Proof of Theorem \ref{psp2}} \label{prooftype2spc}   
Let $\bm W_N:=N^{\frac{1}{6}}(\epm-m_*)$. Fix $M>0$ and let $g:\R^q\to \R$ be a bounded continuous function. Then from similar arguments as \eqref{stp1}, we get that: 
\begin{equation*}
\begin{aligned}
    &q^N Z_N(\beta_N,h_N) \e_{\beta_N,h_N,N}\left[g(\bm W_N) \mathbbm{1}_{\|\bm W_N\| \le M}\right]\\
&=(1+o_N(1)) N^{-\frac{q-1}{2}}  A(\bm m_*) \sum_{\bm v \in \cP_{q,N}} e^{N H_{\beta_N,h_N}\left(\bm m_* + N^{-\frac{1}{6}} \bm w(\bm v)\right)}g(\bm w(\bm v)) \mathbbm{1}_{\|\bm w(\bm v)\| \le M}
\end{aligned}
\end{equation*}
where $\bm w(\bm v) := N^{1/6} (\bm v - \bm m_*)$. By Lemma \ref{specialtype2thm}, we must have $p=4, q=2$, and hence, $\mathcal{H}_q = \mathrm{Span}(\{\bm u\})$. Hence, we can write $\bm w=t(\bm v)\bm u$ and $\bm W_N=T\bm u$. Now, from \eqref{taylorspec_with_o} and Lemma \ref{l2nm}, we get,
\begin{equation*}
    NH_{\beta_{N}, h_{N}}\left(\bm m+N^{-1 / 6} t(\bm v) \bm u\right)=NH_{\beta_N, h_N}\left(\bm m\right)-\frac{32}{15}t(\bm v)^6-\bar{h}t(\bm v)+o_N(1).
\end{equation*}
Hence, under $\p_{\beta_N,h_N,N}$, $T$ conditioned on $|T| \le M$ converges weakly to the density on $\cH_q \bigcap B(0,M)$ with density (with respect to the Lebesgue measure on $\cH_q$) proportional to 
$$t\mapsto \exp\left(-\frac{32}{15}t^6 -\bar{h}t\right).$$
The tightness of $T$ follows from \eqref{tight_tay2}. This completes the proof of Lemma \ref{psp2}.

We now prove a lemma that is necessary for proving asymptotics of the ML estimate of $\beta$.

\begin{lem}\label{cruciallemma8}
    (Asymptotic distribution of $\nrm{\epm}_p^p$ under perturbed $\beta$). Fix $(\beta, h) \in \Theta$, and $\bar{\beta}, \bar{h} \in \R$. Then the following hold:
    \begin{enumerate}[i.]
        \item Suppose $(\beta, h)$ is regular and denote the unique maximizer of $H$ by $\bm m_*= \bm m_*(\beta, h)$. Then, for $\bm X \sim \P_{\beta+N^{-\frac{1}{2}} \bar{\beta}, h}$, as $N \rightarrow \infty$,
              \begin{itemize}
                  \item if $\bm m_* \neq \bm x_0$,
                        \begin{equation*}
                            N^{\frac{1}{2}}\left(\nrm{\epm}_p^p - \nrm{\bm m_*}_p^p\right) \xrightarrow{D} \cN\left(-\frac{\bar{\beta}p^2(q-1)^2}{q^2f^\dprime_{\beta,h}(s)}\left(m_1^{p-1}-m_2^{p-1}\right)^2,-\frac{p^2(q-1)^2}{q^2f^\dprime_{\beta,h}(s)}\left(m_1^{p-1}-m_2^{p-1}\right)^2\right),
                        \end{equation*}
                        where $\bm m_*=(m_1,m_2,\ldots,m_q)$
                  \item if $\bm m_* = \bm x_0$ then $N \left(\nrm{\epm}_p^p - \nrm{\bm m_*}_p^p\right)$ converges to a generalised chi-squared distribution. More specifically,
                        \begin{equation}\label{bp10e}
                            N \left(\nrm{\epm}_p^p - \nrm{\bm m_*}_p^p\right) \xrightarrow{D} \frac{p(p-1)}{2q^{p-2}}\bm W^\top \bm W,
                        \end{equation}
                        where $\bm W \sim \cN_q\left(\boldsymbol{0}, \Sigma\right)$.
              \end{itemize}
        \item Suppose $(\beta, h)$ is critical and denote the $K$ maximizers of $H$ denoted by $\bm m_1:=$ $\bm m_1(\beta, h, p)$, $\ldots$, $\bm m_{K}:=\bm m_{K}(\beta, h, p)$. Then, for $\epm \sim \P_{\beta, h, p}$, as $N \rightarrow \infty$,
              \begin{equation}\label{bel11}
                  \nrm{\epm}_p^p \xrightarrow{D} \sum_{k=1}^K p_k \delta_{\nrm{m_{K}}_p^p},
              \end{equation}
              where
              \begin{equation*}
                  p_k:=\frac{\tau(\bm m_k)}{\sum_{i=1}^{K}\tau(\bm m_i)},
              \end{equation*}
              and $\tau$ is as defined in \eqref{tau_def}.
              Moreover, if $\bm m_*$ is any local maximizer of $H$ contained in the interior of a set $A \subseteq \cP_q $, such that $H(\bm m_*)>H(x)$ for all $x \in$ $A \setminus \{\bm m\}$, then for $\epm \sim \P_{\beta+N^{-\frac{1}{2}} \bar{\beta}, h}$, as $N \rightarrow \infty$,
              \begin{itemize}
                  \item if $\bm x_0 \notin A$,
                        \begin{eqnarray}
                        \label{pertnorm_crit1}
                           && N^{\frac{1}{2}}\left(\nrm{\epm}_p^p - \nrm{\bm m_*}_p^p\right) \mid \{\epm \in A\}\\ \nonumber &\xrightarrow{D}& \cN\left(-\frac{\bar{\beta}p^2(q-1)^2}{q^2f^\dprime_{\beta,h}(s)}\left(m_{(q)}^{p-1}-m_{(1)}^{p-1}\right)^2,-\frac{p^2(q-1)^2}{q^2f^\dprime_{\beta,h}(s)}\left(m_{(q)}^{p-1}-m_{(1)}^{p-1}\right)^2\right),
                        \end{eqnarray}

                        where $m_{(q)}$ and $m_{(1)}$ denote the largest and smallest elements of $\bm m_*$, respectively.
                        
                  \item if $\bm x_0 \in A$ then $N \left(\nrm{\epm}_p^p - \nrm{\bm m_*}_p^p\right)$ converges to a generalised chi-squared distribution. More specifically,
                        \begin{equation*}
                            N \left(\nrm{\epm}_p^p - \nrm{\bm m_*}_p^p\right) \mid \{\epm \in A\} \xrightarrow{D} \frac{p(p-1)}{2q^{p-2}}\bm W^\top \bm W,
                        \end{equation*}
                        where $\bm W \sim \cN_q\left(0, \Sigma\right)$.
              \end{itemize}
        \item Suppose $(\beta, h)$ is type I special and denote the unique maximizer of $H$ by $\bm m_*= \bm m_*(\beta, h, p) := (m_1,\ldots,m_q)$. Let $\bm u=(1-q,1,\ldots,1)$. Then, for $\epm \sim \P_{\beta+N^{-\frac{3}{4}} \bar{\beta}, h}$, as $N \rightarrow \infty$,
              \begin{itemize}
                  \item if $(p,q)\notin \{(2,2)\cup (3,2)\}$,
                        \begin{equation*}
                            N^{\frac{1}{4}} \left(\nrm{\epm}_p^p - \nrm{\bm m_*}_p^p\right) \xrightarrow{D} -T_{\bar{\beta},0}p(q-1)\left(m_1^{p-1}-m_2^{p-1}\right), 
                        \end{equation*}
                  \item if $(p,q)= (2,2)~\text{or}~ (3,2)$ then,
                        \begin{equation*}
                            N^{\frac{1}{2}}\left(\nrm{\epm}_p^p - \nrm{\bm m_*}_p^p\right) \xrightarrow{D} \frac{p(p-1)}{2^{p-2}} T_{0,0}^2.
                        \end{equation*}
                        Here, $\frac{p(p-1)}{2^{p-2}}T_{0,0}^2$ has a density proportional to
$$t^{-1/2} \exp\left(\frac{2^{2p-3} t^2}{3 p^2(p-1)^2} q^4 f_{\beta,h}^{(4)}(0)\right).$$
              \end{itemize}
        \item Suppose $(\beta, h)$ is type II special and denote the unique maximizer of $H$ by $\bm m_*= \bm m_*(\beta, h, p)$. Let $\bm u=(1-q,1,\ldots,1)$. Then, for $\epm \sim \P_{\beta+N^{-\frac{5}{6}} \bar{\beta}, h}$, as $N \rightarrow \infty$,
              \begin{equation*}
                  N^{\frac{1}{3}}\left(\nrm{\epm}_p^p - \nrm{\bm m_*}_p^p\right) \xrightarrow{D} 3F_{0}^2.
              \end{equation*}
               Here, $3F_{0}^2$ has a density proportional to
              $$
                 t^{-1/2} \exp\left(-\frac{32}{405}t^3\right).
              $$
    \end{enumerate}
\end{lem}
\begin{proof}
    \begin{enumerate}[i.]
        \item Let $\phi(t) :=\nrm{\bm m_*+t(\epm-\bm m_*)}_p^p$. By the mean value theorem, there exists $\alpha \in [0,1]$ such that $\phi(1)=\phi(0)+\phi'(\alpha)$ and hence,
        
              \begin{equation}\label{lim_pnorm}
                      \nrm{\epm}_p^p - \nrm{\bm m_*}_p^p = \sum_{r=1}^{q}p(\bar{X}_{\cdot r}-m_r)[m_r+\alpha(\bar{X}_{\cdot r}-m_r)]^{p-1}.
              \end{equation}
              Now, $\bar{X}_{\cdot r}\xrightarrow{P} m_r$. Hence, it follows from Theorem \ref{cltreg} that:
              $$
                  \sqrt N \left(\nrm{\epm}_p^p - \nrm{\bm m_*}_p^p\right) \xrightarrow{D} \ip{\bm W}{p\bm m_*^{p-1}},
              $$
              where $\bm W \sim \cN_q\left(\bar{\beta}p\Sigma \bm m_*^{p-1} , \Sigma\right)$, where $\Sigma$ is as defined in \eqref{sigma_def}.
              Moreover,
              \begin{equation*}
                  \E_{\beta_N,h,N}\left(\sqrt N \left(\nrm{\epm}_p^p - \nrm{\bm m_*}_p^p\right)\right) \to \E\left(\ip{\bm W}{p\bm m_*^{p-1}}\right).
              \end{equation*}
            It is easy to check that,
              $$
                  (p\bm m_*^{p-1})^\top \Sigma (p\bm m_*^{p-1})=-\frac{p^2(q-1)^2}{q^2f^\dprime_{\beta,h}(s)}\left(m_1^{p-1}-m_2^{p-1}\right)^2.
              $$
              Therefore,
              $$
                  \sqrt N \left(\nrm{\epm}_p^p - \nrm{\bm m_*}_p^p\right) \xrightarrow{D} \cN\left(-\frac{\bar{\beta}p^2(q-1)^2}{q^2f^\dprime_{\beta,h}(s)}\left(m_1^{p-1}-m_2^{p-1}\right)^2,-\frac{p^2(q-1)^2}{q^2f^\dprime_{\beta,h}(s)}\left(m_1^{p-1}-m_2^{p-1}\right)^2\right).
              $$
              On the other hand if $\bm m_*=\bm x_0$, then by Taylor's theorem with Lagrange Reminder, $\phi(1)                             =\phi(0)+\phi^\prime(0)+\frac{1}{2}\phi^\dprime(\alpha)$ for some $\alpha \in [0,1]$, and hence,
              \begin{equation}\label{lim_pnorm2}               
                      \nrm{\epm}_p^p - \nrm{\bm m_*}_p^p  = \frac{1}{2}\sum_{r=1}^{q}p(p-1)(\bar{X}_{\cdot r}-m_r)^2[m_r+\alpha(\bar{X}_{\cdot r}-m_r)]^{p-2}
              \end{equation}
              \eqref{bp10e} now follows from \eqref{lim_pnorm2} and Theorem \eqref{cltreg}.
        \item First, note that \eqref{bel11} follows from \eqref{crit_wlln}. Next, if $\bm m=\bm x_0 \notin A$ , then from \eqref{lim_pnorm} we get the same limiting distribution, i. e. ,
              \begin{eqnarray*}
                  &&\sqrt N \left(\nrm{\epm}_p^p - \nrm{\bm m}_p^p\right) \mid \{\epm \in A\}\\ &\xrightarrow{D}& \cN\left(-\frac{\bar{\beta}p^2(q-1)^2}{q^2f^\dprime_{\beta,h}(s)}\left(m_{(q)}^{p-1}-m_{(1)}^{p-1}\right)^2,-\frac{p^2(q-1)^2}{q^2f^\dprime_{\beta,h}(s)}\left(m_{(q)}^{p-1}-m_{(1)}^{p-1}\right)^2\right).
              \end{eqnarray*}
                  
              If $\bm m=\bm x_0 \in A$ , then from \eqref{lim_pnorm2} we get the same limiting distribution, i. e. ,
              $$
                  N \left(\nrm{\epm}_p^p - \nrm{\bm m_*}_p^p\right) \mid \{\epm \in A\} \xrightarrow{D} \frac{p(p-1)}{2q^{p-2}}\bm W^\top \bm W
              $$
              
        \item If $(p,q)\neq (2,2)~\text{or}~(3,2)$ then, $\bm m_* \neq \bm x_0$ by Lemma \ref{strange_spec}. Hence, from \eqref{lim_pnorm} and Theorem \ref{psp1}, we have:
              $$
                  \begin{aligned}
                      N^{\frac{1}{4}} \left(\nrm{\epm}_p^p - \nrm{\bm m_*}_p^p\right) \xrightarrow{D} & T_{\bar{\beta},0}\ip{\bm u}{p\bm m_*^{p-1}}                \\
                                                                                                      & =-T_{\bar{\beta},0} p(q-1)\left(m_1^{p-1}-m_2^{p-1}\right).
                  \end{aligned}
              $$
              On the other hand if $(p,q)= (2,2)~\text{or}~(3,2)$ then, $\bm m_* = \bm x_0$ by Lemma \ref{strange_spec}. It thus follows from  \eqref{lim_pnorm2} and Theorem \ref{psp1}, that:
              $$
                  \begin{aligned}
                      N^{\frac{1}{2}}\left(\nrm{\epm}_p^p - \nrm{\bm m_*}_p^p\right) \xrightarrow{D} & \frac{p(p-1)}{2q^{p-2}} T_{\bar{\beta},0}^2\sum_{r=1}^{q}u_r^2 \\
                                                                                                     & =\frac{p(p-1)}{2^{p-2}} T_{0,0}^ 2.
                  \end{aligned}
              $$
        \item From \eqref{lim_pnorm2} and using the fact that $p=4$ and $q=2$ (by Lemma \ref{specialtype2thm}), we get from Theorem \ref{psp2}
              $$
                  \begin{aligned}
                      N^{\frac{1}{3}}\left(\nrm{\epm}_p^p - \nrm{\bm m_*}_p^p\right) \xrightarrow{D} & \frac{3}{2}F_{h}^2\sum_{r=1}^{q}u_r^2\\
                                                                                                     & =3F_{0}^2.
                  \end{aligned}
              $$
    \end{enumerate}
\end{proof}

\section{Proofs of the Asymptotics of the ML Estimates}
In this section, we prove the results on the asymptotics of the ML estimates of $\beta$ and $h$ stated in Section \ref{sec:asmml}. 

\subsection{Proof of Theorem \ref{hmlreg}}\label{prhmlreg}
    We now prove Theorem \ref{hmlreg}. For any $t \in \R$, we have by \eqref{h_like}, Lemma \ref{like_monotone}, and Theorem \ref{cltreg}, together with the uniform integrability of $\|\bm W\|$,
    $$
        \begin{aligned}
            \P_{\beta, h, p}\left(N^{\frac{1}{2}}\left(\hat{h}_N-h\right) \leq t\right) & =\P_{\beta, h, p}\left(\hat{h}_N \leq h+\frac{t}{N^{\frac{1}{2}}}\right)                                                                                                    \\
                                                                                        & =\P_{\beta, h, p}\left(u_{N, 1}\left(\beta, \hat{h}_N, p\right) \leq u_{N, 1}\left(\beta, h+\frac{t}{N^{\frac{1}{2}}}, p\right)\right)                                      \\
                                                                                        & =\P_{\beta, h, p}\left(\bxd \leq \E_{\beta, h+N^{-\frac{1}{2}} t, p}\left(\bxd\right)\right)                                                                          \\
                                                                                        & =\P_{\beta, h, p}\left(N^{\frac{1}{2}}\left(\bxd-m_1\right) \leq \E_{\beta, h+N^{-\frac{1}{2}} t, p}\left(N^{\frac{1}{2}}\left(\bxd-m_1\right)\right)\right)          \\
                                                                                        & \rightarrow \P_{\beta, h, p}\left(\cN\left(0,-\frac{(q-1)^2}{q^2f_{\beta,h}^\dprime\left(s\right)}\right) \leq-\frac{t(q-1)^2}{q^2f_{\beta,h}^\dprime\left(s\right)}\right) \\
                                                                                        & =\P_{\beta, h, p}\left(\cN\left(0,-\frac{q^2f_{\beta,h}^\dprime\left(s\right)}{(q-1)^2}\right) \leq t\right).
        \end{aligned}
    $$
    This completes the proof of Theorem \ref{hmlreg}.

\subsection{Proof of Theorem \ref{hml77}}\label{hml77pr}
   Let $\bm m_1,\ldots, \bm m_K$ be the $K$ maximizers of $f_{\beta,h}$, ordered in ascending order of their first coordinates. Let us start with disjoint sets $\{A_r\}_{1\le r\le K}$ uniting to $\mathcal{P}_q$, such that $A_i$ contains $\bm m_i$ in its interior, for all $1\le i\le K$.
    Fixing $t \in \mathbb{R}$, we have the following for every $k \in [K]$:
    $$
        \begin{aligned}
            \P_{\beta, h, p}\left(N^{\frac{1}{2}}\left(\hat{h}_N-h\right) \leq t\right) & =\P_{\beta, h, p}\left(N^{\frac{1}{2}}\left(\bar{X}_{\cdot 1} -m_{k,1}\right) \leq \E_{\beta, h+N^{-\frac{1}{2}} t, p}\left(N^{\frac{1}{2}}\left(\bar{X}_{\cdot 1} -m_{k,1}\right)\right)\right) \\
                                                                                        & =: \sum_{i=1}^{K}T_i^k,
        \end{aligned}
    $$
    where
    $$
        T_i^k =\P_{\beta, h, p}\left(N^{\frac{1}{2}}\left(\bar{X}_{\cdot 1}-m_{k,1}\right) \leq \E_{\beta, h+N^{-\frac{1}{2}} t, p}\left(N^{\frac{1}{2}}\left(\bar{X}_{\cdot 1} -m_{k,1}\right)\right) \Big| \epm \in A_i\right) \P_{\beta, h, p}\left(\epm \in A_i\right)
    $$
    Now, by the law of iterated expectations, we have for large $N$,
    \begin{equation*}
        \E_{\beta, h+N^{-\frac{1}{2}} t, p}\left(N^{\frac{1}{2}}\left(\bar{X}_{\cdot 1} -m_{k,1}\right)\right)=\sum_{i=1}^{K}S_i^k,
    \end{equation*}
    where
    \begin{equation*}
        S_i^k:=\E_{\beta, h+N^{-\frac{1}{2}} t, p}\left(N^{\frac{1}{2}}\left(\bar{X}_{\cdot 1} -m_{k,1}\right) \Big| \epm \in A_i\right) \P_{\beta, h+N^{-\frac{1}{2}} t, p}\left(\epm \in A_i\right)
    \end{equation*}

    \noindent (1)~ Suppose that $(\beta,h)\in \cC_{p,q}^1\setminus \{(\beta_c,0)\}$. Then, by Lemma \ref{crit_part} i. and Proposition \ref{propep1}, $H_{\beta,h}$ has exactly two global maximizers $\bm m_1 = \bm x_{s_1}$ and $\bm m_2 = \bm x_{s_2}$ for some $s_2>s_1>0$. By Theorem \ref{pcr},
    $$
        \E_{\beta, h+N^{-\frac{1}{2}} t, p}\left(N^{\frac{1}{2}}\left(\bar{X}_{\cdot 1}-m_{k,1}\right) \mid \epm \in A_k\right) \rightarrow-\frac{t(q-1)^2}{q^2f_{\beta,h}^\dprime\left(s_k\right)},
    $$
   for $k\in \{1,2\}$, as $N \rightarrow \infty$. Suppose that $t>0$. Then, by Lemma \ref{pcl_h}, we know that $$\p_{\beta,h+N^{-\frac{1}{2}} t, p} (\epm \in A_1) \le C_1e^{-C_2\sqrt{N}}$$ for some constants $C_1,C_2>0$. Hence, we have:
   $$\E_{\beta, h+N^{-\frac{1}{2}} t, p}\left(N^{\frac{1}{2}}\left(\bar{X}_{\cdot 1} -m_{2,1}\right)\right) \rightarrow -\frac{t(q-1)^2}{q^2f_{\beta,h}^\dprime\left(s_2\right)}.$$
   Hence, for $t > 0$, we have:
   \begin{equation}\label{1711}
       \P_{\beta, h, p}\left(N^{\frac{1}{2}}\left(\hat{h}_N-h\right) \leq t\right) \rightarrow p_1 + p_2 \p\left(\mathcal{N}\left(0,-\frac{(q-1)^2}{q^2 f_{\beta,h}''(s_2)}\right) \le -\frac{t(q-1)^2}{q^2f_{\beta,h}''(s_2)}\right).
   \end{equation}
 Next, suppose that $t<0$. Then, by Lemma \ref{pcl_h}, we know that $$\p_{\beta,h+N^{-\frac{1}{2}} t, p} (\epm \in A_2) \le C_1e^{-C_2\sqrt{N}},$$ 
for some constant $C_1,C_2>0$. Hence, we have:
   $$\E_{\beta, h+N^{-\frac{1}{2}} t, p}\left(N^{\frac{1}{2}}\left(\bar{X}_{\cdot 1} -m_{1,1}\right)\right) \rightarrow -\frac{t(q-1)^2}{q^2f_{\beta,h}^\dprime\left(s_1\right)}.$$
   Hence, for $t < 0$, we have:
   \begin{equation}\label{1712}
       \P_{\beta, h, p}\left(N^{\frac{1}{2}}\left(\hat{h}_N-h\right) \leq t\right) \rightarrow p_1 \p\left(\mathcal{N}\left(0,-\frac{(q-1)^2}{q^2 f_{\beta,h}''(s_1)}\right) \le -\frac{t(q-1)^2}{q^2f_{\beta,h}''(s_1)}\right)
   \end{equation}
Part (1) now follows from \eqref{1711} and \eqref{1712}.
\vspace{0.15in}

\noindent (2)~Suppose that $(\beta,h) \in \cC_{p,q}^2$. Then, by Lemma \ref{crit_part} ii. (b), $H_{\beta,h}$ has exactly $q$ global maximizers, which are all the possible permutations of $\bm x_s$ for some $s>0$. By Theorem \ref{pcr}, we have:
 $$\E_{\beta, h+N^{-\frac{1}{2}} t, p}\left(N^{\frac{1}{2}}\left(\bar{X}_{\cdot 1}-m_{q,1}\right) \mid \epm \in A_q\right) \rightarrow-\frac{t(q-1)^2}{q^2f_{\beta,h}^\dprime\left(s\right)}\quad\quad\text{and}
    $$

    $$\E_{\beta, h+N^{-\frac{1}{2}} t, p}\left(N^{\frac{1}{2}}\left(\bar{X}_{\cdot 1}-m_{r,1}\right) \mid \epm \in A_r\right) \rightarrow-\frac{t(q-1)}{q^2f_{\beta,h}^\dprime\left(s\right)}\left(1+(q-2) \frac{k^\dprime\left(\frac{1+(q-1)s}{q}\right)}{k^\dprime\left(\frac{1-s}{q}\right)}\right)\quad(1\le r< q).$$
Now, suppose that $t>0$. Then, by Lemma \ref{pcl_h}, we know that for all $r\in [q-1]$, $$\p_{\beta,h+N^{-\frac{1}{2}} t, p} (\epm \in A_r) \le C_1e^{-C_2\sqrt{N}}$$ for some constant $C_1,C_2>0$. Hence, we have:
   $$\E_{\beta, h+N^{-\frac{1}{2}} t, p}\left(N^{\frac{1}{2}}\left(\bar{X}_{\cdot 1} -m_{q,1}\right)\right) \rightarrow -\frac{t(q-1)^2}{q^2f_{\beta,h}^\dprime\left(s\right)}.$$
   Hence, for $t > 0$, we have:
   \begin{equation}\label{1713}
       \P_{\beta, h, p}\left(N^{\frac{1}{2}}\left(\hat{h}_N-h\right) \leq t\right) \rightarrow 1-p_q + p_q \p\left(\mathcal{N}\left(0,-\frac{(q-1)^2}{q^2 f_{\beta,h}''(s)}\right) \le -\frac{t(q-1)^2}{q^2f_{\beta,h}''(s)}\right)
   \end{equation} 
Next, suppose that $t<0$. Then, by Lemma \ref{pcl_h}, we know that $$\p_{\beta,h+N^{-\frac{1}{2}} t, p} (\epm \in A_q) \le C_1e^{-C_2\sqrt{N}}$$ for some constant $C_1,C_2>0$. Hence, we have:
   $$\E_{\beta, h+N^{-\frac{1}{2}} t, p}\left(N^{\frac{1}{2}}\left(\bar{X}_{\cdot 1} -m_{1,1}\right)\right) \rightarrow -\frac{t(q-1)}{q^2f_{\beta,h}^\dprime\left(s\right)}\left(1+(q-2) \frac{k^\dprime\left(\frac{1+(q-1)s}{q}\right)}{k^\dprime\left(\frac{1-s}{q}\right)}\right).$$
   Hence, for $t < 0$, we have:
   \begin{equation}\label{1714}
       \P_{\beta, h, p}\left(N^{\frac{1}{2}}\left(\hat{h}_N-h\right) \leq t\right) \rightarrow (1-p_q) \p\left(\mathcal{N}\left(0,-\frac{q^2f_{\beta,h}^\dprime\left(s\right)}{q-1}\left(1+(q-2) \frac{k^\dprime\left(\frac{1+(q-1)s}{q}\right)}{k^\dprime\left(\frac{1-s}{q}\right)}\right)^{-1}\right) \le t\right)
   \end{equation}
Part (2) now follows from \eqref{1713} and \eqref{1714}.
\vspace{0.15in}

\noindent (3)~Suppose that $(\beta,h) = (\beta_c,0)$ is a critical point. Then, by Lemma \ref{crit_part} ii. (c), $H_{\beta,h}$ has exactly $q+1$ global maximizers, which are all the possible permutations of $\bm x_s$ for some $s>0$, and the vector $(\frac{1}{q},\ldots,\frac{1}{q})$. So, note that $K=q+1$ here, $\bm m_{q+1} = \bm x_s$, $\bm m_q = \bm x_0$ and $\bm m_1,\ldots,\bm m_{q-1}$ are the remaining permutations of $\bm x_s$. By Theorem \ref{pcr}, we have:
 $$\E_{\beta, h+N^{-\frac{1}{2}} t, p}\left(N^{\frac{1}{2}}\left(\bar{X}_{\cdot 1}-m_{q+1,1}\right) \mid \epm \in A_{q+1}\right) \rightarrow-\frac{t(q-1)^2}{q^2f_{\beta,h}^\dprime\left(s\right)},
    $$

 $$\E_{\beta, h+N^{-\frac{1}{2}} t, p}\left(N^{\frac{1}{2}}\left(\bar{X}_{\cdot 1}-m_{q,1}\right) \mid \epm \in A_{q}\right) \rightarrow-\frac{t(q-1)^2}{q^2f_{\beta,h}^\dprime\left(0\right)},\quad \quad \text{and}
    $$

    $$\E_{\beta, h+N^{-\frac{1}{2}} t, p}\left(N^{\frac{1}{2}}\left(\bar{X}_{\cdot 1}-m_{r,1}\right) \mid \epm \in A_r\right) \rightarrow-\frac{t(q-1)}{q^2f_{\beta,h}^\dprime\left(s\right)}\left(1+(q-2) \frac{k^\dprime\left(\frac{1+(q-1)s}{q}\right)}{k^\dprime\left(\frac{1-s}{q}\right)}\right)\quad(1\le r< q).$$
Now, suppose that $t>0$. Then, by Lemma \ref{pcl_h}, we know that for all $r\in [q]$, $$\p_{\beta,h+N^{-\frac{1}{2}} t, p} (\epm \in A_r) \le C_1e^{-C_2\sqrt{N}}$$ for some constant $C_1,C_2>0$. Hence, we have:
   $$\E_{\beta, h+N^{-\frac{1}{2}} t, p}\left(N^{\frac{1}{2}}\left(\bar{X}_{\cdot 1} -m_{q+1,1}\right)\right) \rightarrow -\frac{t(q-1)^2}{q^2f_{\beta,h}^\dprime\left(s\right)}.$$
   Hence, for $t > 0$, we have:
   \begin{equation}\label{1723}
       \P_{\beta, h, p}\left(N^{\frac{1}{2}}\left(\hat{h}_N-h\right) \leq t\right) \rightarrow 1-p_{q+1} + p_{q+1} \p\left(\mathcal{N}\left(0,-\frac{(q-1)^2}{q^2 f_{\beta,h}''(s)}\right) \le -\frac{t(q-1)^2}{q^2f_{\beta,h}''(s)}\right)
   \end{equation} 
Next, suppose that $t<0$. Then, by Lemma \ref{pcl_h}, we know that $$\p_{\beta,h+N^{-\frac{1}{2}} t, p} (\epm \in A_r) \le C_1e^{-C_2\sqrt{N}}\quad\text{for}~r\in \{q,q+1\}$$ for some constant $C_1,C_2>0$. Hence, we have:
   $$\E_{\beta, h+N^{-\frac{1}{2}} t, p}\left(N^{\frac{1}{2}}\left(\bar{X}_{\cdot 1} -m_{1,1}\right)\right) \rightarrow -\frac{t(q-1)}{q^2f_{\beta,h}^\dprime\left(s\right)}\left(1+(q-2) \frac{k^\dprime\left(\frac{1+(q-1)s}{q}\right)}{k^\dprime\left(\frac{1-s}{q}\right)}\right).$$
   Hence, for $t < 0$, we have:
   \begin{eqnarray}\label{1724}
       &&\P_{\beta, h, p}\left(N^{\frac{1}{2}}\left(\hat{h}_N-h\right) \leq t\right)\nonumber\\ &\rightarrow& (1-p_q-p_{q+1}) \p\left(\mathcal{N}\left(0,-\frac{q^2f_{\beta,h}^\dprime\left(s\right)}{q-1}\left(1+(q-2) \frac{k^\dprime\left(\frac{1+(q-1)s}{q}\right)}{k^\dprime\left(\frac{1-s}{q}\right)}\right)^{-1}\right) \le t\right)
   \end{eqnarray}
Part (3) now follows from \eqref{1723} and \eqref{1724}, and the observation that since $p_1=\ldots=p_{q-1}=p_{q+1}$, one must have $q p_{q+1} + p_q = 1$, i.e. $p_{q+1} = (1-p_q)/q$. The proof of Theorem \ref{hml77} is now complete.

\subsection{Proof Theorem \ref{betaml7}}\label{prbetaml7}
    In this section, we prove Theorem \ref{betaml7}. For any $t \in \R$, we have by \eqref{b_like}, Lemma \ref{cruciallemma8}, Lemma \ref{like_monotone}, and \eqref{pertnorm_crit1}, together with uniform integrability of all powers of $\|\bm W\|$,
    
        \begin{eqnarray*}
            &&\P_{\beta, h, p}\left(N^{\frac{1}{2}}\left(\hat{\beta}_N-\beta\right) \leq t\right) \\& =&\P_{\beta, h, p}\left(\hat{\beta}_N \leq \beta+\frac{t}{N^{\frac{1}{2}}}\right)                                                                                                                                                           \\
                                                                                                & =&\P_{\beta, h, p}\left(u_{N, p}\left(\hat{\beta}_N,h, p\right) \leq u_{N, p}\left(\beta+\frac{t}{N^{\frac{1}{2}}}, h, p\right)\right)                                                                                                      \\
                                                                                                & =&\P_{\beta, h, p}\left(\nrm{\epm}_p^p \leq \E_{\beta+N^{-\frac{1}{2}} t, h, p}\left(\nrm{\epm}_p^p\right)\right)                                                                                                                             \\
                                                                                                & =&\P_{\beta, h, p}\left(N^{\frac{1}{2}}\left(\nrm{\epm}_p^p-\nrm{\bm m}_p^p\right) \leq \E_{\beta, h+N^{-\frac{1}{2}} t, p}\left(N^{\frac{1}{2}}\left(\nrm{\epm}_p^p-\nrm{\bm m}_p^p\right)\right)\right)                                             \\
                                                                                                & \rightarrow& \P_{\beta, h, p}\left(\cN\left(0,-\frac{\bar{\beta}p^2(q-1)^2}{q^2f^\dprime_{\beta,h}(s)}\left(m_1^{p-1}-m_2^{p-1}\right)^2\right) \leq-\frac{tp^2(q-1)^2}{q^2f^\dprime_{\beta,h}(s)}\left(m_1^{p-1}-m_2^{p-1}\right)^2\right) \\
                                                                                                & =&\P_{\beta, h, p}\left(\cN\left(0,-\frac{q^2f^\dprime_{\beta,h}(s)}{p^2(q-1)^2}\left(m_1^{p-1}-m_2^{p-1}\right)^{-2}\right) \leq t\right).                
        \end{eqnarray*}
    
    This completes the proof for \eqref{mle_b_reg1}.
    \par
    Next coming to the case where $\bm m_*=x_0$. We get that,
    $$
        \begin{aligned}
            \P_{\beta, h, p}\left(N^{\frac{1}{2}}\left(\hat{\beta}_N-\beta\right) \leq t\right) & =\P_{\beta, h, p}\left(\hat{\beta}_N \leq \beta+\frac{t}{N^{\frac{1}{2}}}\right)                                                                                   \\
                                                                                                & =\P_{\beta, h, p}\left(u_{N, p}\left(\hat{\beta}_N,h, p\right) \leq u_{N, p}\left(\beta+\frac{t}{N^{\frac{1}{2}}}, h, p\right)\right)                              \\
                                                                                                & =\P_{\beta, h, p}\left(\nrm{\epm}_p^p \leq \E_{\beta+N^{-\frac{1}{2}} t, h, p}\left(\nrm{\epm}_p^p\right)\right)                                                     \\
                                                                                                & =\P_{\beta, h, p}\left(N\left(\nrm{\epm}_p^p-\nrm{\bm m}_p^p\right) \leq \E_{\beta+N^{-\frac{1}{2}}t, h, p}\left(N\left(\nrm{\epm}_p^p-\nrm{\bm m}_p^p\right)\right)\right) \\
                                                                                                & \rightarrow \P\left(\bm W^\top \bm W \leq\E\left[ \bm W^\top \bm W  \right]\right)                                                                   \\
                                                                                                & =\P\left(\bm W^\top \bm W\leq\frac{1-q}{k^\dprime\left(\frac{1}{q}\right)}\right)                                                                                \\
                                                                                                & = \gamma_1.
        \end{aligned}
    $$
The proof of Theorem \ref{betaml7} is now complete.

\subsection{Proof of Theorem \ref{btml44}}\label{btml44pr}
   Let $\bm m_1,\ldots, \bm m_K$ be the $K$ maximizers of $f_{\beta,h}$, arranged in ascending order of their $L^p$ norms. Let us start with disjoint sets $\{A_r\}_{1\le r\le K}$ uniting to $\mathcal{P}_q$, such that $A_i$ contains $\bm m_i$ in its interior, for all $1\le i\le K$.
    Fixing $t \in \mathbb{R}$, we have the following for every $k \in [K]$:

        \begin{eqnarray*}
            &&\P_{\beta, h, p}\left(N^{\frac{1}{2}}\left(\hat{\beta}_N-\beta\right) \leq t\right)\\ &=& \P_{\beta, h, p}\left(N^{\frac{1}{2}}\left(\|\epm\|_p^p - \|\bm m_{k}\|_p^p\right) \leq \E_{\beta+N^{-\frac{1}{2}} t, h, p}\left(N^{\frac{1}{2}}\left(\|\epm\|_p^p - \|\bm m_{k}\|_p^p\right)\right)\right) \\& =:& \sum_{i=1}^{K}T_i^k,
        \end{eqnarray*}

    where
    \begin{eqnarray*}
        &&T_i^k \\ &=&\P_{\beta, h, p}\left(N^{\frac{1}{2}}\left(\|\epm\|_p^p - \|\bm m_{k}\|_p^p\right) \leq \E_{\beta+N^{-\frac{1}{2}} t, h, p}\left(N^{\frac{1}{2}}\left(\|\epm\|_p^p - \|\bm m_{k}\|_p^p\right)\right) \Big| \epm \in A_i\right)\\&\times& \P_{\beta, h, p}\left(\epm \in A_i\right)
    \end{eqnarray*}

    Now, by the law of iterated expectations, we have for large $N$,
    \begin{equation*}
        \E_{\beta+N^{-\frac{1}{2}} t, h, p}\left(N^{\frac{1}{2}}\left(\|\epm\|_p^p - \|\bm m_{k}\|_p^p\right)\right)=\sum_{i=1}^{K}S_i^k,
    \end{equation*}
    where
    \begin{equation*}
        S_i^k:=\E_{\beta+N^{-\frac{1}{2}} t, h, p}\left(N^{\frac{1}{2}}\left(\|\epm\|_p^p - \|\bm m_{k}\|_p^p\right) \Big| \epm \in A_i\right) \P_{\beta, h+N^{-\frac{1}{2}} t, p}\left(\epm \in A_i\right)
    \end{equation*}

    \noindent (1)~ Suppose that $(\beta,h)\in \cC_{p,q}^1\setminus \{(\beta_c,0)\}$. Then, by Lemma \ref{crit_part} i. and Proposition \ref{propep1}, $H_{\beta,h}$ has exactly two global maximizers $\bm m_1 = \bm x_{s_1}$ and $\bm m_2 = \bm x_{s_2}$ for some $s_2>s_1>0$. It is easy to check that in this case, $\|\bm m_2\|_p > \|\bm m_1\|_p$. By Lemma \ref{cruciallemma8} ii., 
    $$
        \E_{\beta+N^{-\frac{1}{2}} t, h, p}\left(N^{\frac{1}{2}}\left(\|\epm\|_p^p - \|\bm m_{k}\|_p^p\right) \mid \epm \in A_k\right) \rightarrow -\frac{t p^2(q-1)^2}{q^2f^\dprime_{\beta,h}(s_k)}\left(m_{k,1}^{p-1}-m_{k,2}^{p-1}\right)^2
    $$
   for $k\in \{1,2\}$, as $N \rightarrow \infty$. Suppose that $t>0$. Then, by Lemma \ref{pcl_b}, we know that $$\p_{\beta+N^{-\frac{1}{2}} t,h, p} (\epm \in A_1) \le C_1e^{-C_2\sqrt{N}}$$ for some constants $C_1,C_2>0$. Hence, we have:
   $$\E_{\beta+N^{-\frac{1}{2}} t, h, p}\left(N^{\frac{1}{2}}\left(\|\epm\|_p^p - \|\bm m_{2}\|_p^p\right)\right) \rightarrow -\frac{t p^2(q-1)^2}{q^2f^\dprime_{\beta,h}(s_2)}\left(m_{2,1}^{p-1}-m_{2,2}^{p-1}\right)^2.$$
   Hence, for $t > 0$, we have:
   \begin{equation}\label{1811}
       \P_{\beta, h, p}\left(N^{\frac{1}{2}}\left(\hat{\beta}_N-\beta\right) \leq t\right) \rightarrow p_1 + p_2 \p\left(\mathcal{N}\left(0,-\frac{q^2 f_{\beta,h}''(s_2)}{p^2(q-1)^2}(m_{2,1}^{p-1}-m_{2,2}^{p-1})^{-2}\right) \le t\right)
   \end{equation}
Next, suppose that $t<0$. Then, by Lemma \ref{pcl_b}, we know that $$\p_{\beta+N^{-\frac{1}{2}} t,h, p} (\epm \in A_2) \le C_1e^{-C_2\sqrt{N}}$$ for some constants $C_1,C_2>0$. Hence, we have:
   $$\E_{\beta+N^{-\frac{1}{2}} t, h, p}\left(N^{\frac{1}{2}}\left(\|\epm\|_p^p - \|\bm m_{1}\|_p^p\right)\right) \rightarrow -\frac{t p^2(q-1)^2}{q^2f^\dprime_{\beta,h}(s_1)}\left(m_{1,1}^{p-1}-m_{1,2}^{p-1}\right)^2.$$
   Hence, for $t < 0$, we have:
   \begin{equation}\label{1812}
       \P_{\beta, h, p}\left(N^{\frac{1}{2}}\left(\hat{\beta}_N-\beta\right) \leq t\right) \rightarrow p_1 \p\left(\mathcal{N}\left(0,-\frac{q^2 f_{\beta,h}''(s_1)}{p^2(q-1)^2}(m_{1,1}^{p-1}-m_{1,2}^{p-1})^{-2}\right) \le t\right)
   \end{equation}
Part (1) now follows from \eqref{1811} and \eqref{1812}.
\vspace{0.15in}
    
\noindent (2)~Suppose that $(\beta,h) \in \cC_{p,q}^2$. 
Then $h=0$, and by Lemma \ref{crit_part} ii., all possible permutations of $\bm m := \bm x_s$ for some $s>0$ are precisely the maximizers of $H_{\beta,0}$. Note that the probability measure $\p_{\beta+N^{-\frac{1}{2}}t,0} \epm^{-1}$ is permulation invariant, and hence, assigns equal mass to all these maximizers. By Lemma \ref{cruciallemma8} ii., we have:
    $$
        \E_{\beta+N^{-\frac{1}{2}} t, 0, p}\left(N^{\frac{1}{2}}\left(\|\epm\|_p^p - \|\bm m_{k}\|_p^p\right) \mid \epm \in A_k\right) \rightarrow -\frac{t p^2(q-1)^2}{q^2f^\dprime_{\beta,h}(s)}\left(m_{1}^{p-1}-m_{2}^{p-1}\right)^2
    $$
   for $k\in [q]$, as $N \rightarrow \infty$. Hence,
   $$ \E_{\beta+N^{-\frac{1}{2}} t, 0, p}\left(N^{\frac{1}{2}}\left(\|\epm\|_p^p - \|\bm m_{k}\|_p^p\right) \right) \rightarrow -\frac{t p^2(q-1)^2}{q^2f^\dprime_{\beta,h}(s)}\left(m_{1}^{p-1}-m_{2}^{p-1}\right)^2$$
   Hence, observing that $\p_{\beta,0,p}(\epm \in A_k) = q^{-1}$ for all $k\in [q]$, we have:
   \begin{equation*}
       \P_{\beta, h, p}\left(N^{\frac{1}{2}}\left(\hat{\beta}_N-\beta\right) \leq t\right) \rightarrow \p\left(\mathcal{N}\left(0, \frac{q^2 f_{\beta,h}''(s)}{p^2 (q-1)^2} (m_1^{p-1} - m_2^{p-2})^{-2}\right) \le t\right)
   \end{equation*}
   for all $t\in \mathbb{R}$. This completes the proof of part (2).
   \vspace{0.15in}
    
\noindent (3)~Suppose that $(\beta,h)=(\beta_c,0)$. Then once again, $h=0$, and $H_{\beta,0}$ has $q+1$ global maxizers, which are all permutations of $\bm x_s$ for some $s>0$, and the vector $\bm x_0$. So, $\bm m_1 = \bm x_0$ and without loss of generality, let $\bm m_2 = \bm x_s$. If $t>0$, once again by Lemma \ref{pcl_b}, $$\p_{\beta+N^{-\frac{1}{2}}t,0,p}(\epm \in A_1) \le C_1e^{-C_2\sqrt{N}}$$ for some constant $C>0$, and
 $$
        \E_{\beta+N^{-\frac{1}{2}} t, h, p}\left(N^{\frac{1}{2}}\left(\|\epm\|_p^p - \|\bm m_{2}\|_p^p\right) \mid \epm \in A_2\right) \rightarrow -\frac{t p^2(q-1)^2}{q^2f^\dprime_{\beta,h}(s)}\left(m_{2,1}^{p-1}-m_{2,2}^{p-1}\right)^2
    $$
  Hence, for $t>0$, we have:
  $$\E_{\beta+N^{-\frac{1}{2}} t, h, p}\left(N^{\frac{1}{2}}\left(\|\epm\|_p^p - \|\bm m_{2}\|_p^p\right)\right) \rightarrow -\frac{t p^2(q-1)^2}{q^2f^\dprime_{\beta,h}(s)}\left(m_{2,1}^{p-1}-m_{2,2}^{p-1}\right)^2$$ and hence,
\begin{equation}\label{2011}
       \P_{\beta, h, p}\left(N^{\frac{1}{2}}\left(\hat{\beta}_N-\beta\right) \leq t\right) \rightarrow (1-p_1) \p\left(\mathcal{N}\left(0,-\frac{q^2 f_{\beta,h}''(s_2)}{p^2(q-1)^2}(m_{2,1}^{p-1}-m_{2,2}^{p-1})^{-2}\right) \le t\right) + p_1
   \end{equation}

   Finally, for $t<0$, again we have:
   $$\p_{\beta+N^{-\frac{1}{2}}t,0,p}(\epm \in A_2) \le C_1e^{-C_2\sqrt{N}}$$
   for some constants $C_1,C_2>0$, and
   $$
        \E_{\beta+N^{-\frac{1}{2}} t, h, p}\left(N\left(\|\epm\|_p^p - \|\bm m_{1}\|_p^p\right) \mid \epm \in A_1\right) \rightarrow \frac{p(p-1)}{2q^{p-2}} \e(\bm W^\top \bm W)
    $$
    which is thus also the limiting law of the unconditional expectation, where $\bm W$ is defined as in the statement of Lemma \ref{cruciallemma8} ii.
Therefore, for $t< 0$, one has:
\begin{equation}\label{2211}
    \P_{\beta, h, p}\left(N^{\frac{1}{2}}\left(\hat{\beta}_N-\beta\right) \leq t\right) \rightarrow p_1\p\left(\bm W^\top \bm W \le \e(\bm W^\top \bm W)\right) = p_1\gamma_1.
\end{equation}
    Part (3) now follows from \eqref{2011} and \eqref{2211}, and the proof of Theorem \ref{btml44} is now complete.

\section{Technical Lemmas}
In this section, we prove some technical lemmas necessary for showing the main results of this paper.

\begin{lem}\label{amcard}
    For each $\bm v \in S_N^q \bigcap \cP_q$, we have:
    $$ \exp\left(-N \sum_{r=1}^q v_r \log v_r \right) ~\lesssim_q ~|A_N(\bm v)|~ \lesssim_q~  N^\frac{1}{2} \exp\left(-N \sum_{r=1}^q v_r \log v_r\right)~. $$
\end{lem}
\begin{proof}
    To begin with, let us assume that all entries of $\bm v$ are strictly positive. Note that,
    $$|A_N(\bm v)| = \frac{N!}{\prod_{r=1}^q (Nv_r)!}$$
    Using Stirling's formula, one can easily derive that for every positive integer $k$,
    $$\sqrt{2\pi k}\left(\frac{k}{e}\right)^k ~<~ k!~<~ 2\sqrt{2\pi k}\left(\frac{k}{e}\right)^k~. $$
    Using this bound, one has:
    \begin{eqnarray*}
        &&  \frac{\sqrt{2\pi N}}{2^q \prod_{r=1}^q \sqrt{2\pi Nv_r} v_r^{Nv_r} } ~<~|A_N(\bm v)| ~<~ \frac{2\sqrt{2\pi N}}{ \prod_{r=1}^q \sqrt{2\pi Nv_r} v_r^{Nv_r} }\\&\implies&  2^{-q} \prod_{r=1}^q v_r^{-Nv_r}~\lesssim~ |A_N(\bm v)|~\lesssim~ N^{\frac{1}{2}}  \prod_{r=1}^q v_r^{-Nv_r}\\&\implies & \exp\left(-N\sum_{r=1}^q v_r\log v_r\right) ~\lesssim_q~ |A_N(\bm v)| ~\lesssim~ N^{\frac{1}{2}} \exp\left(-N\sum_{r=1}^q v_r\log v_r\right)~.
    \end{eqnarray*}
    This proves Lemma \ref{amcard} when all entries of $\bm v$ are strictly positive. Now, if $S(q)$ denotes the statement of Lemma \ref{amcard} for $\bm v \in (0,1]^q$, then the case $\bm v$ has some zero entries is essentially same as $S(t)$, where $t$ is the number of non-zero entries of $\bm v$. Since we have proved $S(q)$ for all $q\ge 1$, $S(t)$ should also be valid, which completes the proof of Lemma \ref{amcard}.
\end{proof}

Similar to the spirit of the proof of Theorem \ref{conclem}, we now give an approximation to the non-normalized probability mass function of the empirical magnetization $\epm$. Towards this, for every $\varepsilon \ge 0$, let us define:
$$\cP_{q,\varepsilon} := \{\bm v \in \cP_q~:~ \min v_i \ge \varepsilon\}~,~\cP_{q,\varepsilon,N} := \cP_{q,\varepsilon} \bigcap S_N^q ~,~\cP_{q,0^+} := \bigcup_{\varepsilon > 0} \cP_{q,\varepsilon}~,~ \cP_{q,0^+,N} := \cP_{q,0^+} \bigcap S_N^q~. $$

\begin{lem}\label{densappr}
    For $\bm v \in \cP_{q,0^+,N}$, we have:
    $$q^N Z_N(\beta,h) \p_{\beta,h,N}(\epm = \bm v) = (1+r_{\beta. h,N}(\bm v)) N^{-\frac{q-1}{2}} A(\bm v) e^{N H_{\beta,h}(\bm v)}$$ where $A(\bm v) := (2\pi)^{-(q-1)/2} \prod_{r=1}^q v_r^{-1/2}$ and for any $\varepsilon > 0$,
    $$\lim_{N\rightarrow \infty} ~\sup_{\bm v \in \cP_{q,\varepsilon,N}}~ \sup_{\beta,h} ~|r_{\beta,h,N}(\bm v)| = 0~. $$
\end{lem}

\begin{proof}
    For an $\bm v \in \cP_{q,0^+,N}$, we have:
    \begin{equation}\label{prst}
        q^N Z_N(\beta,h) \p_{\beta,h,N}(\epm = \bm v) = |A_N(\bm v)| \exp \left\{N \left(\beta  \sum_{r=1}^q v_r^p + hv_1\right)\right\}~.
    \end{equation}
    By Stirling's formula, we have:

    \begin{equation*}
        |A_N(\bm v)| = \frac{N!}{\prod_{r=1}^q (Nv_r)!} = (2\pi N)^{\frac{1-q}{2}}\left(\prod_{r=1}^q v_r^{-\frac{1}{2}}\right) e^{-N\sum_{r=1}^q v_r \log v_r}(1+o_{N,\bm v}(1))~.
    \end{equation*}
    where the $o_{N,v}(1)$ term goes to $0$ uniformly over all $\bm v \in \cP_{q,\varepsilon,N}$ for any $\varepsilon > 0$.
    Therefore, we have from \eqref{prst},
    \begin{equation*}
        q^N Z_N(\beta,h) \p_{\beta,h,N}(\epm = \bm v) = (2\pi N)^{\frac{1-q}{2}}\left(\prod_{r=1}^q v_r^{-\frac{1}{2}}\right)e^{N H_{\beta,h}(\bm v)}(1+o_{N,\bm v}(1))~~.
    \end{equation*}
    This completes the proof of Lemma \ref{densappr}.
\end{proof}

\begin{lem}\label{denselem}
    Given any $\bm v \in \cP_q$, there exists a sequence $\bm v_N \in S_N^q\bigcap \cP_q$ such that $\|\bm v_N - \bm v\|_\infty \lesssim_q \frac{1}{N}$, and consequently, $\bm v_N \rightarrow \bm v$.
\end{lem}
\begin{proof}
    For $\bm v \in \cP_q$, define $$\bm v_N := \left( \frac{\lfloor Nv_1 \rfloor}{N},\ldots, \frac{\lfloor Nv_{q-1} \rfloor}{N},1-\frac{\sum_{r=1}^{q-1}\lfloor Nv_{r} \rfloor}{N}\right)^\top~. $$
    Then, $\bm v_N \in S_N^q\bigcap \cP_q$ for each $N \ge 1$. Clearly, $|v_r - v_{N,r}| < \frac{1}{N}$ for all $r\in [q-1]$. Also,
    $$|v_q-v_{N,q}| = \frac{\left|\sum_{r=1}^{q-1} \left(\lfloor N v_r \rfloor - Nv_r\right)\right|}{N} \le \frac{q-1}{N}~. $$ This completes the proof of Lemma \ref{denselem}.
\end{proof}


\begin{lem}\label{sectay}
    If $\bm m_*$ is a global maximizer of the function $H_{\beta,h}$ in the interior of $\cP_q$, and if $\bm m_*+\bm u \in \cP_q$, then:
    $$H_{\beta,h}(\bm m_* + \bm u) = H_{\beta,h}(\bm m_*) + \frac{1}{2}\sum_{r=1}^q\left[\beta p(p-1) (\bm m_*+\alpha \bm u)_r^{p-2} - \frac{1}{(\bm m_*+\alpha \bm u)_r}\right] u_r^2$$
    for some $\alpha \in [0,1]$.
\end{lem}
\begin{proof}
    Define $\phi:[0,1]\to \R$ as $\phi(t) := H_{\beta,h}(\bm m_* + t \bm u)$. Then, convexity of $\cP_q$ implies that $\bm m_* + t \bm u \in \cP_q$ for all $t\in [0,1]$. By a second-order Taylor expansion of the function $\phi$, we have:
    $$\phi(1) = \phi(0) + \phi'(0) + \frac{1}{2}\phi''(\alpha)\quad\text{for some}~\alpha\in [0,1]~. $$
    Note that $\phi(1) = H_{\beta,h}(\bm m_*+\bm u)$ and $\phi(0)= H_{\beta,h}(\bm m_*)$. Also, $\phi'(0)=0$ because $\phi$ is maximized at $0$, and since $\bm m_*$ is in the interior of $\cP_q$, the domain of $\phi$ can be extended to $[-\delta,1]$ within $\cP_q$ for some $\delta>0$. Now, $$\phi'(t) = \sum_{r=1}^q u_r\nabla_r H_{\beta,h}(\bm m_* + t\bm u)\implies \phi''(t) = \sum_{r,s=1}^q u_r u_s \nabla_{r,s}^2 H_{\beta,h}(\bm m_* + t\bm u)~. $$
    Finally, observe that $$\nabla_{r,s}^2 H_{\beta,h}(\bm m_* + t\bm u) = \left[\beta p(p-1)(\bm m_* + t\bm u)_r^{p-2} - \frac{1}{(\bm m_* + t\bm u)_r}\right] \mathbbm{1}_{r=s}~. $$
    This proves Lemma \ref{sectay}.
\end{proof}

\begin{lem}\label{sectay2}
    If $(\beta, h)$ is regular, then for $\varepsilon >0$ small enough, there exists $\alpha >0$ such that for any $N$ large enough and any $\bm w \in \cH_q$ with $K<\|\bm w\|\le \varepsilon N^{\frac{1}{2}}$, we have:
    $$N H_{\beta_N,h_N}\left(\bm m_* + N^{-\frac{1}{2}}\bm w\right) \le N H_{\beta_N,h_N}\left(\bm m_*\right) - \alpha \|\bm w\|^2~. $$
\end{lem}
\begin{proof}
    First, notice that from \eqref{final_reduc_reg}, we have:
    $$NH_{\beta_{N}, h_{N}}\left(\bm m_*+N^{-\frac{1}{2}} \bm w\right)=NH_{\beta_N,h_N}(\bm m_*) + \frac{1}{2} \bm Q_{\bm m_*,\beta}(\bm w)
        +\ip{\bar{\beta}p\bm m_*^{p-1}+\bar{h}\bm e_1}{\bm w}+o\left(1\right). $$
    Now, for $\varepsilon>0$, for all $N$ large enough, one has:
    $$NH_{\beta_{N}, h_{N}}\left(\bm m_*+N^{-\frac{1}{2}} \bm w\right)\leq NH_{\beta_N,h_N}(\bm m_*) + \frac{1}{2} \bm Q_{\bm m_*,\beta}(\bm w)
        +\ip{\bar{\beta}p\bm m_*^{p-1}+\bar{h}\bm e_1}{\bm w}+\epsilon. $$
    It is easy to check by Cauchy-Schwarz inequality that $$\ip{\bar{\beta}p\bm m_*^{p-1}+\bar{h}\bm e_1}{\bm w(\bm v)}\leq D\nrm{\bm w(\bm v)},$$
    for some constant $D>0$.
    Since $\frac{1}{2}\bm Q_{\bm m_*,\beta}$ is negative definite, it is dominated by $-C\|\bm w\|^2$ for some constant $C>0$ (here c depends on $\beta$ and $\bm m_*$). For large $K$, we have $-C^\prime\nrm{\bm w}^2\geq -C\nrm{\bm w}^2+D\nrm{\bm w}$ where $C^\prime>0$. Now, choose $\varepsilon$ small enough such that $C^\prime-\frac{\varepsilon}{K^2}>0$. Let $0<\alpha<C^\prime-\frac{\varepsilon}{K^2}$. Hence, for all $K<\|\bm w\|\le \varepsilon N^{\frac{1}{2}}$, we have Lemma \ref{sectay2}.
\end{proof}
\begin{lem}\label{l2nm}
    For any $\beta,h, \beta_N,h_N$ and any $\bm t \in \cP_q$, the following holds:
    $$H_{\beta_N,h_N}(\bm t) = H_{\beta,h}(\bm t) + (\beta_N-\beta)\sum_{r=1}^q t_r^p + (h_N-h)t_1~. $$
\end{lem}
\begin{proof}
    Note that the right-hand side equals:
    $$\beta\sum_{r=1}^q t_r^p +ht_1 - \sum_{r=1}^q t_r\log t_r + (\beta_N-\beta)\sum_{r=1}^q t_r^p + (h_N-h)t_1 = \beta_N\sum_{r=1}^q t_r^p +h_Nt_1 - \sum_{r=1}^q t_r\log t_r$$ and the last term equals $H_{\beta_N,h_N}(\bm t)$.
\end{proof}

\begin{lem}\label{spectay}
    Let $(\beta,h)\in \cS_{p,q}$ be such that $\bm x_s$ is the unique global maximizer of $H_{\beta,h}$. Let $\bm u :=(1-q,1,\ldots,1)$. Then, for any $t \in \R$ and $\bm v\in \cH_q \cap \bm u^\perp$ such that $\bm x + t\bm u + \bm v \in \cP_{q,0^+}$, there is some $\alpha,\alpha^\prime \in (0,1)$ such that,
    \begin{eqnarray*}
         &&H_{\beta,h}(\bm m+t\bm u+\bm v)\\ &=& H_{\beta,h}(\bm m)+\frac{1}{2}\bm Q_{\bm m+tu+\alpha \bm v,\beta}(\bm v)+\frac{t^4}{24}\sum_{r=1}^{q}\left[\frac{\beta p!}{(p-4)!}(m_r+\alpha^\prime tu_r)^{p-4}-\frac{2}{(m_r+\alpha^\prime tu_r)^3}\right]u_r^4 
    \end{eqnarray*}     
    Furthermore,
    $$
            \frac{1}{2}\bm Q_{\bm x_s,\beta}(\bm v) =\frac{1}{2}k^\dprime\left(\frac{1-s}{q}\right)\nrm{\bm v}^2\quad \text{and} \quad
            \frac{1}{24}\sum_{r=1}^{q}\left[\frac{\beta p!}{(p-4)!}(m_r)^{p-4}-\frac{2}{(m_r)^3}\right]u_r^4 =\frac{1}{24}q^4 f^{(4)}_{\beta,h}(s).
    $$
\end{lem}
\begin{proof}
    Using Taylor expansion on the function $H_{\beta,h}(\bm m)$ we get:
    
        \begin{eqnarray*}
            &&H_{\beta,h}(\bm m+t\bm u+\bm v)\\ & =& H_{\beta,h}(\bm m+t\bm u) + \nabla H_{\beta,h}(\bm m+t\bm u) \cdot \bm v + \frac{1}{2}\bm Q_{\bm m+tu+\alpha \bm v,\beta}(\bm v)                                                                                    \\
                                            & =& H_{\beta,h}(\bm m+t\bm u) + \frac{1}{2}\bm Q_{\bm m+tu+\alpha \bm v,\beta}(\bm v)                                                                                                                                  \\
                                            & = & H_{\beta,h}(\bm m)+\frac{1}{2}\bm Q_{\bm m+tu+\alpha \bm v,\beta}(\bm v)+\frac{t^4}{24}\sum_{r=1}^{q}\left[\frac{\beta p!}{(p-4)!}(m_r+\alpha^\prime tu_r)^{p-4}-\frac{2}{(m_r+\alpha^\prime tu_r)^3}\right]u_r^4,
        \end{eqnarray*}
    
    where $\alpha, \alpha^\prime \in (0,1)$.
    The dot product, $\nabla H_{\beta,h}(\bm m+t\bm u) \cdot \bm v$ is zero as the last $q-1$ coordinates of $\nabla H_{\beta,h}(\bm m+t\bm u)$ are equal and $v_1 = v_2+\ldots+v_q = 0$. We can guarantee that the last $q-1$ coordinates are same because there is a unique maximizer for $(\beta,h) \in \cS_{p,q}$. The last equality is using Taylor expansion again on the function $b(t) := H_{\beta,h}(\bm m+t\bm u)$. $b(t)$ is a function in one variable with maximum at $t=0$. Also, the point is a special point and hence, $b^\dprime(0)$. Therefore, by higher derivative test $b^{(3)}(0)=0$\cite{deriv_test}.
    
    Also, from \eqref{q_form}, we have:
    $$
        \begin{aligned}
            \frac{1}{2}\bm Q_{\bm x_s,\beta}(\bm v) & =\frac{1}{2}k^\dprime\left(\frac{1-s}{q}\right)\sum_{i=2}^{q}v_i^2 \\
                                                  & =\frac{1}{2}k^\dprime\left(\frac{1-s}{q}\right)\nrm{\bm v}^2.
        \end{aligned}
    $$
    Moreover,
    $$
        \begin{aligned}
            \frac{1}{24}\sum_{r=1}^{q}\left[\frac{\beta p!}{(p-4)!}(m_r)^{p-4}-\frac{2}{(m_r)^3}\right]u_r^4 & = \frac{1}{24}\sum_{r=1}^{q}k^{(4)}(m_r)u_r^4                                       \\
                                                                                                             & =\frac{1}{24}(q-1)^4k^{(4)}\left(\frac{1+(q-1)s}{q}\right)+\frac{1}{24}(q-1)k^{(4)}\left(\frac{1-s}{q}\right) \\
                                                                                                             & = \frac{1}{24}q^4 f^{(4)}_{\beta,h}(s).
        \end{aligned}
    $$
    This completes the proof of Lemma \ref{spectay}.
\end{proof}
\begin{lem}
    Consider $(\beta,h)\in \mathcal{S}_{p,q}^1$. Suppose that $\beta_N= \beta+\frac{\bar{\beta}}{N^{\frac{3}{4}}}$ and $h_N= h+\frac{\bar{h}}{N^{\frac{3}{4}}}$, and let $\bm m= \bm x_s \in \cP_q$ be the unique global maximizer of $H_{\beta, h}$.
    \begin{enumerate}[i.]
        \item For any $M>0$,
              \begin{multline}
                  H_{\beta_N, h_N}\left(\bm m+N^{-1 / 4} t \bm u+N^{-1 / 2} \bm v\right)=   H_{\beta, h}(\bm m)+ \frac{\bar{\beta}}{N^{\frac{3}{4}}}\nrm{\bm m}_p^p+\frac{\bar{h}}{N^{\frac{3}{4}}} m_1+(\bar{\beta}p\ip{\bm m^{p-1}}{\bm u}+\bar{h}(1-q))\frac{t}{N}                                     \\
                  \frac{1}{2N} k^{\prime \prime}\left(\frac{1-s}{q}\right)\nrm{\bm v}^2+\frac{1}{24N} q^{4} \f^{(4)}(s) t^4+o(N^{-1})
                  \label{taylor_with_o}
              \end{multline}
              uniformly over $\bm v \in \cH_q \cap \bm u^{\perp} \cap B(0, M)$ and $t \in[-M, M]$.
        \item For large enough $M$, for $N$ large enough, for any $\bm v \in \cH_q \cap \bm u^{\perp}$ and $t \in \R \setminus [-M, M]$, there exists $c_1\geq 0$ and $c_2>0$ such that
              \begin{equation}
                  N H_{\beta_N, h_N}\left(\bm m+N^{-1 / 4} t \bm u+N^{-1 / 2} \bm v\right) \leq  N H_{\beta, h}(\bm m)+ N^{\frac{1}{4}}\bar{\beta}\nrm{\bm m}_p^p+N^{\frac{1}{4}}\bar{h} m_1-c_1\|\bm v\|^2-c_2t^4.
                  \label{tight_tay1}
              \end{equation}
    \end{enumerate}
\end{lem}

\begin{proof}
     To begin with, note that:
    $$
        \begin{aligned}
             & H_{\beta_{N}, h_{N}}(\bm m)=H_{\beta, h}(\bm m)+\frac{\bar{\beta}}{N^{\frac{3}{4}}} \sum_{i=1}^{q} m_{i}^{p}+\frac{\bar{h}}{N^{\frac{3}{4}}} m_{1}
        \end{aligned}
    $$
    Now, let $\bm \omega_{N}=N^{-1 / 4} t \bm u+N^{-1 / 2} \bm v$, whence we have:
    \begin{multline*}
        H_{\beta_{N}, h_{N}}\left(\bm m+\bm \omega_{N}\right)=H_{\beta, h}\left(\bm m+\bm \omega_{N}\right)+\frac{\bar{\beta}}{N^{\frac{3}{4}}}\|\bm m\|_{p}^p+\frac{\bar{h}}{N^{\frac{3}{4}}} m_{1}                                                                                                             \\
        +\frac{\bar{\beta}}{N^{\frac{3}{4}}} \sum_{i=1}^{q}\left[\left(m_{i}+\omega_{N,i}\right)^{p}-m_{i}^{p}\right]+\frac{\bar{h}}{N^{\frac{3}{4}}} \omega_{N 1}                                                                                                                                               \\
    \end{multline*}
    Hence,
    \begin{equation}
        H_{\beta_{N}, h_{N}}\left(\bm m+\bm \omega_{N}\right)=H_{\beta, h}\left(\bm m+\bm \omega_{N}\right)+\frac{\bar{\beta}}{N^{\frac{3}{4}}}\|\bm m\|_{p}^p+\frac{\bar{h}}{N^{\frac{3}{4}}} m_{1}+(\bar{\beta}p\ip{\bm m_*^{p-1}}{\bm u}+\bar{h}(1-q))\frac{t}{N}+o\left(N^{-1}\right)
        \label{taylor_1}
    \end{equation}
    We also have the following:

    \begin{equation}
        \frac{1}{2} \bm Q_{\bm m+N^{-1 / 4} t \bm u+\alpha N^{-1 / 2} \bm v, \beta}\left(N^{-1 / 2} v\right)=\frac{1}{2N} k^{\prime \prime}\left(\frac{1-s}{q}\right)\nrm{\bm v}^{2}+o\left(N^{-1}\right)
        \label{taylor_2}
    \end{equation}
    \begin{equation}
        \frac{t^{4}}{24} \sum_{r=1}^{q}\left[\frac{\beta p !}{(p-4) !}\left(m_{r}+\alpha^{\prime} t N^{-1/4} u_{r}\right)^{p-4}-\frac{2}{\left(m_{r}+\alpha^{\prime} t N^{-1/4} u_{r}\right)^{3}}\right] u_{r}^{4}N^{-1}=\frac{1}{24N} q^{4} f_{\beta, h}^{(4)}(s) t^4+o\left(N^{-1}\right)
        \label{taylor_3}
    \end{equation}
    for any $\alpha'\in (0,1)$.
    Hence, by putting the equations \eqref{taylor_1}, \eqref{taylor_2} and \eqref{taylor_3} together and using Lemma \ref{spectay}, we get \eqref{taylor_with_o}.
    \par

    To prove \eqref{tight_tay1}, note that 
    since $(\beta, h)$ is a type-$I$ special point, we have $f_{\beta, h}^{(4)}(s)<0$. Inequality \eqref{tight_tay1} now follows from \eqref{taylor_with_o}.
\end{proof}

\begin{lem}\label{spec2tay}
    Let $(\beta,h)\in \cS^2_{4,2}$. Theb $\bm m=(1/2,1/2)$ as the unique global maximizer of $H_{\beta,h}$. Let $\bm u :=(-1,1)$. For any $t \in \R$ and such that $\bm m + t\bm u \in \cP_{q,0^+}$, there is some $\alpha \in (0,1)$ such that,
    $$
        H_{\beta,h}(\bm m+t\bm u)=H_{\beta,h}(\bm m)-\frac{t^6}{30}\sum_{r=1}^{2}\left[\frac{u_r^6}{(m_r+\alpha t u_r)^5}\right]
    $$
    Furthermore,
    $$
        \frac{1}{30}\sum_{r=1}^{2}\frac{u_r^6}{m_r^5} =\frac{32}{15}.
    $$
\end{lem}
\begin{proof}
    Again using Taylor expansion with Lagrange reminder, we have:
    $$
        H_{\beta,h}(\bm m+t\bm u)=H_{\beta,h}(\bm m)-\frac{t^6}{30}\sum_{r=1}^{2}\left[\frac{u_r^6}{(x_r+\alpha t u_r)^5}\right].
    $$
   The first to fifth order derivatives vanish by the derivative test\cite{deriv_test}. The second conclusion directly follows from the values of $\bm u$ and $\bm m$
\end{proof}

\begin{lem}\label{dpelv}
    Assume $(\beta,h)$ is a type-$II$ special point. Let $\beta_N-\beta=\frac{\bar{\beta}}{N^{\frac{5}{6}}}$ and $h_N-h=\frac{\bar{h}}{N^{\frac{5}{6}}}$. Now, $\bm m=(1/2,1/2)$ is the global maximizer. Then,
    \begin{enumerate}[i.]
        \item For any $M>0$,
              \begin{equation}
                  H_{\beta_{N}, h_{N}}\left(\bm m+N^{-1 / 6} t \bm u\right)=H_{\beta, h}\left(\bm m\right)-\frac{32}{15N}t^6+\frac{\bar{\beta}}{N^{\frac{5}{6}}}\|\bm m\|_{p}^p+\frac{\bar{h}}{N^{\frac{5}{6}}} m_{1}-\frac{\bar{h} t}{N}+o(N^{-1})
                  \label{taylorspec_with_o}
              \end{equation}
              uniformly over $\bm v \in \cH \cap \bm u^{\perp} \cap B(0, M)$ and $t \in[-M, M]$.
        \item For large enough $M$, for $N$ large enough, for any $t \in \R \setminus [-M, M]$ there exists some $c>0$ such that
              \begin{equation}
                  N H_{\beta_{N}, h_{N}}\left(\bm m+N^{-1 / 6} t \bm u\right) \leq N H_{\beta, h}\left(\bm m\right)-ct^6+ N^{\frac{1}{6}}\bar{\beta}\|\bm m\|_{p}^p+N^\frac{1}{6}\bar{h} m_{1}
                  \label{tight_tay2}
              \end{equation}
    \end{enumerate}
\end{lem}

\begin{proof}
    Let $\bm \omega_{N}=N^{-1 / 6} t \bm u$. Again by Lemma \ref{l2nm},
    \begin{multline*}
        H_{\beta_{N}, h_{N}}\left(\bm m+\bm \omega_{N}\right)=H_{\beta, h}\left(\bm m+\bm \omega_{N}\right)+\frac{\bar{\beta}}{N^{\frac{5}{6}}}\|\bm m\|_{p}^p+\frac{\bar{h}}{N^{\frac{5}{6}}} m_{1}+\frac{\bar{\beta}}{N^{\frac{5}{6}}} \sum_{i=1}^{q}\left[\left(m_{i}+\omega_{N,i}\right)^{p}-m_{i}^{p}\right]+
        \frac{\bar{h}}{N^{\frac{5}{6}}} \omega_{N,1}
    \end{multline*}
    Since $\beta_N-\beta=o(N^{-\frac{5}{6}})$ and $h_N-h=o(N^{-\frac{5}{6}})$, we have:
    \begin{equation}
        H_{\beta_{N}, h_{N}}\left(\bm m+\bm \omega_{N}\right)=H_{\beta, h}\left(\bm m+\bm \omega_{N}\right)+\frac{\bar{\beta}}{N^{\frac{5}{6}}}\|\bm m\|_{p}^p+\frac{\bar{h}}{N^{\frac{5}{6}}} m_{1}- \frac{\bar{h}t}{N}+o\left(N^{-1}\right)
        \label{taylorspec_1}
    \end{equation}
    By Lemma \ref{spec2tay}, we also have:

    \begin{equation}
        H_{\beta, h}\left(\bm m+\bm \omega_{N}\right)=H_{\beta, h}\left(\bm m\right)-\frac{32}{15N}t^6+o(N^{-1})
        \label{taylorspec_2}
    \end{equation}

    Putting the equations \eqref{taylorspec_1} and \eqref{taylorspec_2} together, we get \eqref{taylorspec_with_o}.

    Since $t^6$ dominates $t$ for large $t$,  \eqref{taylorspec_with_o} follows for large $t$ and large $N$,
    $$N H_{\beta_{N}, h_{N}}\left(\bm m+N^{-1 / 6} t \bm u\right) \leq N H_{\beta, h}\left(\bm m\right)-ct^6+ N^{\frac{1}{6}}\bar{\beta}\|\bm m\|_{p}^p+N^\frac{1}{6}\bar{h} m_{1},$$
    for some $c>0$. This proves \eqref{tight_tay2}, and completes the proof of Lemma \ref{dpelv}.
\end{proof}
\section{Perturbative Concentration Lemmas at critical points}\label{sec:partsec}

In this section, we analyse the concentration behavior of $\epm$ at critical points, when the model parameters are perturbed by a factor of $N^{-\frac{1}{2}}$. These results will be crucial in deriving the asymptotics of the ML estimates at the critical points.

\begin{lem}   \label{pcl_h}
    Let $A$ be a set whose interior contains exactly one maximizer $\bm m$ of $H_{\beta,h}$, and whose closure does not include any other maximizer. Let $\bar{h}\ne 0$ be given. Also, suppose that there exists $\bm m^\prime \in \mxs$ with $sgn(m_1-m_1') = -sgn(\bar{h})$. Then there exist positive constants $C_1$ and $C_2$ not depending on $N$, such that:
        \begin{equation*}
        \mathbb{P}_{\beta, h+\frac{\bar{h}}{\sqrt{N}}, N}\left(\epm \in A\right) \leqslant C_1 e^{-C_2\sqrt{N}}.
    \end{equation*}     
\end{lem}
\begin{proof}
Denote $h_N = h+\frac{\bar{h}}{\sqrt{N}}$. It follows from the proof of Theorem \ref{pcr}, that for every $\varepsilon > 0$ sufficiently small,

 $$\frac{\P_{\beta, h_N, N}\left(\epm \in B(\bm m, \varepsilon)\right)}{\P_{\beta, h_N, N}\left(\epm \in B\left(\bm m^\prime, \varepsilon\right)\right)}\sim e^{\bar{h}\sqrt{N}(m_1-m_1^\prime)} \frac{\int_{\cH_q} e^{\bar{h} w_1 +\frac{1}{2} Q_{\bm m, \beta}(\bm w)} dw_1dw_2\ldots dw_q}{\int_{\cH_q} e^{\bar{h} w_1 +\frac{1}{2} Q_{\bm m', \beta}(\bm w)} dw_1dw_2\ldots dw_q},$$
which immediately implies that $\P_{\beta, h_N, N}\left(\epm \in B(\bm m, \varepsilon)\right) \le K_1 e^{-K_2\sqrt{N}}$ for some constants $K_1,K_2>0$. Lemma \ref{pcl_h} now follows from Theorem \ref{conclem}.
\end{proof}

\begin{lem} \label{pcl_b}
   Let $A$ be a set whose interior contains exactly one maximizer $\bm m$ of $H_{\beta,h}$, and whose closure does not include any other maximizer. Let $\bar{\beta}\ne 0$ be given. Also, suppose that there exists $\bm m^\prime \in \mxs$ with $sgn(\|\bm m\|_p^p-\|\bm m'\|_p^p) = -sgn(\bar{\beta})$. Then there exist positive constants $C_1$ and $C_2$ not depending on $N$, such that:
        \begin{equation*}
        \mathbb{P}_{\beta + \frac{\bar{\beta}}{\sqrt{N}}, h, N}\left(\epm \in A\right) \leqslant C_1 e^{-C_2\sqrt{N}}.
    \end{equation*}     
\end{lem}
\begin{proof}
Denote $\beta_N = \beta+\frac{\bar{\beta}}{\sqrt{N}}$. It follows from the proof of Theorem \ref{pcr}, that for every $\varepsilon > 0$ sufficiently small,
    $$\frac{\P_{\beta_N, h, N}\left(\epm \in B(\bm m, \varepsilon)\right)}{\P_{\beta_N, h, N}\left(\epm \in B\left(\bm m^\prime, \varepsilon\right)\right)}\sim e^{\bar{\beta}\sqrt{N}(\nrm{\bm m}_p^p-\nrm{\bm m^\prime}_p^p)}  \frac{\int_{\cH_q} e^{\langle \bar{\beta} p \bm m^{p-1},\bm w\rangle +\frac{1}{2} Q_{\bm m, \beta}(\bm w)} dw_1dw_2\ldots dw_q}{\int_{\cH_q} e^{\langle \bar{\beta} p \bm m^{\prime p-1},\bm w\rangle +\frac{1}{2} Q_{\bm m', \beta}(\bm w)} dw_1dw_2\ldots dw_q},$$
    which immediately implies that $\P_{\beta_N, h, N}\left(\epm \in B(\bm m, \varepsilon)\right) \le K_1 e^{-K_2\sqrt{N}}$ for some constants $K_1,K_2>0$. Lemma \ref{pcl_b} now follows from Theorem \ref{conclem}.
\end{proof}
\section{Maximizers of the Negative Free Energy}\label{apxf}
In this section, we give a detailed analysis of the structure of the maximizers of the negative free energy function $H_{\beta,h}$.
\begin{ppn}\label{propep1}
    Let $\beta, h \geq 0$ and let $\bm m$ be a global maximizer of $H_{\beta, h}$ in $\cP_q$.
    \begin{enumerate}[i.]
        \item The vector $\bm m$ has the coordinate $\min \left(m_i\right)$ repeated $q-1$ times at least. Also, $\min(m_i) \le 1/(\beta p(p-1))^{\frac{1}{p-1}}$.
        \item If $h>0$, then $m_1>m_i$, for all $i \in\{2, \ldots, q\}$.
        \item  The inequality $\min \left(m_i\right)>0$ holds.
        \item  For $q\geq 3$, $\min \left(m_i\right)<1 / (\beta p (p-1))^{\frac{1}{p-1}}$.
    \end{enumerate}
    \label{oned_reduce}
\end{ppn}
\begin{proof}
 By the Lagrange multiplier method, at the global maximum $\bm m$, we have
    \begin{equation}\label{fp1pre3}
                     \frac{\partial}{\partial m_i}H_{\beta, h}(\bm m)-\frac{\partial}{\partial m_j}H_{\beta, h}(\bm m)=0 
            \implies g_{\beta}(m_i)-g_{\beta}(m_j)=h\one\{j=1\}-h\one\{i=1\},
    \end{equation}
    where $g_{\beta}(z)=\beta p z^{p-1}-\log z$. Also, we have:
    \begin{equation}\label{sumder5}
        \frac{\partial^2}{\partial m_i^2}H_{\beta, h}(\bm m)+\frac{\partial^2}{\partial m_j^2}H_{\beta, h}(\bm m)\leq 0 
            \implies  g^\prime_{\beta}(m_i) +  g^\prime_{\beta}(m_j)\leq 0.
    \end{equation}

    Assume, $h=0$ then $g_\beta(m_i)=g_\beta(m_j)$. Now, $g_\beta$ is strictly convex and so, $\{m_i:i=1,\ldots,q\}$ has at most two elements. If it has two elements, then one of these must lie strictly to the left of the global minimizer of $g_{\beta}$, and the other strictly to the right. However, by \eqref{sumder5}, at most one $j$ can satisfy $g_{\beta}'(m_j) > 0$. This forces exactly one entry of $\bm m$ to be equal to the element to the right of the minimizer of $g_{\beta}$ and all the other entries to be equal to the element to the left. Also, note that $m:=\min(m_i)$ must satisfy $g_{\beta}'(m) \le 0$, which implies that $m \le 1/(\beta p(p-1))^{\frac{1}{p-1}}$. This proves $i.$ for the case $h=0$.
    
    Now, suppose that $h> 0$. Then, we have:
    $$
        H_{\beta,h}(\bm m)-H_{\beta,h}(\tilde{\bm m})=(h\one\{i=1\}-h\one\{j=1\})(m_i-m_j)
    $$
    where $\tilde{x}$ is a vector with $i$ and $j$ swapped.
    If $\bm m$ is a maximum, we must then have:
    $$
        \begin{aligned}
            (h\one\{i=1\}-h\one\{j=1\})(m_i-m_j) & \geq 0.
        \end{aligned}
    $$
    This shows that $m_1\ge m_j$ for any $j\neq 1$. By \eqref{fp1pre3}, $m_1 = m_j$ for some $j \neq 1$ is impossible, and hence, we get $ii.$ Now, we complete the proof of $i$. for the case $h>0$. Note that if $m_i \ge 1/(\beta p(p-1))^{\frac{1}{p-1}}$ for some $i\ge 2$, then $g_\beta'(m_i) \ge 0$ which implies that $g_\beta'(m_1) > 0$ (from the strict convexity of $g_\beta$). This contradicts \eqref{sumder5}, thereby implying that $m_i < 1/(\beta p(p-1))^{\frac{1}{p-1}}$ for all $i\ge 2$. Now, \eqref{fp1pre3} forces $g_\beta$ to be constant on the set $\{m_2,\ldots,m_q\}$, and all elements of these set lie to the left of the minimizer of the strictly convex function $g_\beta$, which forces them to be all equal. This completes the proof of $i$.
    
    If $\min_j(m_j)=0$, then there is some $i,j$ such that $m_j>0$ and $m_i=0$.
    \begin{equation*}
        \frac{d}{dt} H_{\beta,h}(\bm m +t(\bm e_i- \bm e_j))=g_\beta(m_i+t)-g_\beta(m_j-t)+h\one\{i=1\}-h\one\{j=1\}.
    \end{equation*}
    The derivative tends to $+\infty$ as $t \rightarrow 0^+$, contradicting that $\bm m$ is a maximizer of $H_{\beta,h}$. Hence, by contradiction $\min_j(m_j)>0$, proving $iii$.

Finally, we prove $iv$. Note that if $h>0$, then $iv$. follows directly from the proof of $i$. So, let us assume that $h=0$. Define a vector $\bm m^\beta$ as:
 $$
        \begin{aligned}
            m_j^\beta & =\frac{1}{(\beta p(p-1))^\frac{
            1
            }{p-1} } \quad \forall j \neq 1                     \\
            m_j^\beta & =\left(1-\frac{q-1}{(\beta p(p-1))^\frac{1}{p-1}}\right) \quad \text{for}~ j =1
        \end{aligned}
    $$
Note that $\bm m^\beta$ is a probability vector for $\beta \ge (q-1)^{p-1}/(p(p-1))$. Suppose that that $m_j = 1/(\beta p(p-1))^{1/(p-1)}$ for all $j> 1$, whence we have for all $j\ne 1$ (by \eqref{fp1pre3}),

$$
        \begin{aligned}
                        & \beta p m_1^{p-1}-\log m_1-\beta p m_j^{p-1}+\log m_j=0                              \\
            \Rightarrow & \beta p m_1^{p-1}-\log m_1- \frac{1}{p-1} -\frac{1}{p-1} \log \beta p(p-1)= 0       \\
            \Rightarrow & \log \left(\beta p(p-1)m_1^{p-1}\right) - \beta p(p-1)m_1^{p-1} +1 = 0 \\
        \end{aligned}
    $$
which implies that $\beta p(p-1)m_1^{p-1} =1 \implies m_1 = 1/(\beta p(p-1))^{1/(p-1)} = m_2=\ldots = m_q = 1/q$. Note that this implies $\beta = q^{p-1}/(p(p-1))$. Since $(q^{-1},\ldots,q^{-1})$ is a minimizer of $H_{\beta,0}$, we must have $f_{\beta,0}'(0) = 0$. Also, we have:
\begin{eqnarray*}
    f_{\beta,0}''(0) &=& \frac{q-1}{q^2} k''(q^{-1}) + \frac{(q-1)^2}{q^2} k''(q^{-1}) = 0
\end{eqnarray*}
since $k''(q^{-1}) =0$. The derivative test now forces $f_{\beta,0}'''(0) = 0$, i.e. 
$$\left[-\frac{q-1}{q^3} + \frac{(q-1)^3}{q^3}\right]k''(q^{-1}) =0.$$
Note that $k''(q^{-1}) = q^2(p-1) \ne 0$. Hence, we must have:
$$(q-1)^3=q-1 \implies q=1,2.$$ 
This completes the proof of $iv$. and the proof of Proposition \ref{propep1}.



    
\end{proof}
\begin{lem}\label{multf1}
    We have the following:
    \begin{enumerate}[i.]
        \item The second derivative of $f_{\beta,h}$ has at most two roots (counting multiplicity) \footnote{A rational function $g(x) := p(x)/q(x)$ where $p$ and $q$ are polynomials, is said to have a root $r$ of multiplicity $k$, if $r$ is a root of $p$ having multiplicity $k$, and $q(r) \ne 0$.} in $(0,1]$.
        \item The second derivative of $f_{\beta,h}$ can have the root $0$ with multiplicity at most four, and this multiplicity is exactly four if and only if $(\beta,h)=(\frac{2}{3},0)$ and $(p,q)=(4,2)$, in which case, $(\frac{1}{2},\frac{1}{2})$ is the unique global maximizer of $H_{\beta,h}$.
         \item Suppose that $0$ is a maximizer of $f_{\beta,h}$, such that $f_{\beta,h}''(0) = 0$. Then, $p \in \{2,3,4\}$, $q = 2$ and $(\beta,h) = (\frac{2^{p-1}}{p(p-1)},0)$. In this case, $f_{\beta,h}$ has $0$ as the unique maximizer.
         \item Suppose that $s>0$ is a maximizer of $f_{\beta,0}$. Then, $f_{\beta,0}''(s)<0$.
    \end{enumerate}
\end{lem}
\begin{proof}
i.~ To start with, note that:
    \begin{equation*}
        \begin{aligned}
             & f_{\beta,h}''(s)                                                                                    \\
             & =\frac{(q-1)^2}{q^2}k^\dprime\left(\frac{1+(q-1)s}{q}\right)+\frac{(q-1)}{q^2}k^\dprime\left(\frac{1-s}{q}\right)     \\
             & =\frac{(q-1)^2\beta p(p-1)}{q^p}(1+(q-1)s)^{p-2}+\frac{(q-1)\beta p(p-1)}{q^p}(1-s)^{p-2}-\frac{q-1}{(1-s)(1+(q-1)s)} \\
             & = p(s)-r(s),
        \end{aligned}
    \end{equation*}
    where $p(s)=\frac{(q-1)^2\beta p(p-1)}{q^p}(1+(q-1)s)^{p-2}+\frac{(q-1)\beta p(p-1)}{q^p}(1-s)^{p-2}$ and $r(s)=\frac{q-1}{(1-s)(1+(q-1)s)}$.
    Now we will prove that for $q\geq 3$, $p(s)$ is a polynomial with positive coefficients such that the coefficients increase and then decrease.
    Let $T_r$ be the coefficient of $s^r$ in $p(s)$.
    Taking $m := p-2$ and $a := q-1$ we get:
    \begin{equation*}
        \frac{T_{r+1}}{T_r}=\frac{m-r}{r+1}\frac{a^{r+2}+(-1)^{r+1}}{a^{r+1}+(-1)^r}.
    \end{equation*}
    Now, $T_r$ is clearly increasing  for $r\leq \frac{m}{2}$. Now, $\frac{T_{r+1}}{T_r}\geq 1$ if and only if,
    \begin{equation*}
        r\leq \frac{ma-x_r}{a+x_r}.
    \end{equation*}
    Let us look at the difference, $\left|\frac{ma-x_r}{a+x_r}-\frac{ma-1}{a+1}\right|$.
    \begin{equation*}
        \begin{aligned}
            \left|\frac{ma-x_r}{a+x_r}-\frac{ma-1}{a+1}\right|&=\frac{a(m+1)\left|1-x_r\right|}{(a+1)(a+x_r)}\\
            &\leq \frac{a(m+1)}{a+x_r}a^{-r-2},
        \end{aligned}
    \end{equation*}
    where $x_r=\frac{1+(-1)^ra^{-r-1}}{1+(-1)^{r+1}a^{-r-2}}$. 
    Since, $\frac{a}{a+x_r}\leq 1$ then for $q\geq 3$,
    \begin{equation*}
        \left|\frac{ma-x_r}{a+x_r}-\frac{ma-1}{a+1}\right|
            \leq \frac{m+1}{2^{r+2}-1}.
    \end{equation*}
    Also, if $r>\frac{m}{2}$ then, $\left|\frac{ma-x_r}{a+x_r}-\frac{ma-1}{a+1}\right|<(m+1)/(2^{\frac{m}{2}+2}-1)$. Moreover, $(m+1)/(2^{\frac{m}{2}+2}-1)< \frac{9}{20}$. If $r\leq \frac{ma-1}{a+1}-\frac{9}{20}$ then $\frac{T_{r+1}}{T_r}\geq 1$. On the other hand, if $r\geq \frac{ma-1}{a+1}+\frac{9}{20}$ then $\frac{T_{r+1}}{T_r}\leq 1$. Now, there can exist at most one integer in between $ \frac{ma-1}{a+1}-\frac{9}{20}$ and $ \frac{ma-1}{a+1}+\frac{9}{20}$. If $i \in [\frac{ma-1}{a+1}-\frac{9}{20},\frac{ma-1}{a+1}+\frac{9}{20}]$ is an integer, then $T_0\leq T_1\leq \ldots \leq T_i$ and $T_{i+1}\geq \ldots \ge T_m$. The same strings of inequalities are true with $i := \lfloor \frac{ma-1}{a+1}-\frac{9}{20} \rfloor$ if the interval $[\frac{ma-1}{a+1}-\frac{9}{20},\frac{ma-1}{a+1}+\frac{9}{20}]$ does not contain an integer. On the other hand if there exists no integer in $[\frac{ma-1}{a+1}-\frac{2}{5},\frac{ma-1}{a+1}+\frac{2}{5}]$, then for $i=\flr{\frac{ma-1}{a+1}-\frac{2}{5}}$, we have $T_0\leq T_1\leq \ldots \leq T_i$ and $T_{i+1}\geq \ldots T_m$.
    This proves that for $q\geq 3$, the coefficients of $p(s)$ increase and then decrease.
    
    Now, consider the polynomial $(1+(q-1)s)p(s)$. The coefficient of $s^r$ of this polynomial is $c_r=(q-1)T_{r-1}+T_r$. Let $i$ be the integer where $T_i$ attains maximum. Then for $r+1\leq i$, $c_{r+1}\geq c_r$ and for $r-1 \geq i$, we have $c_{r+1}\leq c_r$. This shows that, $c_0\leq \ldots \le c_i$ and $c_{i+1}\geq \ldots \ge c_{p-1}$. So, the coefficients of the polynomial $(1+(q-1)s)p(s)$ increase and then decrease. Let,
    $$
        \begin{aligned}
                     & (1+(q-1)s)p(s)=\sum_{r=0}^{p-1}c_rs^r               \\
            \implies & (1-s)(1+(q-1)s)p(s)=\sum_{r=0}^{p}(c_r-c_{r-1})s^r,
        \end{aligned}
    $$
    where $c_{-1}=c_{p}=0$.
    Hence, the coefficients of $(1-s)(1+(q-1)s)p(s)$ has at most one sign change. which implies that the coefficients of $(1-s)(1+(q-1)s)p(s)-(q-1)$ can have at most two sign changes. Hence, by Descartes' rule of signs there are at most two positive roots. Hence, $f^\dprime_{\beta,h}$ can have at most two zeroes in $(0,1]$. 
    \par
    Now, consider the case $q=2$, whence $c_r=T_{r-1}+T_r$. This implies that $c_r-c_{r-1}=T_r-T_{r-2}$, which is $0$ if $r$ is odd, and for $r$ even, equals $\alpha [\binom{p-2}{r} - \binom{p-2}{r-2}]$ for some positive constant $\alpha$, which is non-negative for the first few even values of $r$, and then becomes non-positive for the remaining even values of $r$. Hence, the coefficients of $(1-s)(1+(q-1)s)p(s)$ can have at most one sign change, and the rest of the argument follows exactly as before.
    \vspace{0.1in}
    
\noindent ii.~   Suppose that $f_{\beta,h}''$ has the root $0$ with multiplicity at least $4$. Then, we must have $f^\dprime_{\beta,h}(0)=f^{(3)}_{\beta,h}(0)=f^{(4)}_{\beta,h}(0)=f^{(5)}_{\beta,h}(0)=0$. Now, $f_{\beta,h}''(0) = 0$ implies that
    $$\beta=\frac{q^{p-1}}{p(p-1)}, \text{ and}$$
    $$\bm x_s=\left(\frac{1}{(\beta p(p-1))^{1 / p-1}},\frac{1}{(\beta p(p-1))^{1 / p-1}},\ldots,\frac{1}{(\beta p(p-1))^{1 / p-1}}\right). $$
    Hence, $h=0$ by Proposition \ref{propep1} ii. Since, $f^{(3)}_{\beta, h}(0)=(q-1)(p-1)(q-2)$, we must have $q=2$. Furthermore, $f^{(4)}_{\beta,h}(0)=(p-4)(p-1)$, which immediately gives $p=4$. This also implies that $\beta = 2/3$. It is now easy to check that $f^{(5)}_{\beta,h}(0)=0$ and $f^{(6)}_{\beta,h}(0) = -24 <0$. If $f_{\beta,h}''$ had the root $0$ with multiplicity at least $5$, then $f_{\beta,h}^{(6)}(0)$ would have been $0$, a contradiction. This implies that $f_{\beta,h}''$ has the root $0$ with multiplicity at most $4$, and in this case, $(\beta,h) = (\frac{2}{3},0)$ and $(p,q)=(4,2)$. Conversely, if $(\beta,h) = (\frac{2}{3},0)$ and $(p,q)=(4,2)$, then $f_{\beta,h}''(s) = \frac{s^4}{s^2-1}$, and hence, has the root $0$ with multiplicity exactly four. In this case, $f_{\beta,h}$ is concave on $[0,1]$ and strictly concave on $(0,1]$, hence, any stationary point of $f_{\beta,h}$ in $[0,1]$ must be its unique global maximizer. Clearly, $f_{\beta,h}'(0) = 0$, which now implies that $0$ is the unique global maximizer of $f_{\beta,h}$ and completes the proof of Lemma \ref{multf1} ii.

    \noindent iii.~By the higher derivative test, we must have $f_{\beta,h}^{(3)}(0) =0$ and $f_{\beta,h}^{(4)}(0)\le 0$. Note that $f_{\beta,h}''(0) = 0$ implies that $\beta p (p-1) = q^{p-1}$ and hence, all coordinates of $x_0$ equal $(\beta p (p-1))^{\frac{1}{1-p}}$. This immediately implies that $q=2$ and $h=0$, in view of Proposition \ref{propep1}. Hence, we have:
    $$f_{\beta,h}''(s) = \frac{1}{2}(1+s)^{p-2} + \frac{1}{2}(1-s)^{p-2} - \frac{1}{1-s^2}.$$ Hence, $f_{\beta,h}^{(4)}(0) = (p-1)(p-4)$. Since, $f_{\beta,h}^{(4)}(0) \le 0$, we must have $p\in \{2,3,4\}$. In all these cases,
    $f_{\beta,h}''(s) = - \frac{s^{2\lfloor p/2\rfloor}}{1-s^2} <0$ for all $s\in (0,1]$, which implies that $f_{\beta,h}'$ is strictly decreasing on $[0,1]$. Since $f_{\beta,h}'(0) = 0$, $f_{\beta,h}'$ must be negative on $(0,1]$, hence $f_{\beta,h}$ must be strictly decreasing on $[0,1]$. This completes the proof of Lemma \ref{multf1} iii.

    \noindent iv. Suppose that $s>0$ is a maximizer of $f_{\beta,0}$. If $f_{\beta,0}''(s)=0$, then by the higher derivative test, $f_{\beta,0}^{(3)}(s) =0$, and hence, $s$ is a root of $f_{\beta,0}''$ with multiplicity at least $2$. Since $f_{\beta,0}'(0) = 0$, $f_{\beta,h}''$ must have at least one root in $(0,s)$. This implies that $f_{\beta,h}''$ has at least three positive roots (counting multiplicity), contradicting part i. Hence, $f_{\beta,0}''(s) \ne 0$. Lemma \ref{multf1} iv. now follows from the derivative test.
\end{proof}

\begin{ppn}\label{spcl7}
    Let $\beta, h\geq 0$ and $\bm m_*= \bm x_s$ be some global maximizer of $H_{\beta,h}$.
    \begin{enumerate}[i.]
        \item If $f_{\beta,h}''(s)\ne 0$ then $\bm Q_{\bm m_*, \beta}$ is negative definite in $\cH_q$.
        \item If $f_{\beta,h}''(s) = 0$ then $\bm Q_{\bm m_*, \beta}$ is negative semi-definite in $\cH_q$ with kernel in $\cH_q$ equal to $\mathrm{Span}(\bm u)$ where,
              $$\bm u :=(1-q,1,\ldots,1).$$ 
    \end{enumerate}
\end{ppn}
\begin{proof}
    Now, $\bm m_*=\bm x_s$ for some $s \in [0,1]$.
    $$\bm Q_{\bm x_s,\beta}(\bm t) = \left(k^\dprime\left(\frac{1-s}{q}\right) +(q-1)\alpha(\bm t)k^\dprime\left(\frac{1+(q-1)s}{q}\right)\right) \sum_{r=2}^{q}t_r^2,$$
    where $\alpha(\bm t)=\dfrac{(\sum_{r=2}^{q}t_r)^2}{(q-1)\sum_{r=2}^{q}t_r^2}~.$ Since $0 \le \alpha(\bm t) \le 1$, the following two conditions imply negative definiteness of $\bm Q_{\bm x_s, \beta}$ on $\cH_q$:

    \begin{equation}
        k^\dprime\left(\frac{1-s}{q}\right)<0 \quad \text{and} \quad k^\dprime\left(\frac{1-s}{q}\right) +(q-1)k^\dprime\left(\frac{1+(q-1)s}{q}\right)<0
        \label{cond1}
    \end{equation}
    On the other hand, if $\bm Q_{\bm x_s,\beta}$ is negative definite on $\cH_q$, then setting $\bm t := (1,-1,0,0,\ldots,0)$ and $(\frac{1}{q-1}, \ldots, \frac{1}{q-1})$, one arrives at \eqref{cond1}. Hence, \eqref{cond1} is equivalent to the negative-definiteness of $\bm Q_{\bm x_s,\beta}$ on $\cH_q$. Now, if $f_{\beta,h}''(s) \ne 0$, then since $s$ is a maximizer of $f_{\beta,h}$, one must have $f_{\beta,h}''(s) < 0$. Since $$f_{\beta,h}''(s) = \frac{q-1}{q^2}\left[k^\dprime\left(\frac{1-s}{q}\right) +(q-1)k^\dprime\left(\frac{1+(q-1)s}{q}\right)\right],$$ we conclude that $f_{\beta,h}''(s) \ne 0$ implies the second condition in \eqref{cond1}. It also follows from the proof of Proposition \ref{propep1} i., that if $h>0$, then $k''(\frac{1-s}{q}) < 0$, so we may assume $h=0$. Moreover, if $\frac{1-s}{q} < (\beta p(p-1))^{1/(1-p)}$, then  $k''(\frac{1-s}{q}) < 0$, so by Proposition \eqref{propep1} i., it suffices to assume that $\frac{1-s}{q} = (\beta p(p-1))^{1/(1-p)}$. It now follows from the proof of Proposition \ref{propep1} i. (for the case $h=0$) that $\frac{1-s}{q}$ must be the unique global minimizer of the function $g_\beta$ defined in that proof, and hence, $\bm x_s$ must have all entries equal to $q^{-1}$, so in particular, $s=0$. Therefore, $q=(\beta p(p-1))^{1/(p-1)}$. However, this implies that $f_{\beta,h}''(s) = \frac{q-1}{q} k''(\frac{1}{q}) = 0$, since $k''(q^{-1}) = \beta p(p-1)q^{2-p} - q = 0$. This is a contradiction, thereby completing the proof of part i.

    Since $f_{\beta,h}''(s) = Q_{\bm x_s,\beta}(q^{-1}\bm u)$, we conclude that if $f_{\beta,h}''(s)=0$, then $\mathrm{Span}(\bm u) \subseteq \mathrm{Ker}(Q_{\bm x_s,\beta}) \bigcap \cH_q$. On the other hand, if $\boldsymbol{0} \ne \bm t \in \mathrm{Ker}(Q_{\bm x_s,\beta}) \bigcap \cH_q$, then since $\sum_{r=2}^q t_r^2 > 0$, we must have:
    $$k^\dprime\left(\frac{1-s}{q}\right) +(q-1)\alpha(\bm t)k^\dprime\left(\frac{1+(q-1)s}{q}\right) = 0.$$ Since $f_{\beta,h}''(s) = 0$, we also have:
    $$k^\dprime\left(\frac{1-s}{q}\right) +(q-1)k^\dprime\left(\frac{1+(q-1)s}{q}\right) = 0.$$ Therefore, if $\alpha(\bm t) \ne 1$, then $k^\dprime\left(\frac{1-s}{q}\right) = k^\dprime\left(\frac{1+(q-1)s}{q}\right) = 0$, implying that $\bm Q_{\bm x_s,\beta} \equiv 0$ on $\cH_q$. This implies that $\bm x_s$ must be the constant vector with all entries equal to $(\beta p (p-1))^\frac{1}{1-p} = q^{-1} \implies s=0$ and from Proposition \ref{propep1} ii. and iv., it follows that $h=0$ and $q=2$. 
    Note that for $q=2$, $\cH_q = \mathrm{Span}(\bm u)$. Finally, $\alpha(\bm t) = 1$ implies that $\bm t \in \mathrm{Span}(\bm u)$. This completes the proof of Proposition \ref{spcl7}. 
    \end{proof}

\begin{cor}\label{spder0}
    Let $\bm x_s$ be a global maximizer of $H_{\beta,h}$. If $(\beta,h)$ is a special point, then $f_{\beta,h}''(s) = 0$. On the other hand, if $(\beta,h)$ is a regular or critical point then $f_{\beta,h}''(s) < 0$ .
\end{cor}
    \begin{proof}
   If $(\beta,h)$ is a special point, then by definition, $\bm Q_{\bm x_s,\beta}$ is singular on $\cH_q$, so by Proposition \ref{spcl7} i., one must have $f_{\beta,h}''(s) = 0$. On the other hand, if $(\beta,h)$ is regular or critical, then $\bm Q_{\bm x_s,\beta}$ is negative definite on $\cH_q$. So, by Proposition \ref{spcl7} ii., $f_{\beta,h}''(s) \ne 0$. Since $s$ is a maximizer of $f_{\beta,h}$, we must thus have $f_{\beta,h}''(s) <0$.
    \end{proof}
    

\begin{lem}\label{f2sp}
    $\cS^2_{p,q}$ is non-empty if and only if $p=4$ and $q=2$. Moreover, $\cS^2_{4,2}=\{(2/3,0)\}$.
\end{lem}
\begin{proof}
    $\cS^2_{p,q}$ is non-empty if and only if there exists a special point $(\beta,h)$ satisfying $f^{(4)}_{\beta,h}(s)=0$, where $s$ is the unique maximizer of $f_{\beta,h}(s)$. Since $s$ is a maximizer of $f_{\beta,h}$, we must have $f_{\beta,h}'(s)=0$. Since $(\beta,h)$ is a special point, by Corollary \ref{spder0}, we have $f_{\beta,h}''(s) = 0$.
    Since, $s$ is a maximum of $f_{\beta,h}$, it now follows from the higher derivative test, that $f^{(3)}_{\beta,h}(s)=f^{(5)}_{\beta,h}(s)=0$. Denoting $f_{\beta,h}'' = P/Q$ for polynomials $P$ and $Q$ satisfying $Q(s) \ne 0$, we thus obtain that $P(s) = P'(s)=P''(s)=P^{(3)}(s) = 0$. If $d$ is the degree of $P$, then we have by Taylor expansion of $P$ around $s$:
    $$P(x) = \sum_{k=4}^{d} \frac{(x-s)^k}{k!} P^{(k)}(s) = (x-s)^4 \sum_{k=0}^{d-4} \frac{(x-s)^{k-4}}{k!} P^{(k)}(s).$$
    Thus, $s$ is a root of $f_{\beta,h}''$ of multiplicity at least $4$, and hence, by Lemma \ref{multf1}, one must have $s=0$, $(\beta,h)=(\frac{2}{3},0)$ and $(p,q)=(4,2)$. On the other hand, if $(\beta,h)=(\frac{2}{3},0)$ and $(p,q)=(4,2)$, then $f_{\beta,h}''(s) = \frac{s^4}{s^2-1}$, and hence, $f_{\beta,h}^{(k)}(s) = 0$ for $2\le k\le 5$. Further, by Lemma \ref{multf1}, $\bm x_0 := (\frac{1}{2}, \frac{1}{2})$ is the unique global maximizer of $H_{\beta,h}$. This proves that $(\beta,h)\in \cS_{p,q}^2$, thereby completing the proof of Lemma \ref{f2sp}. 
\end{proof}

\begin{lem}\label{partn}
    The sets $\cR_{p,q}$, $\cC_{p,q}$ and $\cS_{p,q}$ form a partition of the parameter space $\Theta$.
\end{lem}
\begin{proof}
 By definition, the sets $\cR_{p,q}$, $\cC_{p,q}$ and $\cS_{p,q}$ are disjoint. So, in order to prove Lemma \ref{partn}, it suffices to show that if $(\beta,h)\notin \cR_{p,q} \bigcup \cS_{p,q}$, then $(\beta,h)\in \cC_{p,q}$.  In case  $(\beta,h)\notin \cR_{p,q} \bigcup \cS_{p,q}$, the function $H_{\beta,h}$ has at least two different global maximizers. To begin with, assume that $H_{\beta,h}$ has (at least) two distinct global maximizers of the form $\bm x_s$ and $\bm x_t$ for some $0\le s < t$. Then, $s$ and $t$ are global maximizers of $f_{\beta,h}$. In this case, we will show that whenever $\bm x_s$ is a global maximizer of $H_{\beta,h}$, we must have $f_{\beta,h}''(s)\ne 0$. Suppose towards a contradiction, that $f_{\beta,h}''(s) = 0$ for some $s\in [0,1]$ such that $\bm x_s$ is a global maximizer of $H_{\beta,h}$ . By Lemma \ref{multf1} iii., we must have $s > 0$. Since $s$ is a maximizer of $f_{\beta,h}$, we must have $f_{\beta,h}^{(3)}(s) = 0$ by the higher derivatve test. Note that $f_{\beta,h}^{(4)}(s) \ne 0$, since otherwise, $s>0$ would be a root of $f_{\beta,h}^{(2)}$ with multiplicity at least $3$, which is impossible by Lemma \ref{multf1} i. Since $s$ is a maximizer of $f_{\beta,h}$, we must have $f_{\beta,h}^{(4)}(s)<0$, which implies that $f_{\beta,h}''$ attains strict local maximum at $s$. Further, $s$ is a root of $f_{\beta,h}''$ having multiplicity $2$, so by Lemma \ref{multf1} i.,  $f_{\beta,h}''$ cannot have any other root in $(0,1]$. It follows that $f_{\beta,h}''(x) <0$ for all $x\in (0,1]\setminus\{s\}$ and hence, $f_{\beta,h}'$ is strictly decreasing on $(0,1]$, implying that $f_{\beta,h}'$ can have at most one root on $(0,1]$, i.e. $f_{\beta,h}$ can have at most one maximizer on $(0,1]$, a contradiction! So, we must have $f_{\beta,h}''(s) \ne 0$ for all $s\in [0,1]$ such that  $\bm x_s$ is a global maximizer of $H_{\beta,h}$. 

 The only remaining case is when $H_{\beta,h}$ has multiple global maximizers, all of whom are permutations of one another. In this case, the maximizers are trivially not constant vectors, and hence, by Proposition \ref{propep1} ii., we must have $h=0$. Hence, the maximizers of $H_{\beta,h}$ are precisely all permutations of the vector $\bm x_s$ for some $s>0$. By Lemma \ref{multf1} iv., we must have $f_{\beta,h}''(s)<0$. Hence, we have shown that whenever $\bm x_s$ is a global maximizer of $H_{\beta,h}$, $f_{\beta,h}''(s) \ne 0$, and hence, $\bm Q_{\bm x_s,\beta}$ is negative-definite on $\cH_q$ (by Lemma \ref{spcl7}).
\end{proof}

\begin{lem}\label{specialtype2thm}
    $\cS_{p,q}$ is a singleton set for every $p\ge 2, q\ge 2$.
 \end{lem}
\begin{proof}
We claim that if $(\beta,h)\in \cS_{p,q}$,  then $\sup_{x \in [0,1]} f_{\beta,h}''(x) = 0$. Throughout this proof, we will denote the unique maximizer of $f_{\beta,h}$ by $s$. To show the claim, first suppose that $s>0$. Note that by Lemma \ref{multf1} and the higher derivative test, $f_{\beta,h}''$ must have the root $s$ with multiplicity exactly two. Hence, $f_{\beta,h}^{(4)}(s) \ne 0$, and $s$ being the maximizer of $f_{\beta,h}$, one must have $f_{\beta,h}^{(4)}(s) <0$. This immediately shows that $s$ is a local maximizer of $f_{\beta,h}''$, which cannot have any root in $(0,1]$ other than $s$. This forces $f_{\beta,h}''$ to be strictly negative on $(0,1]\setminus \{s\}$, thereby proving the claim for the case $s>0$. Now, if $s=0$, then by Lemma \ref{multf1} iii., we must have $p\in \{2,3,4\}, q=2$ and $(\beta,h) = (\frac{2^{p-1}}{p(p-1)},0)$. In all these cases, $f_{\beta,h}''(s) = -\frac{s^{2\lfloor p/2\rfloor}}{1-s^2} < 0$ on $(0,1]$, completing the proof of the claim.

Now, it is easy to see that the function $w(\beta) := \sup_{x\in [0,1]} f_{\beta,0}''(x)$ is strictly increasing and continuous in $\beta$, with $w(0) \le q^{-1}-1<0$ and $w(\infty) = \infty$. Hence, there exists a unique $\widetilde{\beta}_{p,q}$ such that $w(\wbt) = 0$. In fact, $\wbt$ is given by:
\begin{equation}\label{def_beta_spec}
   \wbt = \inf\{\beta \ge 0: \sup_{x\in [0,1]} f_{\beta,0}''(x) > 0\}. 
\end{equation}
 By the previous paragraph, we have thus shown that if $(\beta,h)\in \cS_{p,q}$, then $\beta = \wbt$. By Lemma \ref{multf1}, $f_{\wbt,0}''$ can have at most three distinct roots in $[0,1]$, and define $s_{p,q}$ to be the largest root. We claim that $s_{p,q}=s$. To see this, note that if $s=0$, then $f_{\wbt,0}''$ is negative on $(0,1]$, and hence, $0$ is its only root, so trivially $s_{p,q}=s=0$. On the other hand, if $s>0$, then $f_{\wbt,0}''$ cannot have any positive root other than $s$ (since the root $s$ has multiplicity two), so once again, $s_{p,q}=s$. Since $f_{\wbt,h}'(s_{p,q}) = 0$, one must have:
 \begin{equation*}
     h = \wht := k_{\wbt,p}'\left(\frac{1-s_{p,q}}{q}\right) - k_{\wbt,p}'\left(\frac{1+(q-1)s_{p,q}}{q}\right)
 \end{equation*}
We have thus shown that $\cS_{p,q}$ has at most one element. We will now show that for every $p\ge 2$, $q\ge 2$, $(\wbt,\wht)\in \cS_{p,q}$. To see this, we will first show that $\wht \ge 0$, for which it suffices to show that the function $v:[0,1]\mapsto \R$ defined as:
$$v(x) := k_{\wbt,p}'\left(\frac{1-x}{q}\right) - k_{\wbt,p}'\left(\frac{1+(q-1)x}{q}\right)$$
is non-decreasing. Note that $v'(x) = \frac{q}{1-q} f_{\wbt,0}''(x) \ge 0$, since $\sup_{x\in [0,1]} f_{\wbt,0}''(x) = 0$. This proves our claim that $\wht \ge 0$. Next, note that by construction, we have:
$$f_{\wbt,\wht}'(s_{p,q}) = f_{\wbt,\wht}''(s_{p,q}) = 0.$$ Also, $f_{\wbt,\wht}''(x) \le 0$ for all $x\in [0,1]$, which implies that $s_{p,q}$ is a global maximizer of $f$. This completes the proof of Lemma \ref{specialtype2thm} in view of Corollary \ref{spder0}.
\end{proof}

\begin{lem}\label{crit_less_spec}
    If $f_{\beta,h}$ has more than one global maximizer in $[0,1]$, then $(\beta,h) \in (\tilde{\beta}_{p,q},\infty)\times [0,\wht)$. 
\end{lem}
\begin{proof}
    If $\beta<\tilde{\beta}_{p,q}$ then it follows from \eqref{def_beta_spec} and the fact that $\sup_{x\in [0,1]} f^\dprime_{\beta,h}(x)$ is strictly increasing in $\beta$, that $\sup_{x\in [0,1]} \f^\dprime(x) < 0$. Hence, $\f$ is strictly concave, so cannot have more than one maximizer. If $\beta = \wbt$, then also, $f_{\beta,h}$ is concave. If it has two distinct global maximizers $s<t$, then $f_{\beta,h}''$ must vanish on the interval $(s,t)$, contradicting Lemma \ref{multf1} i. So, one must have $\beta > \wbt$. Consequently, $f_{\beta,h}''(s_{p,q}) > 0$, where $s_{p,q}$ is the unique global maximizer of $f_{\wbt,\wht}$. Since $\lim_{x\rightarrow 1^-} f_{\beta,h}''(x) = -\infty$, the rational function $f_{\beta,h}''$ must have an odd number of roots (counting multiplicity) in $[s_{p,q},1]$. By Lemma \ref{multf1} i., $f_{\beta,h}''$ has a unique root in $[s_{p,q},1]$.

    Now, suppose that $h \ge \wht$. Then, since $f_{\wbt,\wht}'$ is non-negative on $[0,s_{p,q}]$, one must have $f_{\beta,h}' > 0$ on $(0,s_{p,q}]$. This, in particular, implies that $0$ cannot be a local maximizer of $f_{\beta,h}$. Moreover, $f_{\beta,h}'$ cannot have more than two distinct roots on $[s_{p,q},1]$, since $f_{\beta,h}''$ has a unique root in this interval. Since $f_{\beta,h}$ has at least two distinct global maximizers, $f_{\beta,h}'$ has exactly two distinct roots on $(s_{p,q},1)$, both these roots being global maximizers of $f_{\beta,h}$. Call these two roots $s_1<s_2$. Clearly, $f_{\beta,h}'$ is either stricty positive on $(s_1,s_2)$, or strictly negative on $(s_1,s_2)$, since it cannot have any root in this interval. In the first case, $s_1$ cannot be a local maximizer of $f_{\beta,h}$, and in the second case, $s_2$ cannot be a local maximizer of $f_{\beta,h}$, a contradiction! 
    Hence, $h < \wht$, thereby completing the proof of Lemma \ref{crit_less_spec}.  
\end{proof}

Let $(s_1,s_2) \in U :=\{(x,y):x<y \text{ where } x,y \in [0,1)\}$. Now, if $s_1,s_2$ are the stationary points of $f_{\beta,h}$ then,
$$\begin{aligned}
         & \frac{\beta p}{q^{p-1}}\left(1+(q-1) s_{1}\right)^{p-1}-\frac{\beta p}{q^{p-1}}\left(1-s_1\right)^{p-1}+h=\ln \left(1+(q-1) s_{1}\right)-\ln \left(1-s_{1}\right) \\
         & \frac{\beta p}{q^{p-1}}\left(1+(q-1) s_{2}\right)^{p-1}-\frac{\beta p}{q^{p-1}}\left(1-s_{2}\right)^{p-1}+h=\ln \left(1+(q-1) s_{2}\right)-\ln \left(1-s_{2}\right) \\
    \end{aligned}$$
In other words,
$$
    \left[\begin{array}{ll}
            \frac{p}{q^{p-1}}\left(1+(q-1) s_{1}\right)^{p-1}-\frac{p}{q^{p-1}}\left(1-s_1\right)^{p-1} & 1 \\
            \frac{p}{q^{p-1}}\left(1+(q-1) s_{2}\right)^{p-1}-\frac{p}{q^{p-1}}\left(1-s_{2}\right)^{p-1} & 1
        \end{array}\right]\left[\begin{array}{l}
            \beta \\
            h
        \end{array}\right]=\left[\begin{array}{l}
            \ln \left(1+(q-1) s_{1}\right)-\ln \left(1-s_{1}\right) \\
            \ln \left(1+(q-1) s_{2}\right)-\ln \left(1-s_{2}\right)
        \end{array}\right]
$$
Clearly the matrix is invertible since $a(x)=\frac{p}{q^{p-1}}\left(1+(q-1) x\right)^{p-1}-\frac{p}{q^{p-1}}\left(1-x\right)^{p-1}$ is a strictly increasing function.
Define $G: U \to \R^2$, 
\begin{equation*}
    \begin{bmatrix}
       x \\
       y
    \end{bmatrix}
    \mapsto
    {\left[\begin{array}{ll}
                \frac{p}{q^{p-1}}\left(1+(q-1) x\right)^{p-1}-\frac{p}{q^{p-1}}\left(1-x\right)^{p-1}   & 1 \\
                \frac{p}{q^{p-1}}\left(1+(q-1) y\right)^{p-1}-\frac{p}{q^{p-1}}\left(1-y\right)^{p-1} & 1
            \end{array}\right]^{-1}\left[\begin{array}{l}
                \ln \left(1+(q-1) x\right)-\ln \left(1-x\right) \\
                \ln \left(1+(q-1) y\right)-\ln \left(1-y\right)
            \end{array}\right]}
\end{equation*}
Clearly, $G(U)$ is the set of all $(\beta,h)$ such that $f_{\beta,h}$ has more than one stationary point, and hence, $\cC_{p,q}^1 \subset G(U)$.

\begin{lem}[Properties of strongly critical points]\label{stcrprf}
    \begin{enumerate}[i.]
        \item For any $h\geq 0$ there exists at most one $\beta$ such that $f_{\beta,h}$ has more than one global maximizer.
        \item For any $(\beta_1,h_1),(\beta_2,h_2) \in \cC^1_{p,q}$ such that $h_2>h_1$, we must have $\beta_2<\beta_1$.
        \item $\cC^1_{p,q}$ is a compact set in $G(U)$.
    \end{enumerate}
\end{lem}
\begin{proof}

    \begin{enumerate}[i.]
        \item We will first show that $f_{\beta,h}$ cannot have more than two global maximizers, and in case it has exactly two global maximizers $s<t$, there must exist $u\in (s,t)$ such that $f_{\beta,h}'$ is positive on $[0,s)$, negative on $(s,u)$, positive on $(u,t)$ and negative on $(t,1)$. To see this, let $s$ and $t$ denote the smallest and second smallest global maximizers of $f_{\beta,h}$. Since $f_{\beta,h}'$ has a finite number of roots, it must be negative on some non-empty right neighborhood of $s$ and positive on some non-empty left neighborhood of $t$. Hence, it must vanish at some $u \in (s,t)$. Since $f_{\beta,h}''$ cannot have more than two positive roots, $s,u$ and $t$ must be the only three roots of $f_{\beta,h}'$. Hence, $f_{\beta,h}'$ cannot change sign on each of the intervals $[0,s), (s,u), (u,t)$ and $(t,1]$. Since $s$ and $t$ are global maximizers of $f_{\beta,h}$, $f_{\beta,h}'$ must be positive on some non-empty left neighborhood of $s$ and negative on some non-empty right neighborhood of $t$, too. This proves our claim.
        
        Now, suppose that for some $h \ge 0$, there exist $\beta_1<\beta_2$ such that $f_{\beta_1,h}$ and $f_{\beta_2,h}$ have multiple maximizers. Let $s_1<t_1$ and $s_2<t_2$ be the respective global maximizers. Since $f_{\beta_2,h}' > f_{\beta_1,h}'$ on $(0,1]$, we must have $f_{\beta_2,h}' > 0$ on $[0,s_1]$ and $[u_1,t_1]$. This, coupled with the fact that $f_{\beta_2,h}'(1) = -\infty$, implies that $s_2>s_1$ and $t_2>t_1$. Since $f_{\beta_2,h}''$ already has two roots larger than $s_2$ (one in $(s_2,u_2)$ and the other in $(u_2,t_2)$), it must not change sign in $[0,s_2]$. This sign cannot be positive, because then $f_{\beta_2,h}'$ would be stricty increasing on $[0,s_2]$, which would contradict the fact that $f_{\beta_2,h}' >0$ on $[0,s_2)$ and $0$ at $s_2$. Hence, $f_{\beta_2,h}'' < 0$ on $[0,s_2)$. This shows that $f_{\beta_1,h}'' < 0$ on $[0,s_2)$. Now, suppose that $s_2\ge t_1$. Then, $f_{\beta_1,h}''<0$ on $[0,t_1)$, i.e. $f_{\beta_1,h}'$ is strictly decreasing on $[0,t_1)$, contradicting that it has two roots, $s_1$ and $u_1$ in this interval! Hence, one must have $s_2<t_1$, i.e. $s_1<s_2<t_1<t_2$. 
        
       Now, it is easy to see that for $x>y$,
        \begin{equation}\label{f10y}
            f_{\beta_2,h}(x)-f_{\beta_1,h}(x)>f_{\beta_2,h}(y)-f_{\beta_1,h}(y).
        \end{equation}
        Hence, $f_{\beta_2,h}(s_2)-f_{\beta_1,h}(s_2) < f_{\beta_2,h}(t_1)-f_{\beta_1,h}(t_1) \le f_{\beta_2,h}(s_2) - f_{\beta_1,h}(t_1)$, which implies that $f_{\beta_1,h}(t_1) < f_{\beta_1,h}(s_2)$, contradicting that $t_1$ is a global maximizer of $f_{\beta_1,h}$. The proof of part i. is now complete. 
    
        \item Suppose that $(\beta_1,h_1)$ and $(\beta_2,h_2)$ are two points in $\mathcal{C}_{p,q}^1$ such that $h_2>h_1$ and $\beta_2\ge \beta_1$. Suppose that $s_i<t_i$ are the global maximizers $f_{\beta_i,h_i}$ for $i \in \{1,2\}$. The rest of the proof of part ii. proceeds exactly similarly as the proof of part i, so we just highlight the main steps. By similar logic as in part i, we can argue that $s_1<s_2<t_1<t_2$. Now, we observe that for $x>y$,
        $$f_{\beta_2,h_2}(x) - f_{\beta_1,h_1}(x) > f_{\beta_2,h_2}(y) - f_{\beta_1,h_1}(y)$$ and hence, we have:
        $$f_{\beta_2,h_2}(s_2) - f_{\beta_1,h_1}(s_2) < f_{\beta_2,h_2}(t_1) - f_{\beta_1,h_1}(t_1) \le f_{\beta_2,h_2}(s_2) - f_{\beta_1,h_1}(t_1)$$
        which implies that $f_{\beta_1,h_1}(t_1) < f_{\beta_1,h_1}(s_2)$, contradicting that $t_1$ is a global maximizer of $f_{\beta_1,h_1}$. This completes the proof of part ii.
          
        \item For $(\beta,h) \in \Theta$, let $F(\beta,h,x,y) := f_{\beta,h}(x) - f_{\beta,h}(y)$. Define $\cF: U \to \R$ as:
              $$
                  \cF(s_1,s_2)=F(G(s_1,s_2),s_1,s_2).
              $$
              We claim that $G\left(\mathcal{\cF}^{-1}(\{0\})\right)=\cC_{p, q}^1$. To see this, first note that if $(\beta,h)\in \cC_{p,q}^1$, then the two global maximizers $s_1 < s_2$ of $f_{\beta,h}$ satisfy the stationary equation
              $(\beta,h) = G(s_1,s_2)$, and hence, $\cF(s_1,s_2) = f_{\beta,h}(s_1) - f_{\beta,h}(s_2)=0$. So, $(s_1,s_2) \in \cF^{-1}(\{0\})$ and hence, $(\beta,h) = G(s_1,s_2) \in G(\cF^{-1}(\{0\})$. Conversely, if $(\beta,h) \in G(\cF^{-1}(\{0\})$, then $(\beta,h) = G(s_1,s_2)$ for some $(s_1,s_2)\in U$ satisfying $f_{\beta,h}(s_1) = f_{\beta,h}(s_2)$. Clearly, $s_1,s_2$ are two stationary points of $f_{\beta,h}$, and by Rolle's theorem, $f_{\beta,h}$ has another stationary point $s_3 \in (s_1,s_2)$. By Lemma \ref{multf1} i., $s_1,s_2$ and $s_3$ are the only stationary points of $f_{\beta,h}$. Hence, $f_{\beta,h}'$ must be negative on $(s_2,1]$, since it diverges to $-\infty$ near $1$. If $f_{\beta,h}'$ were negative on $(s_3,s_2)$, then $s_2$ would have been a local maximizer of $f_{\beta,h}'$, which would imply that $f_{\beta,h}'' (s_2) =0$. This would however contradict Lemma \ref{multf1} i., since $f_{\beta,h}''$ has two other roots, one in $(s_3,s_2)$ and the other in $(s_1,s_3)$. Hence, $f_{\beta,h}'$ must be positive on $(s_3,s_2)$. Thus, $s_2$ is a local maximizer of $f_{\beta,h}$. Next, if $f_{\beta,h}'$ were positive on $(s_1,s_3)$, then once again, $f_{\beta,h}''(s_3) = 0$, contradicting that $f_{\beta,h}''$ must have two other roots, one in $(s_1,s_3)$ and the other in $(s_3,s_2)$. Hence, $f_{\beta,h}'$ must be negative on $(s_1,s_3)$. Hence, $s_3$ must be a local minimizer of $f_{\beta,h}$. Finally, by an exactly similar argument, we can derive that $f_{\beta,h}'$ must be positive on $[0,s_1)$. This once again shows that $s_1$ is a local maximizer of $f_{\beta,h}$. Now, it follows from the sign-changing structure of $f_{\beta,h}'$, that $s_1$ and $s_2$ are global maximizers of $f_{\beta,h}$ on the intervals $[0,s_3]$ and $[s_3,1]$, respectively. Since  $f_{\beta,h}(s_1) = f_{\beta,h}(s_2)$, they must be global maximizers of $f_{\beta,h}$ on $[0,1]$, thereby implying that $(\beta,h)\in \cC_{p,q}^1$. This proves our claim, that $G\left(\mathcal{\cF}^{-1}(\{0\})\right)=\cC_{p, q}^1$.

              Now, $\cF$ being a continuous function,  $\mathcal{\cF}^{-1}(\{0\})$ must be closed. Also, it is bounded, since $\mathcal{\cF}^{-1}(\{0\}) \subset [0,1]^2$. Hence, $\mathcal{\cF}^{-1}(\{0\})$ is compact. $G$ being a continuous function on $U$, $G(\mathcal{\cF}^{-1}(\{0\}))$ must be compact in $G(U)$. This proves iii. and completes the proof of Lemma \ref{stcrprf} . 
    \end{enumerate}
\end{proof}

\begin{lem}[Properties of Critical Points]
    \label{crit_part}
    For $(\beta,h)\in \Theta$, we have the following:
    \begin{enumerate}[i.]
        \item If $(\beta,h)\in \cC_{p,q}$ for some $h>0$, then $(\beta,h) \in \cC_{p,q}^1$, and $H_{\beta,h}$ has exactly two global maximizers. Moreover, these maximizers are not permutations of one another.
        
        \item If $h=0$, then there exists $\beta_c = \beta_c(p,q) \in (0,\infty)$ satisfying the following:
        \begin{enumerate}
            \item  If $\beta <\beta_c$, then $\bm x_0 = (\frac{1}{q}, \ldots, \frac{1}{q})$ is the unique global maximizer of $H_{\beta, h}$. Consequently, $(\beta,0) \in \cR_{p,q}$.

            \item If $\beta>\beta_c$, then there are exactly $q$ global maximizers of $H_{\beta, h}$, which are precisely all the $q$ possible permutations of $\bm x_s$ for some appropriate $s \in(0, 1)$. Consequently, $(\beta,0)\in \cC_{p.q}^2$.

            \item If $\beta = \beta_c$, then $\bm x_0 = (\frac{1}{q}, \ldots, \frac{1}{q})$ is a global maximizer of $H_{\beta,0}$. If $q\ne 2$ or $p \ge 5$, then the remaining global maximizers are precisely all the $q$ possible permutations of $\bm x_s$ for some appropriate $s \in(0, 1)$. Otherwise, i.e. if $p\in \{2,3,4\}$ and $q=2$, then $(\beta_c,0)$ is a special point.
        \end{enumerate}
    \end{enumerate}
\end{lem}
\begin{proof}
i.~Assume that $(\beta,h)\in \cC_{p,q}$ for some $h > 0$. It easily follows from Lemma \ref{propep1} ii., that any global maximizer of $H_{\beta,h}$ must be of the form $\bm x_u$ for some $u\in (0,1]$. If $(\beta,h)\in \cC_{p,q}^2$, then $H_{\beta,h}$ must have two distinct global maximizers $\bm x$ and $\bm y$ which are permutations of one another. By Lemma \ref{propep1}, $\bm x = (x_1,x_2,x_2,\ldots,x_2)$ for some $x_1>x_2$. Now, $\bm y$ being a permutation of $\bm x$ distinct from $\bm x$ itself, we must have $y_1 = x_2$ and $y_i = x_1$ for some $i\ge 2$. This implies that $y_1 < y_i$, contradicting Lemma \ref{propep1} ii. Hence, $(\beta,h)\in \cC_{p,q}^1$. So, $f_{\beta,h}$ has at least two global maximizers, and it follows from the arguments in the first paragraph of the proof of Lemma \ref{stcrprf} i., that $f_{\beta,h}$ has exactly two global maximizers $s$ and $t$. Then, $\bm x_s$ and $\bm x_t$ are two distinct global maximizers of $H_{\beta,h}$. If it has a third global maximizer, then this one must be of the form $\bm x_u$ for some $u \in (0,1]$, and hence, $u$ would be a global maximizer of $f_{\beta,h}$. This forces $u \in \{s,t\}$, and hence, $\bm x_u \in \{\bm x_s,\bm x_t\}$, a contradiction. Thus, $H_{\beta,h}$ has exactly two global maximizes $\bm x_s$ and $\bm x_t$. Since the first coordinate of each of these is strictly greater than the remaining $q-1$ coordinates which are all equal, they cannot be permutations of one another. This proves i.

\vspace{0.04in}
\noindent ii.~To begin with, define:
$$\beta_c = \beta_c(p,q) := \inf \left\{\beta \ge 0: f_{\beta,0}~\text{has a positive global maximizer}\right\}.$$
Note that:
$$f_{\beta,0}'(s) = \left(\frac{q-1}{q}\right) \left[\beta p\left(\left(\frac{1+(q-1)s}{q}\right)^{p-1} - \left(\frac{1-s}{q}\right)^{p-1}\right) + \log\left(\frac{1-s}{1+(q-1)s}\right)\right]~.$$
If we define $g:(0,1)\to \mathbb{R}$ as:
$$g(s) := \frac{1}{p}\cdot \frac{\log\left(\frac{1+(q-1)s}{1-s}\right)}{\left(\frac{1+(q-1)s}{q}\right)^{p-1} - \left(\frac{1-s}{q}\right)^{p-1}}~,$$
then $g$ is continuous and positive on $(0,1)$, with
$\lim_{s\rightarrow 0} g(s) = \frac{q^{p-1}}{p(p-1)} \in (0,\infty)$ and $\lim_{s\rightarrow 1} g(s) = +\infty$, which implies that $\inf_{s\in (0,1)} g(s) > 0$. Any $\beta \le \inf_{s\in (0,1)} g(s)$ must satisfy $f_{\beta,0}'(s) \le 0$ for all $s\in [0,1]$. Hence, $0$ is a global maximizer of $f_{\beta,0}$.
 The presence of any other global maximizer $t>0$ of $f_{\beta,0}$ would now imply that $f_{\beta,0}'\equiv 0$ on $[0,t]$, consequently $f_{\beta,0}''\equiv 0$ on $[0,t]$, thereby contradicting Lemma \ref{multf1}. Hence, $0$ is the unique global maximizer of $f_{\beta,0}$, thereby showing that $\beta_c \ge \inf_{s\in (0,1)} g(s) > 0$. Next, note that for $\beta > \sup_{s\in (0,\frac{1}{2}]}g(s)$, we have $f_{\beta,0}'(s) > 0$ for all $s\in (0,\frac{1}{2}]$, implying that any global maximizer of $f_{\beta,0}$ must be greater than or equal to $\frac{1}{2}$. This shows that $\beta_c <\infty$. Now, it trivially follows from the definition of $\beta_c$, that for $\beta<\beta_c$, $f_{\beta,0}$ has $0$ as the only global maximizer, implying that $\bm x_0$ is the only global maximizer of $H_{\beta,0}$. This proves (a).

 We will now show (b), for which it suffices to show that for $\beta>\beta_c$, $f_{\beta,0}$ has a unique global maximizer $s$, and $s>0$. Towards showing this, we first claim that if $\beta>\beta_c$, then $0$ cannot be a global maximizer of $f_{\beta,0}$. To see this, first choose $\beta_1\in (\beta_c,\beta)$ such that $f_{\beta_1,0}$ has a positive global maximizer $s_1$. Now, by \eqref{f10y}, we have
 $$f_{\beta,0}(s_1) - f_{\beta_1,0}(s_1) > f_{\beta,0}(0) - f_{\beta_1,0}(0),\quad \text{i.e.} \quad f_{\beta,0}(s_1) - f_{\beta,0}(0) > f_{\beta_1,0}(s_1) - f_{\beta_1,0}(0) \ge 0,$$ which immediately gives that $f_{\beta,0}(s_1) > f_{\beta,0}(0)$, thereby implying that $0$ cannot be a global maximizer of $f_{\beta,0}$. Suppose that $f_{\beta,0}$ has two distinct global maximizers $s>t>0$. Since $f_{\beta,0}'(0) = 0$, it follows from Lemma \ref{multf1} i., that $f_{\beta,0}$ cannot have any positive stationary point other than $s$ and $t$. So, $f_{\beta,0}'$ cannot change sign in each of the intervals $(0,t), (t,s)$ and $(s,1)$. Now, since $s$ is a global maximizer of $f_{\beta,0}$, the derivative $f_{\beta,0}'$ must be positive on some left neighborhood of $s$, and hence on $(t,s)$. However, since $t$ is also a global maximizer of $f_{\beta,0}$, the derivative $f_{\beta,0}'$ must be negative on some right neighborhood of $s$, and hence on $(t,s)$. This is clearly a contradiction! Hence $f_{\beta,0}$ has a unique global maximizer $s$, and $s>0$. The proof of (b) is now complete.

Finally, suppose that $\beta=\beta_c$. Note that for every $\beta_1<\beta$, $0$ is a global maximizer of $f_{\beta_1,0}$, i.e. $f_{\beta_1,0}(0) \ge f_{\beta_1,0}(x)$ for all $x\in [0,1]$. Taking limit as $\beta_1\uparrow \beta$, we have $f_{\beta,0}(0) \ge f_{\beta,0}(x)$ for all $x\in [0,1]$, thereby implying that $0$ is a global maximizer of $f_{\beta,0}$, i.e. $\bm x_0$ is a global maximizer of $H_{\beta,0}$. Now, note that if $f_{\beta,0}''(0)>0$, then $f_{\beta,0}' >0$ on some right neighborhood of $0$, which is not possible, since $0$ is a global maximizer of $f_{\beta,0}$. Hence, we must have $f_{\beta,0}''(0)\le 0$.

First, consider the case $q\ne 2$ or $p\ge 5$. Lemma \ref{multf1} iii. immediately gives us $f_{\beta,0}''(0)<0$. Now, take a sequence $\beta_n \downarrow \beta$, and let $s_n$ be the unique positive global maximizer of $f_{\beta_n,0}$. Since $\{s_n\}_{n\ge 1}$ is a bounded sequence, it has a convergent subsequence $\{s_{n_k}\}_{k\ge 1}$, converging to some $s\in [0,1]$. By uniform convergence of $f_{\beta_n,0}$ to $f_{\beta,0}$, we can conclude that $s$ is a global maximizer of $f_{\beta,0}$. Now, there exists $\varepsilon>0$ and $\delta>0$, such that $f_{\beta,0}''<-\delta$ on $[0,\varepsilon]$. By uniform convergence of $f_{\beta_{n_k},0}''$ to $f_{\beta,0}''$, we conclude that there exists $K\ge 1$ such that for all $k\ge K$, $f_{\beta_{n_k},0}''< -\frac{\delta}{2}$ on $[0,\varepsilon]$. Hence, $f_{\beta_{n_k},0}' < 0$ on $(0,\varepsilon]$, i.e. $f_{\beta_{n_k}}$ is strictly decreasing on $[0,\varepsilon]$ for all $k\ge K$. So, $s_{n,k} > \varepsilon$ for all $k\ge K$, thereby implying that $s\ge \varepsilon>0$. We have already shown that $f_{\beta,0}$ cannot have any positive global maximizer other than $s$. This shows that $H_{\beta,0}$ has precisely the global maximizers $\bm x_0$ and all permutations of $\bm x_s$. 

Finally, suppose that $p\in \{2,3,4\}$ and $q=2$, We will first show that $\beta_c = \frac{2^{p-1}}{p(p-1)}$ in this case. To see this, note that:
$$f_{\beta_1,0}''(s) = \frac{\beta_1 p(p-1)}{2^p}\left[(1+s)^{p-2} + (1-s)^{p-2}\right] - \frac{1}{1-s^2}\quad \implies \quad f_{\beta_1,0}''(0) = \frac{\beta_1 p(p-1)}{2^{p-1}} -1$$ and hence for $\beta_1 > \frac{2^{p-1}}{p(p-1)}$, $f_{\beta_1,0}''(0) > 0$. Since $f_{\beta_1,0}'(0) = 0$, this immediately implies that $f_{\beta_1,0}'$ is strictly positive on some right neighborhood of $0$, and hence it has some positive global maximizer. So, we must have $\beta_c\le \frac{2^{p-1}}{p(p-1)}$. Next, if $\beta_p := \frac{2^{p-1}}{p(p-1)}$, then $f_{\beta_p,0}''(s) = \frac{-s^{2\lfloor p/2\rfloor}}{1-s^2} < 0$ for all $s>0$, which, coupled with the fact that $f_{\beta_p,0}'(0) = 0$, implies that $f_{\beta_p,0}' <0$ on $(0,1]$, and hence, $f_{\beta_p,0}$ cannot have a positive global maximizer. This proves our claim that $\beta_c = \frac{2^{p-1}}{p(p-1)}$, and hence, $f_{\beta_c,0}' < 0$ on $(0,1]$. The last conclusion implies that $0$ is the unique global maximizer of $f_{\beta_c,0}$. Since $f_{\beta_c,0}''(0) =0$, we conclude that $(\beta_c,0)$ is a special point. This completes the proof of Lemma \ref{crit_part}.
\end{proof}

\begin{lem}
    If $q\ne 2$ or $p\ge 5$, then the special point $(\wbt,\wht)$ satisfies $\wht >0$, and there exists a strictly decreasing, smooth function $\phi_{p,q}: [0,\tilde{h}_{p,q})\rightarrow [0,\infty)$ such that
    \begin{equation*}
        \cC^1_{p,q}= \{(\phi_{p,q}(h),h): h\in [0,\wht)\}.
    \end{equation*}
    Further, $\phi_{p,q}(0) = \beta_c(p,q)$ and $\lim_{x\rightarrow \wht}\phi_{p,q}(x) = \wbt$. 
    Otherwise, i.e. if $q=2$ and $p\in \{2,3,4\}$, then $\cC^1_{p,q}$ is empty.
\end{lem}
\begin{proof}
If $q\neq 2$ or $p\geq 5$ then by Lemma \ref{crit_part} ii. (c), we get that $(\beta_c,0)$ is a strong critical point. This shows that $\cC^1_{p,q}$ is non-empty. Assume there exists $\bm s=(s_1,s_2)^\top$ and $\bm t=(t_1,t_2)^\top$ such that $G(\bm s)=G(\bm t)=(\beta,h)^\top$, such that $(\beta,h) \in \cC^1_{p,q}$. Hence, $f_{\beta, h}^{\prime}\left(s_{1}\right)=f_{\beta, h}^{\prime}\left(s_{2}\right)=f_{\beta, h}^{\prime}\left(t_{1}\right)=f_{\beta, h}^{\prime}\left(t_{2}\right)$. $\f^\prime$ has at most three roots. Also, $s_1\neq s_2$ and $t_1\neq t_2$. If $\bm s \neq \bm t$, then $\f^\prime$ has three distinct roots. Let $s_3$ be a root of $\f^\prime$ such that $\f^\dprime(s_3)>0$. If $s_1=s_3$ then $s_2>s_3$. Now, $\f^\prime>0$ in $(s_3,s_2)$ hence $\f^\prime(s_2)>\f^\prime(s_3)=\f^\prime(s_1)$. So, $s_1\neq s_3$ and similarly $s_2 \neq s_3$. So, $s_1$ and $s_2$ are the maximizers of $\f$. Also, $t_1$ and $t_2$ are the maximizers of $\f$. $\f$ can have at most two maximizers from Lemma \ref{crit_part}. This shows that $s_1=t_1$ and $s_2=t_2$. This is a contradiction.
    
    Therefore, for any $(\beta, h) \in \cC^1_{p,q}$ there exists unique $s_{1}, s_{2}$ such that,
    $$f_{\beta, h}\left(s_{1}\right)=f_{\beta, h}\left(s_{2}\right) \text{ and } f_{\beta, h}^{\prime}\left(s_{1}\right)=f_{\beta, h}^{\prime}\left(s_{2}\right)=0$$
    So, $G^{-1}$ exists on $C^1_{p,q}$. Now, $\cF^{-1}(0)$ is a compact which shows that $G^{-1}: \cC_{p, q}^{1} \rightarrow \cF^{-1}(0)$ is continuous. Also, $G$ is smooth and $J_G$ is invertible in $\cC^1_{p,q}$. Hence, by inverse function theorem, $G^{-1}$ is smooth. It is easy to check that $G\left(\mathcal{\cF}^{-1}(0)\right)=\cC_{p, q}^1$. Also, $G(U)$ is a connected set. Hence, $\Pi_1(G(U))$ and $\Pi_2(G(U))$ are connected sets too where $\Pi_1$ and $\Pi_2$ are projection onto $x$ and $y$-axes. Moreover, $\Pi_{2}\left(\cC_{p,q}^1\right)$ is compact in  $\Pi_{2}(G(U))$ and hence closed in $\Pi_{2}(G(U))$. Note that $\Pi_{1}(G(U))$ and $\Pi_{2}(G(U))$ are non-degenerate intervals. For any  $\beta \in \Pi_1(C_{p, q}^1)$ there exists unique $h$ such that $(\beta, h) \in C_{p, q}^1$. Hence, there exists $\phi_{p, q}: \Pi_{1}\left(C_{p, q}^2 \right) \rightarrow \Pi_{2}\left(C_{p, q}^{1}\right)$
    such that $(\phi_{p,q}(h), h) \in C_{p, q}^{1}$.
    Let $h \in \Pi_{2}\left(\cC_{p,q}^1\right)$. By implicit function theorem there is a neighborhood $U_{h}$ where there exists a smooth function $g$ with,
    \begin{equation*}
        \begin{aligned}
             & F\left(G^{-1}(g(h), h))=0 \quad \forall  h \in U_{h}\right. \\
             \implies & U_{h} \subseteq \Pi_{2}\left(\cC_{p,q}^1\right)
        \end{aligned}
    \end{equation*}
    Hence, $\Pi_{2}\left(\cC_{p,q}^1\right)=\bigcup_{h} U_{h}$ which shows that $\Pi_{2}\left(\cC_{p,q}^1\right)$ is open. Since $\Pi_{2}(G(U))$ is connected and $\Pi_{2}\left(\cC_{p,q}^1\right)$ is clopen and non-empty so $\Pi_{2}\left(C^1_{p,q}\right)=\Pi_{2}(G(U))$. So, $ \Pi_{2}\left(\cC_{p,q}^{1}\right)$ is connected and hence an interval. Similarly, $ \Pi_{1}\left(\cC_{p,q}^{1}\right)$ is interval. Here $\phi_{p, q}$ is monotonic in an interval with image as interval. Hence, $\phi_{p, q}$ is continuous. Implicit function theorem also shows that $\phi_{p, q}$ is smooth.

    Define $\psi_1(h)=\Pi_1(G^{-1}(\phi_{p,q}(h),h))$ and $\psi_2(h)=\Pi_2(G^{-1}(\phi_{p,q}(h),h))$. It is clear that $\psi_{1}(h)=s_{1}$ and $\psi_{2}(h)=s_{2}$, where $s_1$ and $s_2$ are the maximizers of $f_{\phi_{p,q}(h),h}$. Let $\check{\beta}_{p q}=$ $\sup \Pi_{1}\left(C^1_{p,q}\right)$ and $\check{h}_{p, q}=$ $\sup \Pi_{2}\left(C^1_{p,q}\right)$. Now, $\Pi_{1}\left(C^1_{p,q}\right)=\Pi_1(G(U))$ is a non-degenerate interval and hence,  $\left(\check{\beta}_{p, q}, \check{h}_{p, q}\right) \neq\left(\beta_{c}, 0\right)$. Therefore, the curve $\left(\psi_{1}(h), \psi_{2}(h)\right)$ therefore has two end points. By Lemma \ref{crit_part} [ii.(c)], $h=0$ is a strongly critical point and so a boundary point of the interval $\Pi_1(\cC^1_{p,q})$. Hence, $\left(\psi_1(0), \psi_2(0)\right)$ is one of the end points. Let $(x, y)\neq \left(\psi_1(0), \psi_2(0)\right)$ be another end point. If $x=0$ then $h=0$ which further implies that $y=\psi_{2}\left(0\right)$. Hence, a contradiction.
 On the other hand, if $(x, y) \in U^{\circ}$, then by implicit function theorem, there is a neighborhood $N(x, y) \subset U$ such that there exists $d \in C^1(N(x, y))$ with
    \begin{equation*}
        F\left(d(s_{1}), s_{2}\right)=0 . \quad \forall s_{1}, s_{2} \in N(x, y).
    \end{equation*}
    Hence, $N(x,y) \subset U$. This is again a contradiction as $(x,y)$ is a boundary point. Also, $x \neq 1 $ and $ y \neq 1$. Hence, $(x, y) \in\left\{\left(s_{1}, s_{2}\right) \mid s_{1}=s_{2}\right\}$. Therefore, $(x, y) \notin U$. So, $\left(\check{\beta}_{p,q}, \check{h}_{p, q}\right) \notin C^1_{p,q}$.
    Let $\left(\beta_{n}, h_{n}\right) \rightarrow\left(\check{\beta}_{p, q}, \check{h}_{p, q}\right)$ be a sequence such that $\left(\beta_{n}, h_{n}\right) \in C^1_{p,q}$. Hence, $\left(s_{1, n}, s_{2, n}\right) \rightarrow(x, y)=(s, s)$ as $G^{-1}$ is continuous where $x=y=s$ is the other end point.
    Now,
    \begin{equation*}
        \begin{aligned}
             & f_{\beta_{n}, h_{n}}'\left(s_{1,n}\right)-f_{\beta_{n}, h_{n}}'\left(s_{2,n}\right)=0                                                      \\
             & \Rightarrow \frac{f_{\beta_{n}, h_n}'\left(s_{1,n}\right)-f_{\beta_{n}, h_{n}}'\left(s_{2,n}\right)}{s_{1,n}-s_{2,n}}=0 \\
             & \Rightarrow f_{\check{\beta}_{p.q},\check{h}_{p, q}}^\dprime(s)=0
        \end{aligned}
    \end{equation*}
   Therefore $\check{\beta}_{p, q}, \check{h}_{p, q}$ is a special point, where $\check{h}_{p, q}>0$. Hence, by Lemma \ref{specialtype2thm}, $(\check{\beta}_{p, q}, \check{h}_{p, q})=(\wbt,\wht)$.
    
    If $q=2$ and $p=2,3,4$ then by Lemma \ref{crit_part} ii. (c), we get that $(\beta_c,0)$ is a special point. This further suggest that $(\wbt,\wht)=(\beta_c,0)$ as special points are unique from Lemma \ref{specialtype2thm}. Now, by Lemma \ref{crit_less_spec} we get that $\cC^1_{p,q}$ is empty.
\end{proof}
\begin{lem}\label{strange_spec}
    Let $(\beta,h) \in \cS_{p,q}$. Then $\bm x_0=(\frac{1}{q},\ldots,\frac{1}{q})$ is a global maximizer of $\hf$ if and only if $p\in \{2,3,4\}$ and $q=2$.
\end{lem}
\begin{proof}
    The proof of Lemma \ref{strange_spec} follows from Proposition \ref{spcl7} i., Lemma \ref{multf1} iii. and its proof. 
\end{proof}

\section{Technical Lemmas Relevant to Maximum Likelihood Estimation}
In this section, we collect some technical results that are relevant to maximum likelihood estimation of $\beta$ and $h$. 
\begin{lem}
\label{ml_exp}
    $\hat{\beta}_N$ is a solution of the equation (in $\beta$),
\begin{equation*}
    \mathbb{E}_{\beta, h, p}\left(\|\epm\|_p^p\right)=\|\epm\|_p^p~,
\end{equation*}
and for fixed $\beta\in \mathbb{R}$, $\hat{h}_N$ is a solution of the equation (in $h$),
\begin{equation*}
   \mathbb{E}_{\beta, h, p}\left(\bar{X}_{\cdot 1}\right)=\bar{X}_{\cdot 1}.
\end{equation*}
\end{lem}
\begin{proof}
    The log-likelihood function is given by,
\begin{equation*}
    \ell_N(\beta,h,\bm x)=N\beta \nrm{\bar{\bm x}_N}_p^p + Nh\bm \bar{x}_{\cdot 1}-F_N(\beta,h),
\end{equation*}
where $F_N :=\log(q^N Z_N(\beta,h))$. Hence, $$\frac{\partial}{\partial \beta}\ell_N(\beta,h,\bm x)=N \nrm{\bar{\bm x}_N}_p^p-\frac{\partial}{\partial \beta}F_N(\beta,h) = N\|\bar{\bm x}_N\|_p^p - N\e_{\beta,h,p}(\|\epm\|_p^p)$$ 
and
$$\frac{\partial}{\partial h}\ell_N(\beta,h,\bm x)=Nh\bm \bar{x}_{\cdot 1}-\frac{\partial}{\partial h}F_N(\beta,h) = N\bar{x}_{\cdot 1} - N\e_{\beta,h,p}(\bar{X}_{\cdot 1})~.$$
The proof of Lemma \ref{ml_exp} is now complete.
\end{proof}

\begin{lem}
    For every fixed $h$, the function $\beta \mapsto F_N(\beta, h, p)$ is strictly convex, and for every fixed $\beta$, the function $h \mapsto F_N(\beta, h, p)$ is strictly convex. Consequently, the maps $u_{N, 1}$ and $u_{N, p}$ are strictly increasing in both $\beta$ and $h$.
    \label{like_monotone}
\end{lem}

\begin{proof}
    Let $\psi_N(\beta, h):=F_N(\beta, h, p) - N \log q=\log \sum_{\epm \in \mathcal{C}_N} e^{N \beta \nrm{\epm}_p^p +N h \bar{X}_{\cdot 1}}$. Then for every $\beta_1, \beta_2, h$ and $\lambda \in(0,1)$, we have by Hölder's inequality,
    $$
        \begin{aligned}
            \psi_N\left(\lambda \beta_1+(1-\lambda) \beta_2, h\right) & =\log \sum_{\epm \in \mathcal{C}_N} e^{N \lambda\left(\beta_1 \nrm{\epm}_p^p +h \bar{X}_{\cdot 1}\right)} e^{N(1-\lambda)\left(\beta_2 \nrm{\epm}_p^p +h \bxd\right)}                                                  \\
                                                                      & <\log \left[\left(\sum_{\epm \in \mathcal{C}_N} e^{N \beta_1 \nrm{\epm}_p^p +N h \bxd}\right)^\lambda\left(\sum_{\epm \in \mathcal{C}_N} e^{N \beta_2 \nrm{\epm}_p^p +N h \bxd}\right)^{1-\lambda}\right] \\
                                                                      & =\lambda \psi_N\left(\beta_1, h\right)+(1-\lambda) \psi_N\left(\beta_2, h\right)
        \end{aligned}
    $$
    Similarly, for every $h_1, h_2, \beta$ and $\lambda \in(0,1)$, we have by Hölder's inequality,
    $$
        \begin{aligned}
            \psi_N\left(\beta, \lambda h_1+(1-\lambda) h_2\right) & =\log \sum_{\epm \in \mathcal{C}_N} e^{N \lambda\left(\beta \nrm{\epm}_p^p +h_1 \bxd\right)} e^{N(1-\lambda)\left(\beta \nrm{\epm}_p^p +h_2 \bxd\right)}                                                  \\
                                                                  & <\log \left[\left(\sum_{\epm \in \mathcal{C}_N} e^{N \beta \nrm{\epm}_p^p +N h_1 \bxd}\right)^\lambda\left(\sum_{\epm \in \mathcal{C}_N} e^{N \beta \nrm{\epm}_p^p +N h_2 \bxd}\right)^{1-\lambda}\right] \\
                                                                  & =\lambda \psi_N\left(\beta, h_1\right)+(1-\lambda) \psi_N\left(\beta, h_2\right) .
        \end{aligned}
    $$
    This shows strict convexity of the functions $\beta \mapsto F_N(\beta, h, p)$ and $h \mapsto F_N(\beta, h, p)$. Now, note that
    $$
        \frac{\partial}{\partial \beta} F_N(\beta, h, p)=N u_{N, p}(\beta, h, p) \quad \text { and } \quad \frac{\partial}{\partial h} F_N(\beta, h, p)=N u_{N, 1}(\beta, h, p) .
    $$
    Lemma \ref{like_monotone} now follows from the fact that the first derivative of a differentiable, strictly convex function is strictly increasing.
\end{proof}
\end{document}